\documentclass{amsart}
\usepackage{palatino,amssymb,epsfig,amsmath,latexsym,amsthm}

\input rlepsf.sty

\theoremstyle{plain}
\newtheorem{theorem}{Theorem}[section]
\newtheorem{lemma}[theorem]{Lemma}
\newtheorem*{claim}{Claim}
\newtheorem{proposition}[theorem]{Proposition}
\newtheorem{corollary}[theorem]{Corollary}
\newtheorem{question}[theorem]{Question}

\newtheorem{problem}[theorem]{Problem}
\newtheorem{Bounded Diameter Lemma}[theorem]{Bounded Diameter Lemma}
\theoremstyle{definition}
\newtheorem{definition}[theorem]{Definition}
\newtheorem*{thin_problem}{Thin Obstacle Problem}
\newtheorem{remark}[theorem]{Remark}
\newtheorem{acknowledgements}[theorem]{Acknowledgements}

\newtheorem{remarks}[theorem]{Remarks}

\newtheorem{notation}[theorem]{Notation}

\DeclareMathOperator{\axis}{axis}

\DeclareMathOperator{\length}{length}
\DeclareMathOperator{\area}{area}

\DeclareMathOperator{\rank}{rank}
\DeclareMathOperator{\genus}{genus}
\DeclareMathOperator{\kernal}{ker}
\DeclareMathOperator{\rel}{rel}

\DeclareMathOperator{\textint}{int}
\DeclareMathOperator{\Bag}{Bag}

\newcommand{\Sinfty}{S^2_{\infty}}

\newcommand{\finv}{f^{-1}}

\newcommand{\piinv}{\pi^{-1}}

\newcommand{\BH}{\mathbb H}
\newcommand{\BN}{\mathbb N}
\newcommand{\BR}{\mathbb R}
\newcommand{\BZ}{\mathbb Z}
\newcommand{\BZt}{{\mathbb Z}/2{\mathbb Z}}

\newcommand{\mE}{\mathcal{E}}
\newcommand{\mF}{\mathcal{F}}
\newcommand{\mG}{\mathcal{G}}
\newcommand{\mW}{\mathcal{W}}
\newcommand{\Nthick}{N_{\ge\epsilon}}
\newcommand{\Nthin}{N_{\le\epsilon}}
\newcommand{\pC}{\partial_\mE C}

\def\CAT{{\text{CAT}}}
\def\inte{{\text{int}}}
\def\dist{{\text{dist}}}
\def\diam{{\text{diam}}}
\def\tlength{{\text{length}_t}}
\def\tarea{{\text{area}_t}}
\def\new{\text{new}}
\def\E{\mathcal{E}}
\def\C{\mathcal{C}}
\def\H{\mathbb{H}}
\def\R{\mathbb{R}}
\def\Z{\mathbb{Z}}
\def\K{\mathbb{K}}

\begin{document}

\title{Shrinkwrapping and the Taming of Hyperbolic 3--Manifolds}
\author{Danny Calegari}
\address{Department of Mathematics \\Caltech\\ Pasadena, CA 91125}
\thanks{partially supported by Therese Calegari and NSF grant DMS-0405491.}
\author{David Gabai}\address{Department of Mathematics\\Princeton
University\\Princeton, NJ 08544}
\thanks{partially supported by NSF grant DMS-0071852}
\thanks{Version 1.00, October 22, 2005}

\maketitle

\setcounter{section}{-1}

\section{Introduction}\label{S0}

During the period 1960--1980,
Ahlfors, Bers, Kra, Marden, Maskit,
Sullivan, Thurston and many others developed the theory of
{\em geometrically finite ends} of hyperbolic $3$--manifolds.
It remained to understand those ends which are not geometrically finite; such
ends are called {\em geometrically infinite}.

Around 1978 William Thurston gave a conjectural description of geometrically
infinite ends of complete hyperbolic $3$--manifolds.
An example of a geometrically infinite end is given by an infinite
cyclic covering space of a closed
hyperbolic 3-manifold which fibers over the circle.
Such an end has cross sections of uniformly bounded area. By contrast, the
area of sections of geometrically finite ends
grow exponentially in the distance from the convex core.

For the sake of clarity we will assume throughout this introduction that
$N=\BH^3/\Gamma$ where $\Gamma$ is parabolic free.  Precise
statements of the parabolic case will be
given in \S 7.

Thurston's idea was formalized by Bonahon \cite{Bo} and Canary
\cite{Ca} with the following.

\begin{definition}  An end $\mE$ of a hyperbolic 3-manifold $N$ is
\emph{simply degenerate} if it has a closed neighborhood 
of the form $S\times [0,\infty)$ where $S$ is a closed
surface, and there exists a
sequence $\{S_i\}$ of $\CAT(-1)$ surfaces exiting $\mE$ which are
homotopic to $S\times 0$ in
$\mE$.  This means that there exists a sequence of maps $f_i:S\to N$
such that the induced path metrics
induce $\CAT(-1)$ structures on the $S_i$'s, $f(S_i)\subset
S\times [i,\infty)$ and $f$ is homotopic to a homeomorphism onto
$S\times 0$ via a homotopy supported
in $S\times [0,\infty)$.
\end{definition}

Here by $\CAT(-1)$, we mean as usual a geodesic metric space for which geodesic
triangles are ``thinner'' than comparison triangles in hyperbolic space. If the
metrics pulled back by the $f_i$ are smooth,
this is equivalent to the condition that the Riemannian
curvature is bounded above by $-1$. See \cite{BriHae} for
a reference. Note that by Gauss--Bonnet, the area of a
$\CAT(-1)$ surface can be estimated from its Euler characteristic; it follows
that a simply degenerate end has cross sections of uniformly bounded area,
just like the end of a cyclic cover of a manifold fibering over the circle.

Francis Bonahon \cite{Bo} observed that geometrically infinite ends are
exactly those ends possessing an exiting sequence of closed geodesics.
This will be our working definition of such ends throughout
this paper.

The following is our main  result.

\begin{theorem}\label{clean}  An end $\mE$ of a complete  hyperbolic
3-manifold $N$ with finitely generated fundamental group is simply
degenerate if there exists a sequence of closed geodesics exiting
$\mE$.\end{theorem}

Consequently we have,

\begin{theorem}\label{geometrically tame}  Let $N$ be a complete hyperbolic
3-manifold with finitely generated fundamental group. Then every end of $N$ is
geometrically tame, i.e. it is either geometrically finite or simply
degenerate. \end{theorem}

In 1974 Marden \cite{Ma} showed that a geometrically finite hyperbolic
3-manifold is \emph{topologically tame}, i.e. is the interior of a compact
3-manifold.  He asked whether all complete hyperbolic 3-manifolds
with finitely generated fundamental group are
topologically tame.  This question is now
known as the \emph{Tame Ends Conjecture} or \emph{Marden Conjecture}.

\begin{theorem} \label{marden}  If $N$ is a complete hyperbolic 3-manifold with
finitely generated fundamental group, then $N$ is topologically tame.
\end{theorem}

Ian Agol \cite{Ag} has independently proven Theorem \ref{marden}.

There have been many important steps towards Theorem \ref{clean}.
The seminal result
was obtained by Thurston
(\cite{T}, Theorem 9.2)  who proved Theorems \ref{geometrically tame}
and \ref{marden} for certain algebraic
limits of quasi Fuchsian groups.  Bonahon
\cite{Bo} established Theorems \ref{clean} and \ref{marden} when
$\pi_1(N)$ is freely
indecomposible
and Canary \cite{Ca} proved that topological tameness  implies
geometrical tameness.  Results in
the direction of \ref{marden} were also obtained by Canary-Minsky
\cite{CM},  Kleineidam--Souto \cite{KS}, Evans
\cite{Ev}, Brock--Bromberg--Evans--Souto \cite{BBES}, Ohshika,
Brock--Souto \cite{BS} and Souto \cite{So}.  

Thurston first discovered how to obtain analytic conclusions from
the existence of exiting sequences of
$\CAT(-1)$ surfaces.  Thurston's work as generalized by
Bonahon \cite{Bo} and Canary \cite{Ca}
combined with Theorem \ref{clean} yields a positive proof of the
Ahlfors' Measure Conjecture  \cite{A2}.

\begin{theorem}\label{ahlfors}  If $\Gamma$ is a finitely generated
Kleinian group, then the limit set $L_\Gamma$ is either $\Sinfty$ or has Lebesgue
measure  zero.  If $L_\Gamma=\Sinfty$, then
$\Gamma$ acts ergodically on $\Sinfty$.
\end{theorem}

Theorem \ref{ahlfors} is one of the many analytical consequences of
our main result.
Indeed Theorem \ref{clean} implies that a complete hyperbolic
3-manifold $N$ with finitely generated
fundamental group is \emph{analytically tame} as defined by Canary
\cite{Ca}.  It follows from Canary that the various results of \S 9
\cite{Ca} hold  for $N$.

Our main result is the last step needed to  prove
the following monumental result, the other parts being established by
Alhfors, Bers, Kra, Marden, Maskit, Mostow, Prasad, Sullivan,
Thurston, Minsky, Masur--Minsky, Brock--Canary--Minsky, Ohshika,
Kleineidam--Souto, Lecuire, Kim--Lecuire--Ohshika, Hossein--Souto and Rees. See
\cite{Mi} and \cite{BCM}.

\begin{theorem}[Classification Theorem] If $N$ is a complete
hyperbolic 3-manifold with finitely
generated fundamental group, then $N$ is determined up to isometry by
its topological type, the
conformal boundary of its geometrically finite ends and the ending
laminations of its geometrically
infinite ends. \end{theorem}

The following result was conjectured by Bers, Sullivan and Thurston.
Theorem \ref{marden} is one of many results, many of them recent,
needed to build a proof.  Major contributions were made by Alhfors, Bers, Kra,
Marden, Maskit, Mostow, Prasad, Sullivan, Thurston, Minsky, Masur--Minsky,
Brock--Canary--Minsky, Ohshika, Kleineidam--Souto, Lecuire,
Kim--Lecuire--Ohshika, Hossein--Souto, Rees, Bromberg and Brock--Bromberg.

\begin{theorem}[Density Theorem]\label{density}  If $N=\BH^3/\Gamma$ is a
complete finitely generated 3-manifold with finitely generated fundamental
group, then
$\Gamma$ is the algebraic limit of geometrically finite Kleinian groups.\end
{theorem}

The main technical
innovation of this paper is a new technique called \emph{shrinkwrapping} for
producing
$\CAT(-1)$ surfaces in hyperbolic 3-manifolds.  Historically, such
surfaces have been immensely important in the study of hyperbolic
3-manifolds, e.g. see \cite{T}, \cite{Bo}, \cite{Ca} and \cite{CM}.

        Given a locally finite
set $\Delta$ of pairwise disjoint simple closed curves in
the 3-manifold $N$, we say that the embedded surface $S\subset N$ is
\emph{2-incompressible} rel. $\Delta$ if every compressing disc for
$S$ meets $\Delta$ at
least twice. Here is a sample theorem.

\begin{theorem}[Existence of shrinkwrapped surface]\label{sample theorem}
Let $M$ be a complete, orientable, parabolic free hyperbolic
$3$--manifold, and let
$\Gamma$ be a finite collection of pairwise disjoint simple
closed geodesics in $M$. Further, let
$S \subset M \backslash \Gamma$ be a closed embedded $2$--incompressible
surface rel. $\Gamma$ which is either nonseparating in $M$ or separates some
component of $\Gamma$ from another. Then $S$ is homotopic to a
$\CAT(-1)$ surface $T$ via a homotopy
$$F:S \times [0,1] \to M$$
such that
\begin{enumerate}
\item{$F(S \times 0) = S$}
\item{$F(S \times t) = S_t$ is an embedding disjoint from $\Gamma$ for
$0 \le t < 1$}
\item{$F(S \times 1) = T$}
\item{If $T'$ is any other surface with these properties,
then $\text{area}(T) \le \text{area}(T')$}
\end{enumerate}
We say that $T$ is obtained from $S$ by {\em shrinkwrapping rel. $\Gamma$},
or if $\Gamma$ is understood, $T$ is obtained from $S$ by
{\em shrinkwrapping}.
\end{theorem}

In fact, we prove the stronger result that $T$ is {\em $\Gamma$--minimal}
(to be defined in \S 1) which implies in particular that it
is intrinsically $\CAT(-1)$

Here is the main technical result of this paper.\vskip 12 pt

\begin{theorem}\label{main} Let $\mE$ be an end
of the complete orientable hyperbolic 3-manifold $N$ with
finitely generated fundamental group.  Let $C$ be a 3-dimensional
compact core of $N$, $\partial_{\mE}C$ the component of
$\partial C$ facing $\mE$ and $g=genus(\partial_{\mE}C)$.  If there
exists a sequence of closed
geodesics exiting $\mE$, then there exists a sequence $\{S_i\}$ of
$\CAT(-1)$ surfaces of genus $g$ exiting
$\mE$ such that each $\{S_i\}$ is homologically separating in $\mE$.
That is, each $S_i$ homologically separates $\partial_{\mE}C$ from
$\mE$.
\end{theorem}

Theorem \ref{marden} can now be deduced from Theorem \ref{main} and Souto 
\cite{So}; however, we prove that Theorem \ref{main} implies Theorem
\ref{marden} using only 3-manifold topology and elementary hyperbolic
geometry.  

\vskip 12pt

The proof of Theorem \ref{main} blends elementary aspects
of minimal surface theory, hyperbolic geometry, and 3-manifold topology.
The method will be demonstrated in \S4
where we give a proof of Canary's theorem.  The first time reader is urged to
begin with that section.

This paper is organized as follows.  In \S1 and \S2 we establish
the shrinkwrapping technique for finding
$\CAT(-1)$ surfaces in hyperbolic 3-manifolds.    In \S3
we prove the existence of $\epsilon$-separated simple geodesics
exiting the end of parabolic free manifolds.  In \S4 we  prove
Canary's theorem.  This proof will  model  the proof of the general case.
The general strategy will be outlined at the end of that section.  In
\S5 we develop the topological theory of end reductions in
3-manifolds.  In \S6 we give the proofs of our main results.
In \S7 we give the necessary embellishments of our
methods to state and prove our results in the case of manifolds with
parabolic cusps.

\begin{notation} If $X\subset Y$, then $N(X)$ denotes a regular
neighborhood of $X$ in $Y$ and
$\inte(X)$ denotes the interior of $X$.  If
$X$ is a topological space, then $|X|$ denotes the
number of components of $X$. If $A,B$ are topological subspaces of a third
space, then $A\backslash B$ denotes the intersection of $A$ with the
complement of $B$.
\end{notation}

\begin{acknowledgements}  The first author is grateful to Nick Makarov
for some useful analytic discussions.
The second author is grateful to Michael Freedman
for many long conversations in Fall 1996 which introduced him to the Tame Ends
conjecture.  He thanks Francis Bonahon, Yair
Minsky and Jeff Brock for their interest and helpful comments.  Part
of this research was carried out while he was visiting Nara Women's
University, the Technion and the Institute for Advanced Study.  He
thanks them for their hospitality.  We thank the referees for their many thoughtful
suggestions and comments.  
\end{acknowledgements}
\section{Shrinkwrapping}\label{shrink}

In this section, we introduce a new technical tool for finding $\CAT(-1)$
surfaces in hyperbolic $3$--manifolds, called {\em shrinkwrapping}. Roughly
speaking, given a collection of simple closed geodesics $\Gamma$
in a hyperbolic $3$--manifold $M$ and an embedded surface
$S \subset M \backslash \Gamma$, a surface $T \subset M$ is
obtained from $S$ by {\em shrinkwrapping $S$ rel. $\Gamma$} if it
homotopic to $S$, can be approximated by an isotopy from $S$ supported in
$M \backslash \Gamma$, and is least area subject to these constraints.

Given mild topological conditions on $M,\Gamma,S$
(namely $2$--incompressibility,
to be defined below) the shrinkwrapped surface exists, and is $\CAT(-1)$ with
respect to the path metric induced by the Riemannian metric on $M$.

We use some basic analytical tools throughout this section, including the
Gauss--Bonnet formula, the coarea formula, and the Arzela--Ascoli theorem.
At a number of points we must invoke results from the literature to establish
existence of minimal surfaces (\cite{MSY}), existence of limits with area and
curvature control (\cite{CiSc}), and regularity of the shrinkwrapped surfaces
along $\Gamma$ (\cite{Rich},\cite{Frehse}). General references are \cite{CoMi}, \cite{Jost}
\cite{Fed} and \cite{Ball}.

\subsection{Geometry of surfaces}

For convenience, we state some elementary but
fundamental lemmas concerning curvature of (smooth) surfaces in Riemannian
$3$-manifolds.  

We use the following standard terms to refer to different kinds of minimal surfaces:

\begin{definition}
A smooth surface $\Sigma$ in a Riemannian $3$-manifold is {\em minimal} if
it is a critical point for area with respect to all smooth compactly
supported variations. It is {\em locally least area} (also called {\em stable})
if it is a local minimum for area with respect to all smooth, compactly supported variations.
A closed, embedded surface is {\em globally least area} if it is an absolute
minimum for area amongst all smooth surfaces in its isotopy class.
\end{definition} 

Note that we do not require that our minimal or locally least area surfaces are complete.

Any subsurface of a globally least area surface is locally least area, and a locally least area
surface is minimal. A smooth surface is minimal iff its mean curvature vector field vanishes
identically. For more details, consult \cite{CoMi}, especially chapter 5.

The intrinsic curvature of a minimal surface is controlled by the geometry of the
ambient manifold. The following lemma is formula 5.6 on page 100 of \cite{CoMi}.

\begin{lemma}[Monotonicity of curvature]\label{curvature_comparison}
Let $\Sigma$ be a minimal surface in a Riemannian manifold $M$. Let $K_\Sigma$
denote the curvature of $\Sigma$, and $K_M$ the sectional curvature of $M$.
Then restricted to the tangent space $T\Sigma$,
$$K_\Sigma = K_M - \frac 1 2 |A|^2$$
where $A$ denotes the second fundamental form of $\Sigma$.
\end{lemma}

In particular, if the Riemannian curvature on $M$ is bounded from above by some
constant $K$, then the curvature of a minimal surface $\Sigma$ in $M$ is
also bounded above by $K$.

The following lemma is just the usual Gauss--Bonnet formula:

\begin{lemma}[Gauss--Bonnet formula]\label{Gauss_Bonnet}
Let $\Sigma$ be a $C^3$ Riemannian surface with (possibly empty) $C^3$ boundary
$\partial \Sigma$. Let $K_\Sigma$ denote the Gauss curvature of $\Sigma$, and
$\kappa$ the geodesic curvature along $\partial \Sigma$. Then
$$\int_\Sigma K_\Sigma = 2\pi\chi(\Sigma) - \int_{\partial \Sigma} 
\kappa \, dl$$
\end{lemma}

Many simple proofs exist in the literature. For example, see \cite{Jost}.

If $\partial \Sigma$ is merely piecewise $C^3$, with finitely many corners
$p_i$ and external angles $\alpha_i$, the Gauss--Bonnet formula must be
modified as follows:

\begin{lemma}[Gauss--Bonnet with corners]\label{Gauss_Bonnet_corners}
Let $\Sigma$ be a $C^3$ Riemannian surface with boundary $\partial \Sigma$
which is piecewise $C^3$, and has external angles $\alpha_i$ at finitely many
points $p_i$. Let $K_\Sigma$ and $\kappa$ be as above. Then
$$\int_\Sigma K_\Sigma = 2\pi\chi(\Sigma) - \int_{\partial \Sigma} \kappa \, dl
- \sum_i \alpha_i$$
\end{lemma}

Observe for $abc$ a geodesic triangle with external angles
$\alpha_1,\alpha_2,\alpha_3$
that Lemma~\ref{Gauss_Bonnet_corners} implies
$$\int_{abc} K = 2\pi - \sum_i \alpha_i$$

Notice that the geodesic curvature $\kappa$ vanishes precisely when 
$\partial \Sigma$
is a geodesic, that is, a critical point for the length functional. 
More generally,
let $\nu$ be the normal bundle of $\partial \Sigma$ in $\Sigma$, 
oriented so that
the inward unit normal is a positive section. The exponential
map restricted to $\nu$ defines a map
$$\phi:\partial \Sigma \times [0,\epsilon] \to \Sigma$$
for small $\epsilon$, where $\phi(\cdot,0) = \text{Id}|_{\partial \Sigma}$, and
$\phi(\partial \Sigma, t)$ for small $t$ is the boundary in $\Sigma$ of the
tubular $t$ neighborhood of $\partial \Sigma$. Then
$$\int_{\partial \Sigma} \kappa \, dl = - \frac d {dt}\Bigl|_{t=0}
\length(\phi_t(\partial \Sigma))$$

Note that if $\Sigma$ is a surface with sectional curvature bounded above by
$-1$, then by integrating this formula we see that the ball $B_t(p)$
of radius $t$ in $\Sigma$ about a point $p \in \Sigma$ satisfies
$$\area(B_t(p)) \ge 2\pi(\cosh(t) - 1) > \pi t^2$$
for small $t>0$.
\subsection{Comparison geometry}

For basic elements of the theory of comparison geometry, see \cite{BriHae}.

\begin{definition}[Comparison triangle]
Let $a_1a_2a_3$ be a geodesic triangle in a geodesic metric space 
$X$. Let $\kappa \in \R$
be given. A {\em $\kappa$--comparison triangle} is a geodesic triangle
$\overline{a_1}\overline{a_2}\overline{a_3}$ in the complete simply--connected
Riemannian $2$--manifold of
constant sectional curvature $\kappa$, where the edges $a_ia_j$ and
$\overline{a_i}\overline{a_j}$ satisfy
$$\length(a_ia_j) = \length(\overline{a_i}\overline{a_j})$$
Given a point $x \in a_1a_2$ on one of the edges of $a_1a_2a_3$, 
there is a corresponding
point $\overline{x} \in \overline{a_1}\overline{a_2}$ on one of the 
edges of the
comparison triangle, satisfying $$\length(a_1x) = 
\length(\overline{a_1}\overline{x})$$ and
$$\length(xa_2) = \length(\overline{x}\overline{a_2})$$
\end{definition}

\begin{remark}
Note that if $\kappa > 0$, the comparison triangle might not exist if 
the edge lengths
are too big, but if $\kappa \le 0$ the comparison triangle always 
exists, and is
unique up to isometry.
\end{remark}

There is a slight issue of terminology to be aware of here. In a {\em surface},
a triangle is a polygonal disk with $3$ geodesic edges. In a {\em 
path metric space},
a triangle is just a union of $3$ geodesic segments with common endpoints.

\begin{definition}[CAT($\kappa$)]
Let $S$ be a closed surface with a path metric $g$. Let $\tilde{S}$ denote
the universal cover of $\tilde{S}$, with path metric induced by the pullback
of the path metric $g$. Let $\kappa \in \R$ be given.
$S$ is said to be {\em $\CAT(\kappa)$} if for every geodesic triangle $abc$ in
$\tilde{S}$, and every point $z$ on the edge $bc$, the distance in $\tilde{S}$
from $a$ to $z$ is no more than the distance from $\overline{a}$ to 
$\overline{z}$ in
a $\kappa$--comparison triangle.
\end{definition}

By Lemma~\ref{Gauss_Bonnet_corners} applied to geodesic triangles, one can
show that a $C^3$ surface $\Sigma$ with sectional curvature $K_\Sigma$ satisfying
$K_\Sigma \le \kappa$ everywhere is $\CAT(\kappa)$ with respect to the
Riemannian path metric. This fact is essentially due to Alexandrov; see \cite{Ball}
for a proof.

More generally, suppose $\Sigma$ is a surface which is $C^3$ outside a
closed, nowhere dense subset $X \subset \Sigma$. Furthermore, suppose that 
$K_\Sigma \le \kappa$ holds in $\Sigma \backslash X$, and suppose that 
the formula from Lemma~\ref{Gauss_Bonnet_corners} holds for every geodesic triangle with
vertices in $\Sigma \backslash X$ (which is a dense set of geodesic triangles).
Then the same argument shows that $\Sigma$ is $\CAT(\kappa)$. See e.g.
\cite{Resh} \S 8, pp. 135--140 for more details, and a general discussion of metric surfaces
with (integral) curvature bounds.

\begin{definition}[$\Gamma$--minimal surfaces]
Let $\kappa \in \R$ be given. Let $M$ be a complete Riemannian $3$--manifold
with sectional curvature bounded above by $\kappa$, and let
$\Gamma$ be an embedded collection of simple closed
geodesics in $M$. An immersion
$$\psi:S \to M$$
is {\em $\Gamma$--minimal} if it is smooth with mean curvature $0$ in 
$M \backslash \Gamma$,
and is metrically $\CAT(\kappa)$ with respect to the path metric 
induced by $\psi$
from the Riemannian metric on $M$.
\end{definition}

Notice by Lemma~\ref{curvature_comparison} that a smooth surface $S$ 
with mean curvature
$0$ in $M$ is $\CAT(\kappa)$, so a minimal surface (in the usual sense) is an
example of a $\Gamma$--minimal surface.

\subsection{Statement of shrinkwrapping theorem}

\begin{definition}[$2$--incompressibility]
An embedded surface $S$ in a $3$--manifold $M$ disjoint from a 
collection $\Gamma$
of simple closed curves is said to be {\em $2$--incompressible rel. 
$\Gamma$} if
any essential compressing disk for $S$ must intersect $\Gamma$ in at 
least two points. If
$\Gamma$ is understood, we say $S$ is {\em $2$--incompressible}.
\end{definition}

\begin{theorem}[Existence of shrinkwrapped surface]\label{shrinkwrap_exists}
Let $M$ be a complete, orientable, parabolic free hyperbolic 
$3$--manifold, and let
$\Gamma$ be a finite collection of pairwise disjoint simple closed geodesics in
$M$. Further, let $S \subset M \backslash \Gamma$ be a closed 
embedded $2$--incompressible
surface rel. $\Gamma$ which is either nonseparating in $M$ or separates some
component of $\Gamma$ from another. Then $S$ is homotopic to a 
$\Gamma$--minimal surface
$T$ via a homotopy
$$F:S \times [0,1] \to M$$
such that
\begin{enumerate}
\item{$F(S \times 0) = S$}
\item{$F(S \times t) = S_t$ is an embedding disjoint from $\Gamma$ 
for $0 \le t < 1$}
\item{$F(S \times 1) = T$}
\item{If $T'$ is any other surface with these properties,
then $\text{area}(T) \le \text{area}(T')$}
\end{enumerate}
We say that $T$ is obtained from $S$ by {\em shrinkwrapping rel. 
$\Gamma$}, or if
$\Gamma$ is understood, $T$ is obtained from $S$ by {\em shrinkwrapping}.
\end{theorem}

The remainder of this section will be taken up with the proof of
Theorem~\ref{shrinkwrap_exists}.

\begin{remark}
In fact, for our applications, the property we want to use of our surface $T$
is that we can estimate its diameter (rel. the thin part of $M$) from its
Euler characteristic. This follows from a Gauss--Bonnet estimate and the bounded diameter
lemma (Lemma~\ref{bounded_diameter_lemma}, to be proved below). In fact, our
argument will show directly that the surface $T$ satisfies Gauss--Bonnet; the
fact that it is $\CAT(-1)$ is logically superfluous for the purposes 
of this paper.
\end{remark}

\subsection{Deforming metrics along geodesics}

\begin{definition}[$\delta$--separation]
Let $\Gamma$ be a collection of disjoint simple geodesics in a 
Riemannian manifold
$M$. The collection $\Gamma$ is {\em $\delta$--separated} if any path
$\alpha:I \to M$ with endpoints on $\Gamma$ and satisfying
$$\length(\alpha(I)) \le \delta$$
is homotopic rel. endpoints into $\Gamma$. The supremum of such $\delta$
is called the {\em separation constant} of $\Gamma$.
The collection $\Gamma$ is
{\em weakly $\delta$--separated} if
$$\dist(\gamma,\gamma') > \delta$$
whenever $\gamma,\gamma'$ are distinct components of $\Gamma$.
The supremum of such $\delta$ is called the {\em weak separation 
constant} of $\Gamma$.
\end{definition}

\begin{definition}[Neighborhood and tube neighborhood]
Let $r > 0$ be given. For a point $x \in M$, we let $N_r(x)$ denote the closed
ball of radius $r$ about $x$, and let $N_{<r}(x),\partial N_r(x)$ 
denote respectively
the interior and the boundary of $N_r(x)$. For a closed geodesic 
$\gamma$ in $M$,
we let $N_r(\gamma)$ denote the closed tube of radius $r$ about $\gamma$, and
let $N_{<r}(\gamma),\partial N_r(\gamma)$ denote respectively the interior and
the boundary of $N_r(\gamma)$. If $\Gamma$ denotes a union of 
geodesics $\gamma_i$,
then we use the shorthand notation
$$N_r(\Gamma) = \bigcup_{\gamma_i} N_r(\gamma_i)$$
\end{definition}

\begin{remark}
Topologically, $\partial N_r(x)$ is a sphere and $\partial 
N_r(\gamma)$ is a torus,
for sufficiently small $r$. Similarly, $N_r(x)$ is a closed ball, and 
$N_r(\gamma)$
is a closed solid torus. If $\Gamma$ is $\delta$--separated, then
$N_{\delta/2}(\Gamma)$ is a union of solid tori.
\end{remark}

\begin{lemma}[Bounded Diameter Lemma]\label{bounded_diameter_lemma}
Let $M$ be a complete hyperbolic $3$--manifold. Let $\Gamma$ be a disjoint
collection of $\delta$--separated embedded geodesics. Let $\epsilon > 0$ be
a Margulis constant for dimension $3$, and let $M_{\le \epsilon}$ denote
the subset of $M$ where the injectivity radius is at most $\epsilon$.
If $S \subset M\backslash \Gamma$
is a $2$--incompressible $\Gamma$--minimal surface, then there is a constant
$C = C(\chi(S),\epsilon,\delta) \in \R$ and $n = 
n(\chi(S),\epsilon,\delta) \in \Z$
such that for each component $S_i$ of $S \cap (M \backslash M_{\le \epsilon})$,
we have
$$\diam(S_i) \le C$$
Furthermore, $S$ can only intersect at most $n$ components of $M_{\le \epsilon}$.
\end{lemma}
\begin{proof}
Since $S$ is $2$--incompressible, any point $x \in S$ either lies
in $M_{\le \epsilon}$, or is the center of an embedded $m$--disk in $S$,
where
$$m = \min(\epsilon/2,\delta/2)$$

Since $S$ is $\CAT(-1)$, Gauss--Bonnet implies that the area of an embedded
$m$--disk in $S$ has area at least $2\pi(\cosh(m)-1) > \pi m^2$.

This implies that if $x \in S \cap M \backslash M_{\le\epsilon}$ then
$$\area(S \cap N_m(x)) \ge \pi m^2$$
The proof now follows by a standard covering argument.
\end{proof}

A surface $S$ satisfying the conclusion of the Bounded Diameter Lemma 
is sometimes
said to have diameter bounded by $C$ modulo $M_{\le \epsilon}$.

\begin{remark}
Note that if $\epsilon$ is a Margulis constant, then $M_{\le 
\epsilon}$ consists of
Margulis tubes and cusps. Note that the same argument shows that, away from the
thin part of $M$ and an $\epsilon$ neighborhood of $\Gamma$, the 
diameter of $S$ can
be bounded by a constant depending only on $\chi(S)$ and $\epsilon$.
\end{remark}

The basic idea in the proof of Theorem~\ref{shrinkwrap_exists} is
to search for a least area representative of the isotopy class of the surface
$S$, subject to the constraint that the track of this isotopy does not
cross $\Gamma$. Unfortunately, $M\backslash \Gamma$ is not complete, so
the prospects for doing minimal surface theory in this manifold are remote.
To remedy this, we deform the metric on $M$ in a neighborhood of $\Gamma$ in
such a way that we can guarantee the existence of a least area surface
representative with respect to the deformed metric, and then take a limit
of such surfaces under a sequence of smaller and smaller such metric deformations.
We describe the deformations of interest below.

In fact, for technical reasons which will become apparent in 
\S \ref{CAT_property_subsection}, the deformations described below are not
quite adequate for our purposes, and we must consider metrics which are
deformed {\em twice} --- firstly, a mild deformation which satisfies
curvature pinching $-1 \le K \le 0$, and which is totally Euclidean in a neighborhood
of $\Gamma$, and secondly a deformation analogous to the kind described below in
Definition~\ref{deformed_metric_definition}, which is supported in this
totally Euclidean neighborhood. Since the reason for this ``double perturbation"
will not be apparent until \S \ref{CAT_property_subsection}, we postpone
discussion of such deformations until that time.

\begin{definition}[Deforming metrics]\label{deformed_metric_definition}
Let $\delta>0$ be such that $\Gamma$ is $\delta$--separated. Choose some small
$r$ with $r < \delta/2$.
For $t \in [0,1)$ we define a family of Riemannian metrics $g_t$ on $M$ in the
following manner. The metrics $g_t$ agree with the hyperbolic metric away from
some fixed tubular neighborhood $N_r(\Gamma)$.

Let $$h:N_{r(1-t)}(\Gamma) \to [0,r(1-t)]$$
be the function whose value at a point $p$ is the hyperbolic distance from $p$
to $\Gamma$. We define a metric $g_t$ on $M$ which agrees with the hyperbolic
metric outside $N_{r(1-t)}(\Gamma)$, and on $N_{r(1-t)}(\Gamma)$ is conformally
equivalent to the hyperbolic metric, as follows.
Let $\phi:[0,1] \to [0,1]$ be a $C^\infty$ bump function, which is equal to $1$ on the
interval $[1/3,2/3]$, which is equal to $0$ on the intervals $[0,1/4]$ and $[3/4,1]$,
and which is strictly increasing on $[1/4,1/3]$ and strictly decreasing on
$[2/3,3/4]$. Then define the ratio
$$\frac {g_t \text{ length element}} {\text{hyperbolic length element}} =
1 + 2\phi \Biggl( \frac {h(p)} {r(1-t)} \Biggr)$$
\end{definition}

We are really only interested in the behaviour of the metrics $g_t$ 
as $t \to 1$.
As such, the choice of $r$ is irrelevant. However, for convenience, we will
fix some small $r$ throughout the remainder of \S 1.

The deformed metrics $g_t$ have the following properties:

\begin{lemma}[Metric properties]\label{properties_of_deformed_metrics}
The $g_t$ metric satisfies the following properties:
\begin{enumerate}
\item{For each $t$ there is an $f(t)$ satisfying $r(1-t)/4 < f(t) < 3r(1-t)/4$
such that the union of tori
$\partial N_{f(t)}(\Gamma)$ are totally geodesic for the $g_t$ metric}
\item{For each component $\gamma_i$ and each $t$, the metric $g_t$ restricted
to $N_r(\gamma_i)$ admits a family of isometries which preserve $\gamma_i$ and
acts transitively on the unit normal bundle (in $M$) to $\gamma_i$}
\item{The area of a disk cross--section on $N_{r(1-t)}$ is
$O((1-t)^2)$.}
\item{The metric $g_t$ dominates the hyperbolic metric on $2$--planes.
That is, for all $2$--vectors $\nu$, the $g_t$ area of $\nu$ is at least as
large as the hyperbolic area of $\nu$}
\end{enumerate}
\end{lemma}
\begin{proof}
Statement (2) follows from the fact that the definition of $g_t$ has
the desired symmetries. Statements (3) and (4) follow from the fact that
the ratio of the $g_t$ metric to the hyperbolic metric is pinched 
between $1$ and $3$. Now, a radially symmetric circle linking $\Gamma$
of radius $s$ has length $2\pi\cosh(s)$ in the hyperbolic metric, and
therefore has length
$$2\pi\cosh(s)(1+2\phi(s/r(1-t)))$$
in the $g_t$ metric. For sufficiently small (but fixed) $r$, this function
of $s$ has a local minimum on the interval $[r(1-t)/4,3r(1-t)/4]$.
It follows that the family of radially symmetric tori linking a component of
$\Gamma$ has a local minimum for area in the interval $[r(1-t)/4,3r(1-t)/4]$.
By property (2), such a torus must be totally geodesic for the $g_t$ metric.
\end{proof}

\begin{notation}
We denote length of an arc $\alpha:I \to M$
with respect to the $g_t$ metric as $\tlength(\alpha(I))$, and area of a
surface $\psi:R \to M$ with respect to the $g_t$ metric as $\tarea(\psi(R))$.
\end{notation}

\subsection{Constructing the homotopy}

As a first approximation, we wish to construct surfaces in 
$M\backslash \Gamma$ which are globally least area with respect
to the $g_t$ metric. There are various tools for constructing least area
surfaces in Riemannian $3$-manifolds under various conditions, and subject
to various constraints. Typically one works in closed $3$-manifolds, but
if one wants to work in $3$-manifolds with boundary, the ``correct" boundary
condition to impose is {\em mean convexity}. A co-oriented
surface in a Riemannian $3$-manifold is said to be {\em mean convex} if 
the mean curvature vector of the surface always points to the negative side of
the surface, where it does not vanish. Totally geodesic surfaces and other
minimal surfaces are examples of mean convex surfaces, with respect to any
co-orientation. Such surfaces act
as {\em barriers} for minimal surfaces, in the following sense: suppose that
$S_1$ is a mean convex surface, and $S_2$ is a minimal surface.
Suppose further that $S_2$ is on the negative side of $S_1$. Then
if $S_2$ and $S_1$ are tangent, they are equal. 
One should stress that this barrier property is {\em local}. 
See \cite{MSY} for a more thorough discussion of barrier surfaces.

\begin{lemma}[Minimal surface exists]\label{surface_exists_for_each_t}
Let $M,\Gamma,S$ be as in the statement of Theorem~\ref{shrinkwrap_exists}.
Let $f(t)$ be as in Lemma~\ref{properties_of_deformed_metrics}, so that
$\partial N_{f(t)}(\Gamma)$ is totally geodesic with respect to the $g_t$ metric.
Then for each $t$, there exists an embedded surface $S_t$ isotopic in 
$M\backslash N_{f(t)}(\Gamma)$ to $S$, and which is globally
$g_t$--least area among all such surfaces.
\end{lemma}
\begin{proof}
Note that with respect to the $g_t$ metrics, the surfaces 
$\partial N_{f(t)}(\Gamma)$ described in
Lemma~\ref{properties_of_deformed_metrics} are totally geodesic, and 
therefore act as barrier surfaces. We remove the tubular neighborhoods
of $\Gamma$ bounded by these totally geodesic surfaces, and denote the
result $M\backslash N_{f(t)}(\Gamma)$ by $M'$ throughout the remainder of this
proof. We assume, after a small isotopy if necessary, 
that $S$ does not intersect $N_{f(t)}$ for any $t$, and therefore we can (and do)
think of $S$ as a surface in $M'$.
Notice that $M'$ is a complete Riemannian manifold with totally geodesic
boundary. We will construct the surface $S_t$ in $M'$, in the same isotopy class
as $S$ (also in $M'$).

If there exists a lower bound on the injectivity
radius in $M'$ with respect to the $g_t$ metric, then the
main theorem of \cite{MSY} implies that either such a globally least area
surface $S_t$ can be found, or $S$ is the boundary
of a twisted $I$--bundle over a closed surface in $M'$,
or else $S$ can be homotoped off every compact set in $M'$.

First we show that these last two possibilities cannot occur.
If $S$ is nonseparating in $M$, then it intersects some 
essential loop $\beta$ with algebraic intersection number
$1$. It follows that $S$ cannot be homotoped off $\beta$, and does not bound an
$I$--bundle. Similarly, if $\gamma_1,\gamma_2$ are distinct geodesics of
$\Gamma$ separated from each other by $S$, then the $\gamma_i$'s can be
joined by an arc $\alpha$ which has algebraic intersection
number $1$ with the surface $S$. The same is true of any $S'$ homotopic to $S$;
it follows that $S$ cannot be homotoped off the arc $\alpha$, nor does it bound
an $I$--bundle disjoint from $\Gamma$, and therefore does not bound an
$I$--bundle in $M'$.

Now suppose that the injectivity radius on $M'$ is not bounded below. We
use the following trick. Let $g_t'$ be obtained from the metric $g_t$ 
by perturbing
it on the complement of some enormous compact region $E$ so that it 
has a flaring end
there, and such that there is a barrier $g_t'$-minimal surface close
to $\partial E$, separating the complement of $E$ in $M'$ from $S$. Then by
\cite{MSY} there is a globally $g_t'$ least area surface $S_t'$, contained in
the compact subset of $M'$ bounded by this barrier surface. Since
$S_t'$ must either intersect $\beta$ or $\alpha$, by the Bounded Diameter
Lemma~\ref{bounded_diameter_lemma},
unless the hyperbolic area of $S_t' \cap E$ is very large, the diameter of
$S_t'$ in $E$ is much smaller than the distance from $\alpha$ or $\beta$ to
$\partial E$. Since by hypothesis, $S_t'$ is least area for the $g_t'$ metric,
its restriction to $E$ has hyperbolic area less than the hyperbolic 
area of $S$,
and therefore there is an {\it a priori} upper bound on its diameter in $E$.
By choosing $E$ big enough, we see that $S_t'$ is contained in the 
interior of $E$,
where $g_t$ and $g_t'$ agree. Thus $S_t'$ is globally least area
for the $g_t$ metric in $M'$, and therefore $S_t = S_t'$ exists for any $t$.
\end{proof}

The bounded diameter lemma easily implies the following:

\begin{lemma}[Compact set]\label{surfaces_in_compact}
There is a fixed compact set $E \subset M$ such that
the surfaces $S_t$ constructed in Lemma~\ref{surface_exists_for_each_t}
are all contained in $E$.
\end{lemma}
\begin{proof}
Since the hyperbolic areas of the $S_t$ are all uniformly bounded (by e.g. the
hyperbolic area of $S$) and are $2$--incompressible rel. $\Gamma$, they have
uniformly bounded diameter away from $\Gamma$ outside of Margulis tubes. 
Since for homological reasons they must intersect the
compact sets $\alpha$ or $\beta$, they can intersect
at most finitely many Margulis tubes. It follows that they are all contained in
a fixed bounded neighborhood $E$ of $\alpha$ or $\beta$, containing $\Gamma$.
\end{proof}

To extract good limits of sequences of minimal surfaces, one generally needs
{\it a priori} bounds on the area and the total curvature of the limiting
surfaces. Here for a surface $S$, the {\em total curvature} of $S$ is just the
integral of the absolute value of the (Gauss) curvature over $S$.
For minimal surfaces
of a fixed topological type in a manifold with sectional curvature 
bounded above, a
curvature bound follows from an area bound by Gauss--Bonnet. However, our
surfaces $S_t$ are minimal with respect to the $g_t$ metrics, which have no
uniform upper bound on their sectional curvature, so we must work 
slightly harder
to show that the the $S_t$ have uniformly bounded
total curvature. More precisely, we show that their
restrictions to the complement of any fixed
tubular neighborhood $N_\epsilon(\Gamma)$
have uniformly bounded total curvature.

\begin{lemma}[Finite total curvature]\label{finite_total_curvature}
Let $S_t$ be the surfaces constructed in Lemma~\ref{surface_exists_for_each_t}.
Fix some small, positive $\epsilon$.  Then the subsurfaces
$$S_t' := S_t \cap M \backslash N_\epsilon(\Gamma)$$
have uniformly bounded total curvature.
\end{lemma}
\begin{proof}
Having chosen $\epsilon$, we choose $t$ large enough so that $r(1-t) 
< \epsilon/2$.

Observe firstly that each $S_t$ has $g_t$ area less that the $g_t$ 
area of $S$, and
therefore hyperbolic area less that the hyperbolic area of $S$ for sufficiently
large $t$.

Let $\tau_{t,s} = S_t \cap \partial N_s(\Gamma)$ for small $s$.
By the coarea formula (see \cite{Fed},\cite{CoMi} page 8) we can estimate
$$\area(S_t \cap (N_{\epsilon}(\Gamma)\backslash N_{\epsilon/2}(\Gamma)))
\ge \int_{\epsilon/2}^\epsilon \length(\tau_{t,s}) \, ds$$
If the integral of geodesic curvature along a component $\sigma$ of 
$\tau_{t,\epsilon}$
is large, then the length of the curves obtained by isotoping $\sigma$
into $S_t \cap N_\epsilon(\Gamma)$ grows very rapidly, by the definition of
geodesic curvature.

Since there is an {\it a priori} bound on the hyperbolic area of 
$S_t$, it follows that there cannot be any long components of $\tau_{t,s}$ with
big integral geodesic curvature. More precisely, consider a long component $\sigma$
of $\tau_{t,s}$. For $l\in [0,\epsilon/2]$ the boundary $\sigma_l$ of the
$l$-neighborhood of $\sigma$ in $S_t \cap N_{\epsilon}(\Gamma)$ is contained
in $N_\epsilon(\Gamma)\backslash N_{\epsilon/2}(\Gamma)$. 
If the integral of the geodesic curvature along
$\sigma_l$ were sufficiently large for {\em every} $l$, 
then the derivative of the length of the $\sigma_l$ would be large for every $l$, and
therefore the lengths of the $\sigma_l$ would be large for all 
$l \in [\epsilon/4,\epsilon/2]$.
It follows that the hyperbolic 
area of the $\epsilon/2$ collar neighborhood of $\sigma$ in
$S_t$ would be very large, contrary to existence of an {\it a priori} upper bound on
the total hyperbolic area of $S_t$. 

This contradiction implies that for some $l$,
the integral of the geodesic curvature along $\sigma_l$ can be bounded from above.
To summarize, for each constant $C_1>0$ there is a constant $C_2>0$, such that for each
component $\sigma$ of $\tau_{t,\epsilon}$ which has length $\ge C_1$ there
is a loop
$$\sigma' \subset S_t \cap (N_\epsilon(\Gamma) \backslash 
N_{\epsilon/2}(\Gamma))$$
isotopic to $\sigma$ by a short isotopy, satisfying
$$\int_{\sigma'} \kappa \, dl \le C_2$$

On the other hand, since $S_t$ is $g_t$ minimal,
there is a constant $C_1>0$ such that each component $\sigma$ of 
$\tau_{t,\epsilon}$
which has length $\le C_1$ bounds a hyperbolic globally least area disk which
is contained in $M\backslash N_{\epsilon/2}(\Gamma)$.
For $t$ sufficiently close to $1$, such a disk is contained in 
$M\backslash N_{r(1-t)}(\Gamma)$, and therefore must actually be a subdisk of $S_t$.

By the coarea formula above, we can choose $\epsilon$ so that
$\length(\tau_{t,s})$ is {\it a priori} bounded. It follows that if $S_t''$
is the subsurface of $S_t$ bounded by the components of $\tau_{t,s}$ of
length $> C_1$ then we have {\it a priori} upper bounds on the
area of $S_t''$, on $\int_{\partial S_t''} \kappa \, dl$, and on 
$-\chi(S_t'')$.
Moreover, $S_t''$ is contained in $M \backslash N_{r(1-t)}$ where the metric
$g_t$ agrees with the hyperbolic metric, so the curvature $K$ of $S_t''$ is
bounded above by $-1$ pointwise, by Lemma~\ref{curvature_comparison}.
By the Gauss--Bonnet formula, this gives an {\it a priori} upper bound on
the total curvature of $S_t''$, and therefore on $S_t' \subset S_t''$.
\end{proof}

\begin{remark}
A more highbrow proof of Lemma~\ref{finite_total_curvature}
follows from Theorem 1 of \cite{S}, using the fact that the surfaces $S_t'$ are
locally least area for the hyperbolic metric, for $t$ sufficiently close to $1$
(depending on $\epsilon$).
\end{remark}

\begin{lemma}[Limit exists]\label{limit_exists}
Let $S_t$ be the surfaces constructed in Lemma~\ref{surface_exists_for_each_t}.
Then there is an increasing sequence
$$0 < t_1 < t_2 < \cdots$$
such that $\lim_{i \to \infty} t_i = 1$, and the $S_{t_i}$ converge on compact
subsets of $M \backslash \Gamma$ in the $C^\infty$ topology
to some $T' \subset M \backslash \Gamma$ with closure $T$ in $M$.
\end{lemma}
\begin{proof}
By definition, the surfaces $S_t$ have $g_t$ area bounded above
by the $g_t$ area of $S$. Moreover, since $S$ is disjoint from 
$\Gamma$, for sufficiently
large $t$, the $g_t$ area of $S$ is equal to the hyperbolic area of $S$. Since
the $g_t$ area dominates the hyperbolic area, it follows that the $S_t$ have
hyperbolic area bounded above, and by Lemma~\ref{finite_total_curvature}, for
any $\epsilon$, the restrictions of $S_t$ to $M \backslash N_\epsilon(\Gamma)$
have uniformly bounded finite total curvature.

Moreover, by Lemma~\ref{surfaces_in_compact}, each
$S_t$ is contained in a fixed compact subset of $M$. By standard
compactness theorems (see e.g. \cite{CiSc}) any infinite sequence 
$S_{t_i}$ contains a subsequence which converges on compact subsets of
$E \backslash \Gamma$, away from finitely many points where some subsurface
with nontrivial topology might collapse. That is, there might be isolated points
$p$ such that for any neighborhood $U$ of $p$, the intersection of $S_{t_i}$ with $U$
contains loops which are essential in $S_{t_i}$ for all sufficient large $i$.

But $S$ is $2$--incompressible rel. $\Gamma$, so in particular it is 
incompressible in $M \backslash \Gamma$, and no such collapse can take place. So
after passing to a subsequence, a limit $T'\subset M\backslash \Gamma$ 
exists (compare \cite{MSY}). Since
each $S_t$ is a globally least area surface in $M\backslash N_{f(t)}(\Gamma)$
with respect to the $g_t$ metric, it is a locally least area surface
with respect to the hyperbolic metric on $M \backslash N_{r(1-t)}(\Gamma)$.
It follows that $T'$ is locally least area in the hyperbolic metric,
properly embedded in $M \backslash \Gamma$, and we can define $T$ to 
be the closure of $T'$ in $M$.
\end{proof}

\begin{lemma}[Interpolating isotopy]\label{approximate_by_isotopies}
Let $\lbrace t_i \rbrace$ be the sequence as in Lemma~\ref{limit_exists}. Then
after possibly passing to a subsequence, there is
an isotopy $F:S \times [0,1) \to M\backslash \Gamma$ such that
$$F(S,t_i) = S_{t_i}$$
and such that for each $p \in S$ the track of the isotopy $F(p,[0,1))$ either
converges to some well--defined limit $F(p,1) \in M \backslash \Gamma$
or else it is eventually contained in $N_\epsilon(\Gamma)$ for any 
$\epsilon > 0$.
\end{lemma}
\begin{proof}
Fix some small $\epsilon$. Outside $N_\epsilon(\Gamma)$, the surfaces $S_{t_i}$
converge uniformly in the $C^\infty$ topology to $T'$. It follows that for
any $\epsilon$, and for $i$ sufficiently large (depending on $\epsilon$), 
the restrictions of $S_{t_i}$ and $S_{t_{i+1}}$ to the complement of $N_\epsilon(\Gamma)$
are both sections of the exponentiated 
unit normal bundle
of $T' \backslash N_\epsilon(\Gamma)$, and therefore we can isotope these subsets of
$S_{t_i}$ to $S_{t_{i+1}}$ along the fibers of the normal bundle. We wish to patch
this partial isotopy together with a partial isotopy supported in a small neighborhood
of $N_\epsilon(\Gamma)$ to define the correct isotopy from $S_{t_i}$ to $S_{t_{i+1}}$.

Let $Z$ be obtained from $N_\epsilon(\Gamma)$ by isotoping it slightly into
$M \backslash N_\epsilon(\Gamma)$ so that it is transverse to $T$, 
and therefore
also to $S_{t_i}$ for $i$ sufficiently large. For each $i$, we consider
the intersection
$$\tau_i = S_{t_i} \cap \partial Z$$ and observe that the limit satisfies
$$\lim_{i \to \infty} \tau_i = \tau = T \cap \partial Z$$

Let $\sigma$ be a component of $\tau$ which is inessential in 
$\partial Z$. Then
for large $i$, $\sigma$ can be approximated by $\sigma_i \subset 
\tau_i$ which are
inessential in $\partial Z$. Since the $S_{t_i}$ are 
$2$--incompressible rel. $\Gamma$,
the loops $\sigma_i$ must
bound subdisks $D_i$ of $S_{t_i}$. Since $\partial Z$ is a convex surface with
respect to the hyperbolic metric, and the $g_t$ metric agrees with 
the hyperbolic metric outside $Z$ for large $t$, 
it follows that the disks $D_i$ are actually contained in
$Z \backslash \Gamma$ for large $i$.
It follows that $D_i$ and $D_{i+1}$ are isotopic by an isotopy supported in
$Z \backslash \Gamma$, which restricts to a very small isotopy of $\sigma_i$ to
$\sigma_{i+1}$ in $\partial Z$.

Let $\sigma$ be a component of $\tau$ which is essential in $\partial 
Z$. Then so is
$\sigma_i$ for large $i$. Again, since $S$, and therefore $S_{t_i}$ is
$2$--incompressible
rel. $\Gamma$, it follows that $\sigma_i$ cannot be a meridian of 
$\partial Z$, and
must actually be a longitude. It follows that there is another essential
curve $\sigma'_i$ in each $\tau_i$,
such that the essential curves $\sigma'_i$ and $\sigma_i$ cobound a
subsurface $A_i$ in $S_{t_i} \cap Z \backslash \Gamma$. After passing 
to a diagonal
subsequence, we can assume that the $\sigma'_i$ converge to some 
component $\sigma'$
of $\tau$.

By $2$--incompressibility,
the surfaces $A_i$ are annuli. Note that there are {\em two} relative isotopy
classes of such annuli. By passing to a further diagonal subsequence, 
we can assume
$A_i$ and $A_{i+1}$ are isotopic in $Z \backslash \Gamma$ by an isotopy which
restricts to a very small isotopy of $\sigma_i \cup \sigma_i'$ to 
$\sigma_{i+1} \cup
\sigma_{i+1}'$ in $\partial Z$.

We have shown that for any small $\epsilon$ and any sequence $S_{t_i}$,
there is an arbitrarily large index
$i$ and infinitely many indices $j$ with $i<j$
so that the surfaces $S_{t_i}$ and $S_{t_j}$ are isotopic, and the isotopy can
be chosen to have the following properties:
\begin{enumerate}
\item{The isotopy takes $N_\epsilon(\Gamma) \cap S_{t_i}$ to
$N_\epsilon(\Gamma) \cap S_{t_j}$ by an isotopy supported in 
$N_\epsilon(\Gamma)$.}
\item{Outside $N_\epsilon(\Gamma)$, the tracks of the isotopy are contained in
fibers of the exponentiated normal bundle of $T'\backslash 
N_\epsilon(\Gamma)$.}
\end{enumerate}
Choose a sequence $\epsilon_i \to 0$,
and pick a subsequence of the $S_{t_i}$'s and relabel
so that $S_{t_i},S_{t_{i+1}}$ satisfy the properties above with respect to
$N_{\epsilon_i}(\Gamma)$. Then the composition of
this infinite sequence of isotopies is $F$.
\end{proof}

\begin{remark}
The reason for the circumlocutions in the statement of
Lemma~\ref{approximate_by_isotopies}
is that we have not yet proved that $T$ is a limit of the $S_t$ as 
maps from $S$
to $M$. This will follow in \S \ref{cone_subsection}, where we 
analyze the structure of
$T$ near a point $p \in \Gamma$, and show it has a well--defined tangent cone.
\end{remark}

\subsection{Existence of tangent cone}\label{cone_subsection}

We have constructed $T$ as a subset of $M$, and have observed that away from
$\Gamma$, $T$ is a minimal surface for the hyperbolic metric. We refer to the
intersection $T \cap \Gamma$ as the {\em coincidence set}. In general, one cannot
expect $T$ to be smooth along the coincidence set. However, we show that it
does have a well-defined {\em tangent cone} in the sense of Gromov, and this
tangent cone is in fact of a very special form. In particular, this is enough
to imply that $T$ exists as the image of a
map from $S$ to $M$, and we may extend the isotopy
$F:S \times [0,1) \to M$ to a homotopy $F:S \times [0,1] \to M$ with $T = F(S,1)$.

By a {\em tangent cone} we mean the following: at each point $p \in T\cap \Gamma$, 
consider the pair of metric spaces $(N_s(p),T_s(p))$ where $T_s(p)$
is the intersection $T_s = T \cap N_s(p)$. We rescale the metric on this pair by the factor
$1/s$. Then we claim that this sequence of (rescaled) pairs of metric spaces converge
in the Gromov-Hausdorff sense to a limit $(B,C)$ where $B$ is the unit ball in
Euclidean $3$-space, and $C$ is the cone (to the origin) over a great bigon in
the unit sphere. Here by a {\em great bigon} we mean the union of two spherical geodesics
joining antipodal points in the sphere.
In fact we do not quite show that $T$ has this structure, but rather that each
{\em local branch} of $T$ has this structure. Here we are thinking of the
{\em map} $F(\cdot,1):S \to M$ whose image is $T$, and by ``local branch" we mean the image
of a regular neigborhood of a point preimage.

\begin{lemma}[Tangent cone]\label{tangent_cone_exists}
Let $T$ be as constructed in Lemma~\ref{limit_exists}. Let
$p \in T \cap \Gamma$. Then near $p$, $T$ is a (topologically immersed) surface, 
each local branch of which has a well--defined tangent cone, which is the cone on
a great bigon.
\end{lemma}
\begin{proof}
We use what is essentially a curve--shortening argument. For each small
$s$, define
$$T_s = \partial N_s(p) \cap T$$
For each point $q \in T \backslash \Gamma$, we define $\alpha(q)$ to be the
angle between the tangent space to $T$ at $q$, and the radial geodesic through
$q$ emanating from $p$. By the coarea formula, we can calculate
$$\area(T \cap N_s(p)) = \int_0^s \int_{T_t} \frac 1 {\cos(\alpha)} dl \, dt
\ge \int_0^s \length(T_t) dt$$
where $dl$ denotes the length element in each $T_t$. Note that this estimate
implies that $T_t$ is rectifiable for a.e. $t$. We choose $s$ to be such a
rectifiable value.

Now, each component $\tau$ of $T_s$ is a limit of components $\tau_i$ of
$S_{t_i} \cap \partial N_s(p)$ for large $i$. By $2$--incompressibility of
the $S_{t_i}$, each $\tau_i$ is a loop bounding a subdisk $D_i$ of $S_{t_i}$
for large $i$.

Now, $\partial N_s(p)$ is convex in the hyperbolic metric, though not necessarily
in the $g_t$ metric. By cutting out the
disks $\partial N_s(p) \cap N_{r(1-t)}(\Gamma)$ and replacing them with the
disks $D^\pm$ orthogonal to $\Gamma$ which are totally geodesic in both the 
$g_t$ and the hyperbolic metrics, we can approximate $\partial N_s(p)$ by a surface 
$\partial B$ bounding a ball $B \subset N_s(p)$ which is convex in the $g_s$ metric for
all $s\ge t$. The ball $B$ is illustrated in figure~\ref{convex}.

\begin{figure}[ht]
\centerline{\relabelbox\small \epsfxsize 3.0truein
\epsfbox{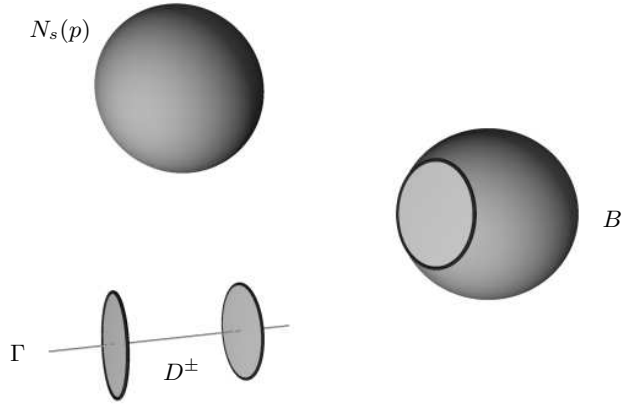}
\extralabel <-3.0truein, 2truein> {$N_s(p)$}
\extralabel <-2.3truein, 0.2truein> {$D^\pm$}
\extralabel <-3.1truein, 0.3truein> {$\Gamma$}
\extralabel <0.0truein, 1truein> {$B$}
\endrelabelbox}
\caption{The ball $B$ has boundary which is convex in both the hyperbolic
and the $g_s$ metrics for all $s\ge t$}
\label{convex}
\end{figure}

Note that after lifting $B$ to the universal cover, 
there is a retraction onto $B$ which
is length non-increasing, in both the $g_t$ and the hyperbolic metric.
This retraction projects along the fibers of the product structure on
$N_{r(1-t)}(\Gamma)$ to $D^\pm$, and outside $N_{r(1-t)}(\Gamma)$, it
is the nearest point projection to $\partial B\backslash D^\pm$.

Let $\tau_i'$ be the component of $S_{t_i} \cap \partial B$ approximating
$\tau_i$, and let $D_i'$ be the subdisk of $S_{t_i}$ which it bounds.

Then the disk $D_i'$ must be contained in $B$, or else we could decrease its
$g_t$ and hyperbolic area by the retraction described above.
The disks $D_i'$ converge to the component
$D \subset T'$ bounded by $\tau$, and the
hyperbolic areas of the $D_i'$  converge to the hyperbolic area of $D$.

Note that $B$ as above is really shorthand for $B_t$, since it depends on
a choice of $t$. Similarly we have $\tau_t$ and $D_t$. Since the component
$D_t \subset T'$ bounded by $\tau_t$ is contained in $B_t$ for all $t$,
the component $D \subset T'$ bounded by $\tau \subset \partial N_s$ is
contained in $N_s$, since $B_t \to N_s$ as $t \to 1$. So we can, and do,
work with $N_s(p)$ instead of $B$ in the sequel.

Now, let $D'$ be the cone on $\tau$ to the point $p$. $D'$ can be
perturbed an arbitrarily small amount to an embedded disk $D''$, and therefore
by comparing $D''$ with the $D_i'$, we see that the
hyperbolic area of $D'$ must be at least as large as that of $D$.
Note that this perturbation can be taken to move $D'$ off $\Gamma$, and
can be approximated by perturbations which miss $\Gamma$. Similar facts
are true for all the perturbations we consider in the sequel.

Since this is true for each component $\tau$ of $T_s$, by abuse of notation we
can replace $T$ by the component of $T \cap N_s(p)$ bounded by a 
single mapped in circle $\tau$. This will be the local ``branch" of the topologically
immersed surface $T$.
We use this notational convention for the remainder of the proof of the lemma.
Note that the inequality above still holds. It follows that we must have
$$\area(T \cap N_s(p)) \le \int_0^s \length(T_s) \frac {\sinh(t)} 
{\sinh(s)} dt =
\area(\text{cone on }T_s)$$

Now, for each sphere $\partial N_s(p)$, we let $\phi$ be the projection, along 
hyperbolic geodesics, to the unit sphere $S^2$ in the tangent space at $p$.
For each $t \in (0,1]$, define 
$$\|T_t\| = \length(\phi(T_t)) = \frac {\length(T_t)} {\sinh(t)}$$

It follows from the inequalities above that for some intermediate 
$s'$ we must have
$$\|T_{s'}\| \le \|T_s\|$$
with equality iff $T \cap N_s(p)$ is equal to the cone on $T_s$.

Now, the cone on $T_s$ is not locally least area for the
hyperbolic metric in $N_s(p) \backslash \Gamma$ 
unless $T_s$ is a great circle or geodesic bigon in $\partial N_s(p)$ 
(with endpoints on $\partial N_s(p) \cap \Gamma$), in which case the lemma is
proved. To see this, just observe that a cone has vanishing principal curvature
in the radial direction, so its mean curvature vanishes iff it is totally
geodesic away from $\Gamma$.

So we may suppose that for any $s$ there is some $s' < s$ such that
$\|T_{s'}\| < \|T_s\|$. Therefore we choose a sequence of values $s_i$ with
$s_i \to 0$ such that $\|T_{s_i}\| > \|T_{s_{i+1}}\|$, such that
$\|T_{s_i}\|$ converges to the infimal value of $\|T_t\|$ with $t \in 
(0,s]$, and
such that $\|T_{s_i}\|$ is the minimal value of $\|T_t\|$ on the interval
$t \in [s_i,1]$. Note that for any small $t$, the cone on $T_t$ has area
$$\area(\text{cone on } T_t) = 
\frac t 2 \length(T_t) + O(t^3) = \frac {t^2} 2 \|T_t\| + O(t^3)$$

The set of loops in the sphere with length bounded above by some constant,
parameterized by arclength, is {\em compact}, by the Arzela--Ascoli
theorem, and so we can suppose that the $\phi(T_{s_i})$
converge in the Hausdorff sense to a loop $C \subset S^2$.

\begin{claim}
$C$ is a geodesic bigon.
\end{claim}
\begin{proof}
We suppose not, and will obtain a contradiction.

We fix notation: for each $i$, let $C_i$ denote the inverse image $\phi^{-1}(C)$
under $\phi:\partial N_{s_i}(p) \to S^2$. So $C_i$ is a curve in $\partial N_{s_i}(p)$.
By the cone on $C_i$ we mean the union of the hyperbolic geodesic segments in $N_{s_i}(p)$
from $C_i$ to $p$. By the cone on $C$ we mean the union of the geodesic segments in
the unit ball in Euclidean $3$-space from $C \subset S^2$ to the origin. For each
$i$, we have an estimate
$$\area(\text{cone on } C_i) = s_i^2\area(\text{cone on } C) + O(s_i^3)$$

For each $i$, let $T^i$ denote the surface obtained from $T \cap N_{s_i}(p)$
by rescaling metrically by $1/{s_i}$. Then $T^i$ is a surface with boundary
contained in a ball of radius $1$ in a space of constant curvature $-s_i^2$.
Moreover, it enjoys the same least area properties as $T \cap N_{s_i}(p)$.

By the monotonicity property of the $\|T_{s_i}\|$ and the
coarea formula, we have an inequality
$$\lim_{i \to \infty}\area(T^i) \ge \area(\text{cone on } C)$$
On the other hand, since each $T^i$ is least area, we have an
estimate
$$\frac 1 2 \|T_{s_i}\|+ O(s_i) = 
\frac{\area(\text{cone on } T_{s_i})} {s_i^2} \ge \area(T^i)$$
It follows that the limit of the area of the $T^i$ is actually equal to
the area of the cone on $C$.

On the other hand, since the $\phi(T_{s_i})$ converge
to $C$, for sufficiently large $i$ we can find an immersed
annulus $A_i$ in $S^2$ with area $\le \kappa$ for any positive $\kappa$,
which is the track of a homotopy (in $S^2$) from $\phi(T_{s_i})$ to $C$. 
We let $\phi^{-1}(A_i)$ denote the corresponding annulus in $\partial N_{s_i}(p)$

We can build a new immersed surface
bounded by $T_{s_i}$ which is the union of this annulus $\phi^{-1}(A_i)$
with the cone on $C_i$. This surface can be perturbed an arbitrarily small amount,
away from $\Gamma$, to an embedded surface $F^i$. 
After rescaling $F^i$ by $1/s_i$, we get a surface
$G^i$ with the same boundary as $T^i$ of area equal to 
$$\area(G^i) = \area(\text{cone on } C) + \kappa + O(s_i)$$
Since $T^i$ is locally least area, it follows that for any $\kappa>0$,
for sufficiently large $i$ (depending on $\kappa$),
$$\area(\text{cone on } C) + 2\kappa \ge \area(G^i) \ge 
\area(T^i) \ge \area(\text{cone on } C) - \kappa$$

The surface $F^i$ contains a subsurface which is the cone on $C_i$.
Since by hypothesis, $C$ is not a geodesic bigon, the cone on $C$ can be
perturbed by a compactly supported perturbation to a surface whose area is smaller
than that of the cone on $C$ by some definite amount $\epsilon$. Similarly,
the cone on $C_i$ can be
perturbed by a compactly supported perturbation to a surface whose area is smaller than
the cone on $C_i$ by $\epsilon(s_i)^2$ where $\epsilon$ is independent of $i$. 
After rescaling by $1/s_i$, It follows that $G^i$ can
be perturbed by a compactly supported perturbation to $H^i$ with the
same boundary as $G^i$ and $T^i$, for which
$$\area(H^i) \le \area(G^i)-\epsilon$$
where $\epsilon$ is {\em independent} of $i$. Since $\kappa$ may be chosen as
small as we like, we choose $3\kappa< \epsilon$. Then for 
sufficiently large $i$ we get
$$\area(H^i) < \area(T^i)$$
which contradicts the least area property of $T^i$.
This contradiction shows that $C$ is actually a geodesic bigon, and completes
the proof of the claim.
\end{proof}

We now complete the proof of Lemma~\ref{tangent_cone_exists}.

Let $C \subset S^2$ be this geodesic bigon. Then inside an
$\epsilon$ neighborhood of $C$ in $S^2$, we can find a 
pair of curves
$C^\pm$, where $C^+$ is convex, and $C^-$ is convex except for two acute angles
on $T\Gamma \cap S^2$. For each $i$, let $C^\pm_i \subset \partial N_{s_i}(p)$
be the inverse image of $C^\pm$ under $\phi$.

The cone on $C^\pm_i$ is a pair of barrier surfaces in $N_s(p)$. 
In particular, once $\phi(T_{s_i})$ and $\phi(T_{s_{i+1}})$ are
both trapped between $C^+$ and $C^-$, the same is true of
$\phi(T_{s'})$ for all $s' \in [s_{i+1},s_i]$. This is enough to establish the
existence of the tangent cone.
\end{proof}

Notice that Lemma~\ref{tangent_cone_exists} actually implies that $T$ 
is a rectifiable
surface in $M$, which is a local (topological) embedding. In 
particular, this shows
that the isotopy $F:S \times [0,1) \to M$ constructed in
Lemma~\ref{approximate_by_isotopies} can be chosen to limit to a homotopy
$F:S \times [0,1] \to M$ such that $F(S,1) = T$.

\subsection{The thin obstacle problem}

From the proof of Lemma~\ref{tangent_cone_exists}, we see that $T$ exists as a
$C^0$ map which by abuse of notation we denote $u:T \to M$. 
One may immediately improve the regularity of $u$. From the construction of $T$,
it is standard to show that $u$ is actually in the Sobolev space $H^{1,2}$ 
--- that is, the
derivative $du$ is defined, and is $L^2$, in the sense of distribution;
see \cite{Morrey} for a rigorous definition of Sobolev spaces, 
and basic properties. 

To see this, observe that $u$ is a limit of maps $F(\cdot,t_{i}):S \to M$ which
are minimal for the $g_t$ metric, and therefore are $L^2$ energy minimizers
for the conformal structure on $S$ pulled back by $F(\cdot,t_{i})$.
If the set of conformal structures obtained in this way is precompact, 
one may extract a limit and therefore bound the $L^2$ norm of $du$ in terms
of the $L^2$ norms of the derivatives of any $F(\cdot,t_{i})$. How
can the sequence of conformal structures fail to be precompact? If and only
if the conformal structures degenerate by a neck pinch. But the
$2$--incompressibility of $S$ rel. $\Gamma$ implies that there is a lower
bound on the length of the image of any essential curve in $S$. It follows
that the $L^2$ norm of the derivative blows up along such a pinching neck,
contrary to the energy minimizing property. So no such degeneration can occur,
and $u$ is in $H^{1,2}$ as claimed. This argument is essentially contained in
\cite{SchoenYau} (see e.g. Lemma 3.1, page 134), 
and one may consult this paper for details.

We need to establish further regularity of $du$ along $\Gamma$ in the 
following sense.
Recall that we are calling $L := T \cap \Gamma$ the {\em coincidence set}. For 
each local sheet of $T$, we want $u$ to be $C^1$
along the interior of $L$ from either side, and at a 
non--interior point of $L$, we want $u$ to be $C^1$ on the nose.

Now, if $I \subset L$ is an interval, then the reflection principle 
(see \cite{Oss})
implies that each local sheet $T^+$ of $T$ with $\partial T^+ = I$
can be analytically continued to
a minimal surface across $I$, by taking another copy of $T^+$,
rotating it through angle $\pi$ along the axis $I$ and gluing it to 
the original
$T^+$ along $I$. It follows
that $du$ is real analytic from either side along the interior of 
$L$. Note that
if the tangent cone at a point $p$ is not literally a tangent {\em plane}, 
then an easy
comparison argument implies that $p$ is an interior point of the
coincidence set. See \cite{Nitsche} page 90 for a fuller discussion.

Non--interior points of $L$ are more difficult to deal with, and we 
actually want
to conclude that $du$ is continuous at such points.
Fortunately, this is a well--known problem in the theory of 
variational problems,
known as the {\em Signorini problem}, or the (two dimensional) {\em 
thin obstacle problem}.

In the literature, this problem is usually formulated in the following terms:

\begin{thin_problem}\label{thin_obstacle}
Let $\Omega$ be a bounded open subset of $\R^2$, and $A$ an oriented 
line contained
in $\Omega$. Let $\psi:A \to \R$ and $g:\partial \Omega \to \R$ be given, with
$g \ge \psi$ on $\partial \Omega \cap A$. Define
$$\K = \lbrace v \in g + H^{1,p}_0 \, | \, v \ge \psi \text{ on } A \rbrace$$
Minimize $$J(u) = \int F(x,u,\nabla u) dx$$
over $u \in \K$.
\end{thin_problem}
Here $H^{1,p}$ denotes the usual Sobolev space over $\Omega$ for the 
$L^p$ norm,
with zero boundary conditions.

The main conditions typically imposed on $F$ are sufficient regularity of $F$
and its partial derivatives (Lipschitz is usually enough) and {\em ellipticity},
meaning that the matrix $(F_{ik})_{i,k=1,2}$ of the second partial
derivatives of $F(x,u,\eta)$
with respect to $\eta \in \R^2$ is uniformly positive definite on compact
subsets of $\overline{\Omega} \times \R^{2+1}$
(see \cite{Frehse} page 281 for details). Roughly speaking, $F$ is elliptic
if the critical functions of the functional $J$ satisfy a
``mean value property" --- the value at each point is a weighted average of
the value at nearby points.

For example, if we want the graph of $u$ to be a
(Euclidean) minimal surface away from $\psi(A)$, then the formula for $F$ is
$F = (1 + |\nabla u|^2)^{1/2}$ which is certainly real analytic and elliptic.
The definition of $F$ for a nonparametric
minimal surface in exponential co--ordinates on hyperbolic space is 
more complicated,
but certainly $F$ is real analytic and elliptic in the sense of Frehse.

\begin{figure}[ht]
\centerline{\relabelbox\small \epsfxsize 5.0truein
\epsfbox{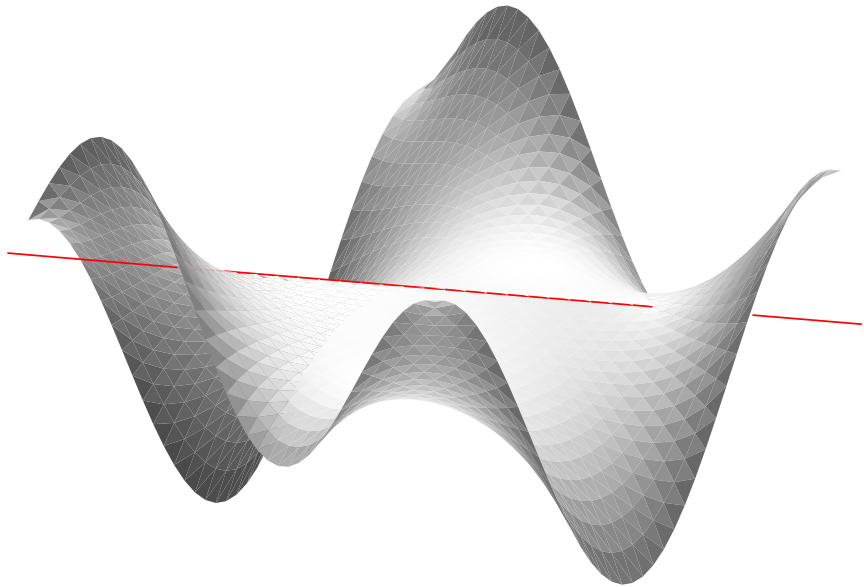}
\endrelabelbox}
\caption{The graph of a function solving the thin obstacle problem}
\label{thin_obstacle_figure}
\end{figure}

See Figure~\ref{thin_obstacle_figure} for an example of the graph
of a function solving the Dirichlet thin obstacle problem, where 
$\psi|_A$ is constant.
This surface is visually indistinguishable from the
graph of the function solving the unparameterized minimal surface thin obstacle
problem with the same boundary and obstacle data, but for computer 
implementation,
the Dirichlet problem is less computationally costly.

The next theorem establishes not only the desired continuity
of $\partial u$, but actually gives an estimate for the modulus of continuity.
The following is a restatement of Theorem 1.3 on page 26 of \cite{Rich}
in our context:

\begin{theorem}[Richardson \cite{Rich} Regularity of
thin obstacle]\label{obstacle_regularity}
Let $u$ be a solution to the thin obstacle problem for $F$
elliptic in the sense of Frehse and $p \in [1,\infty]$, and suppose
that $\partial \Omega,\psi,g$ are smooth. Then
$\partial u$ is continuous along $A$ in the tangent direction, 
one--sided continuous
in the normal direction on either side, and continuous in the normal 
direction at
a non--interior point. Furthermore, $\partial u$ is H\"older 
continuous, with exponent
$1/2$; i.e. the modulus of continuity of $\partial u$ is $O(t^{1/2})$.
\end{theorem}

\begin{remark}
Note that $C^{1+1/2}$ is actually best possible. Consider the function
$u: z \to \text{Im}(z^{3/2})$ for $u \in \mathbb{C}$ slit along the positive
real axis, where we take the branch which is
negative sufficiently close to the slit. This solves a thin obstacle problem for
the Dirichlet integral, and is only $C^{1+1/2}$ at $z=0$.
\end{remark}

\begin{remark}
For our applications, the fact 
that $u$ is $C^{1+1/2}$ is more than necessary. In fact, all we use 
is that $u$ is $C^1$. This is proved (with a logarithmic
modulus of continuity for $du$) by \cite{Frehse}, and (with a H\"older 
modulus of continuity for $du$) in arbitrary dimension by \cite{Kind}.
\end{remark}

We apply this theorem to our context:

\begin{lemma}[Regularity along coincidence set]\label{continuity_of_du}
For $u:T \to M$ defined as above, the derivative $du$ along local sheets of $T$
is continuous from each side along the coincidence set $L$, and continuous at
non--interior points.
\end{lemma}
\begin{proof}
If $p$ is an interior point of $L$, this follows by the reflection principle.
Otherwise, by Lemma~\ref{tangent_cone_exists} and the discussion above,
the tangent cone is a plane $\pi$ in the tangent space at $p$.

We show how to choose local co--ordinates in a ball
$B$ near each point $p \in L$ such that $B \cap \Gamma$ is the $x$--axis,
each local sheet of $T$ is the graph of a function $u: \Omega \to 
\R$, and $u$ is
non--negative along the $x$--axis. Let $\gamma = B \cap \Gamma$, and let
$\gamma'$ be another geodesic through $p$ orthogonal to $\gamma$ and tangent
to $\pi$. Let $\mF$ and $\mG$ be foliations of $B$ by totally geodesic
planes orthogonal to $\gamma$ and $\gamma'$ respectively. Then each 
leaf of $\mF$ is totally
geodesic for both the hyperbolic and the $g_t$ metric for all $t$, 
and each leaf of
$\mG$ is totally geodesic for the $g_t$ metric for sufficiently large $t$.
It follows that $T$ has no source or sink singularities with respect to
either foliation. Since $T \cap B$ is a (topological) disk, by
reasons of Euler characteristic it can have no saddle singularities either, and
therefore no singularities at all. We let $\mF$ and $\mG$ be level sets of
two co--ordinate functions on $B$. Define a third co--ordinate 
function to be (signed)
hyperbolic distance to the plane containing $\gamma$ and $\gamma'$,
and observe that $u$ is a graph in these co--ordinates.

It follows that $u$ solves an instance of
the thin obstacle problem, and by Theorem~\ref{obstacle_regularity} 
or by \cite{Frehse} or \cite{Kind}
the desired regularity of $du$ follows.
\end{proof}

\begin{remark}
The structure of the 
coincidence set is important to understand, and it has been
studied by various authors. Hans Lewy \cite{Lewy} showed that for $J$ 
the Dirichlet integral and $\psi$ analytic, 
the coincidence set is a finite union of points and intervals. 
Athanasopoulos \cite{Ath} proved the same result for the minimal 
surface question, for symmetric domain $\Omega$ and obstacle $A$, 
but his (very short and elegant) proof relies fundamentally on the 
symmetry of the problem, and we do not see how it applies in our context.

Note that if the Hausdorff dimension of the coincidence set is strictly 
$<1$, then since $T$ is $C^{1+1/2}$ (and therefore Lipschitz) along 
this coincidence set, the theory of removable singularities implies 
that $T$ is actually real analytic along $\Gamma$. It follows in this 
case that the coincidence set consists of a finite union of isolated 
points, and that $T$ is actually a minimal surface. See e.g. 
\cite{Car} for details.
\end{remark}

\begin{remark}
Existence results for the thin obstacle problem for minimal surfaces 
with analytic obstacles (see e.g. \cite{Rich}, \cite{Kind}, 
\cite{Nitsche}) gives an alternative
proof of the existence of the 
limit $T$. Given $S$, we can shrinkwrap $S$
near $\Gamma$ in small balls by using existence for the thin obstacle problem,
and away from $\Gamma$ by replacing small disks with least area embedded disks
with the same boundary. The argument of \cite{Hass_Scott} implies that for 
$S$ $2$--incompressible, this converges to a surface $T$.
\end{remark}

\subsection{$\CAT(-1)$ property}\label{CAT_property_subsection}

We have shown that $T$ satisfies all the properties of the
conclusion of Theorem~\ref{shrinkwrap_exists},
except that we have not yet shown 
that it is intrinsically $\CAT(-1)$.
In this subsection we show that after possibly replacing $T$ by
a new surface with the same properties, we can insist that $T$ 
is $\CAT(-1)$ with respect to
the path metric induced from $M$.

\begin{lemma}[CAT(-1) property]\label{intrinsically_negatively_curved}
After possibly replacing $T$ by a new immersed surface with 
the same properties, $T$ 
is $\CAT(-1)$ with respect to the path metric induced from $M$.
\end{lemma}
\begin{proof}
To show that $T$ is $\CAT(-1)$ we will 
show that there is no distributional
positive curvature concentrated 
along the coincidence set $L$. Since 
$T\backslash \Gamma$ is a 
minimal surface, the curvature of $T$ is bounded above
by $-1$ on this subset. It will follow by Gauss--Bonnet 
that $T$ is $\CAT(-1)$.

We first treat a simpler problem in Euclidean 
$3$--space, which we denote by $\R^3$. 
Let $\Sigma$ be an embedded 
surface in $\R^3$ which is $C^3$ outside a 
subset $X$ which is 
contained in a geodesic $\gamma$ in $\R^3$, and
which is $C^1$ along 
$X$ from either side along the interior of $X$, and
$C^1$ at non--interior points of $X$. Then we claim,
for each subsurface $R 
\subset \Sigma$ with $C^3$ boundary
$\partial R \subset \Sigma 
\backslash X$, that 
$$\int_{R\backslash X} K_\Sigma = 2\pi\chi(R) - 
\int_{\partial R} \kappa \, dl$$
Compare Lemma~\ref{Gauss_Bonnet}.

In other words, we want to show that $X$ 
is a ``removable singularity'' for $R$,
at least with respect to the 
Gauss--Bonnet formula.

Let $\phi:R\backslash X \to S^2$ denote the 
Gauss map, which takes each point
$p \in R$ to its unit normal, in 
the unit sphere of $S^2$. Then $K_\Sigma$ is
the pullback of the area 
form by $\phi$. Let $\overline{R}$ denote the
completion of 
$R\backslash X$ with respect to the path metric. Then
$\overline{R}$ is obtained from $R$ by cutting it open along 
each interval in $R \cap X$, and sewing in two copies 
of the interval thereby removed. Notice
that there is a natural forgetful map 
$\overline{R} \to R$.

By the assumptions about the regularity of $R$, the Gauss map 
$\phi$ actually extends to a {\em continuous} map 
$\phi:\overline{R} \to S^2$. Moreover, since $X$ is
contained in a geodesic $\gamma$ of $\R^3$, the image
$\phi(\overline{R} \backslash R)$ is contained in a 
great circle $C$ in $S^2$.

For each boundary component $\tau$ of 
$\overline{R}$, we claim that the map $\phi: \tau \to C$ has degree 
zero. For, otherwise, by a degree argument,
there are points $p^\pm \in \overline{R}$ 
which map to the same point in $p$,
for which 
$\phi(p^+) = -\phi(p^-)$, and the graphs of $\phi(p^+)$ and 
$-\phi(p^-)$ locally have a nonzero algebraic intersection number. 
It follows that the local sheets of $R$ from either side must 
actually intersect along $p$, contrary to the fact that
$\Sigma$ is embedded. It follows that we can sew in a 
disk to $\overline{R}$ along each boundary component to get a surface $\overline{R}'$ 
homeomorphic to $R$, with $\partial \overline{R}' = \partial R$, and 
extend $\phi$ to $\phi:\overline{R}' \to S^2$
by mapping each such disk into $C$.

Now, the surface $\overline{R}'$ can be perturbed 
slightly in a neighborhood of $X$ to a new surface $\overline{R}''$ 
which is $C^3$ in $\R^3$, in such a way that the Gauss map of 
$\overline{R}''$ is a perturbation of $\phi$.
So the usual Gauss--Bonnet formula (Lemma~\ref{Gauss_Bonnet}) shows 
that
$$\int_{\overline{R}''} K = 2\pi\chi(R) - \int_{\partial R} \kappa \, dl$$
But $\int_{\overline{R}''} K$ is just the integral of 
the area form on $S^2$ pulled back by the Gauss map; it follows that 

$$\int_{\overline{R}''} K = \int_{S^2} \text{degree}(\phi(\overline{R}'))$$
and
$$\int_{S^2} \text{degree}(\phi(\overline{R}'')) = \int_{S^2} 
\text{degree}(\phi(\overline{R}'))$$
since one map is obtained from the other by a small perturbation supported away
from the boundary. Since the measure of $C$ is zero, this last integral is just
equal to
$$\int_{S^2\backslash C} \text{degree}(\phi(\overline{R}')) = \int_{R\backslash X} K$$
and the claim is proved.

Now we show how to apply this to our shrinkwrapped surface $T$.
We use the following trick. Let $j_t$ with $t \in [0,1)$ 
be a family of metrics on $M$, 
conformally equivalent to the hyperbolic metric, which
agree with the hyperbolic metric outside $N_{r(1-t)}$, which are Euclidean on 
$N_{r(1-t)/2}$, and which have curvature pinched between $-1$ and 
$0$, and are rotationally and translationally symmetric along the 
core geodesic. Then we let $T_t$
be the surface obtained by 
shrinkwrapping $S$ {\em with respect to the $j_t$ metric}.
I.e. we let $g_{s,t}$ be a family of metrics as in 
Definition~\ref{deformed_metric_definition}
which agree with the 
$j_t$ metric outside $N_{r(1-t)(1-s)}$, construct minimal 
surfaces $S_{s,t}$ as in Lemma~\ref{surface_exists_for_each_t}, 
and so on, limiting to
the immersed surface $T_t$ which is minimal for 
the $j_t$ metric on $M \backslash \Gamma$,
and $C^{1+1/2}$ along $T_t \cap \Gamma$. 
Arguing locally as above, we see that small 
subsurfaces of $T_t$ contained in the Euclidean tubes $N_{r(1-t)/2}$ 
satisfy Gauss--Bonnet in the complement of the coincidence set. By 
Lemma~\ref{curvature_comparison},
the surfaces $T_t$ all have 
curvature bounded above by $0$, and bounded above by $-1$
outside $N_{r(1-t)}$. By Gauss--Bonnet for geodesic triangles, 
$T_t$ is $\CAT(0)$, and
actually $\CAT(-1)$ outside $N_{r(1-t)}$.

Now take 
the limit as $t \to 1$. Some subsequence of the surfaces $T_t$ 
converges to a
limit which by abuse of notation we denote $T$. Note that this
is not necessarily the same as the surface $T$ constructed in previous
sections, but it enjoys the same properties. Again, by Gauss--Bonnet for geodesic 
triangles, the limit is actually $\CAT(-1)$, and the lemma is 
proved.
\end{proof}

This completes the proof of Theorem~\ref{shrinkwrap_exists}.

\begin{problem}
Develop a simplicial or PL theory of 
shrinkwrapping.
\end{problem}

\begin{remark}
Since this paper appeared in preprint form, Soma has developed some elements
of a PL theory of shrinkwrapping; see \cite{Soma_shrink}. This theory proves
a PL analogue of Theorem~\ref{shrinkwrap_exists}. 
\end{remark}\section{The Main Construction Lemma}

The purpose of this section is to state the main construction 
Lemma~\ref{construction_lemma}
and show how it follows easily from Theorem~\ref{shrinkwrap_exists}. 

\subsection{Shrinkwrapping in covers}\label{cover_subsection}

Let $N$ be a complete, orientable, parabolic free hyperbolic
$3$--manifold, and let $\Gamma$ be a finite collection of pairwise disjoint
simple closed geodesics in $N$, just as in 
the statement of Theorem~\ref{shrinkwrap_exists}.
For the purposes of introducing the Main Construction Lemma, we will assume
that $N$ has a single end $\E$.
We consider the family of $g_t$ and $g_{s,t}$ metrics, as in
Definition~\ref{deformed_metric_definition} and 
Lemma~\ref{intrinsically_negatively_curved}.

Suppose there is an embedded surface $\partial W$ in $N\backslash \Gamma$ which
separates off the end of $N$ from a compact submanifold
$W \subset N$. Let $X$ be a covering space
of $W$ (possibly infinite). The preimage of the geodesics $\Gamma$ are
a collection of locally finite geodesics $\hat{\Gamma} \subset X$,
some of which might be finite, and some infinite. Let $\Delta \subset \hat{\Gamma}$
be some nonempty collection, consisting entirely of simple closed geodesics.
Then we can consider a second surface
$S \subset X\backslash \Delta$ and can ask whether it is possible to
shrinkwrap $S$ rel. $\Delta$. Notice that we cannot directly apply
Theorem~\ref{shrinkwrap_exists} because the hyperbolic manifold $X$ is not
{\em complete}, and therefore a shrinkwrap representative of $S$ might not
exist. However, we note that for each metric $g_t$ on $N$,
we get a $g_t$ locally least area representative $\partial W_t$ isotopic to $\partial W$.
The submanifold $W_t$ of $N$ bounded by $\partial W_t$ lifts to a covering space $X_t$
which is homeomorphic to $X$.
The metric $g_t$ pulls back to a metric on $X_t$ which by abuse of notation
we also refer to as $g_t$. Then $\partial X_t$ which is a lift of $\partial W_t$
is $g_t$ locally least area, and therefore acts as a barrier surface.
It follows that we can find, for each $t$, a surface
$S_t$ in the isotopy class of $S$ in $X_t\backslash N_{f(t)}(\Delta)$,
which is globally $g_t$-least area among all such surfaces (compare
with the statement of Lemma~\ref{surface_exists_for_each_t}).

The theory of shrinkwrapping developed in \S 1 goes through almost identically
for the surfaces $S_t$ with one important exception: the metric $g_t$
on $X_t$ does {\em not} agree with the hyperbolic metric away from $\Delta$
and $\partial X_t$, but rather is deformed along the other geodesics 
$\hat{\Gamma}\backslash \Delta$. It follows that we should take care to
analyze the quality of the surfaces $S_t$ and their limit $S'$ near components
of $\hat{\Gamma}\backslash \Delta$.

Fortunately the situation is as simple as it could be:

\begin{lemma}[Superfluous geodesics invisible]\label{geodesics_invisible}
With notation and definitions as above, in a neighborhood of a point $p$ on
$\hat{\Gamma}\backslash \Delta$, the surface $S'$ is a locally least
area surface for the hyperbolic metric.
\end{lemma}
\begin{proof}
Suppose not. Then there is a compactly supported perturbation $F$ of $S'$
which agrees with $S'$ outside a fixed neighborhood of $p$, and which has
strictly less {\em hyperbolic} area than $S'$, so that
$$\area(S') - \area(F) \ge \epsilon$$
for some positive constant $\epsilon$. After another small perturbation of
$F$ to $F'$,
which can be taken to increase the hyperbolic area as little as required,
we can assume that $F'$ is transverse to $\Gamma$ near $p$, intersecting it
in $n$ points for some finite $n$, and satisfying
$$\area(S') - \area(F') \ge \epsilon/2$$
By property (3) of the $g_t$ metric
(see Lemma~\ref{properties_of_deformed_metrics}) the $g_t$ area of
$F'$ is at most equal to the hyperbolic area plus $nC(1-t)^2$, for some constant
$C$ independent of $t$. For sufficiently small $t$, 
$$nC(1-t)^2 < \epsilon/2$$
and therefore the $g_t$
area of $F'$ is less than the {\em hyperbolic} area of $S'$, which is less than
the $g_t$ area of $S'$, thereby contradicting the {\em global} $g_t$ minimality
of $S'$ in its isotopy class in $X_t\backslash N_{f(t)}(\Delta)$.

This contradiction proves the lemma.
\end{proof}

A similar argument holds for the $g_{t,s}$ metric in place of the $g_t$ metric,
and therefore by means of Lemma~\ref{geodesics_invisible} we can shrinkwrap
in covers, obtaining $\CAT(-1)$ surfaces in the limit.

\begin{remark}
One should think of Lemma~\ref{geodesics_invisible} as a kind of ``removable
singularity" theorem for transverse obstacles. Compare with the following
physical experiment: one knows from experience that very
thin needles can be pushed through 
soap bubbles without popping them or distorting their geometry. (Try it!)
\end{remark}

\subsection{The main construction lemma}

We now state and prove the main construction lemma. The context of this Lemma
is the same as that of \S \ref{cover_subsection}: we want to shrinkwrap a certain surface
in a cover, using the boundary of that cover as a barrier surface. See Figure~\ref{21}
for an idealized depiction of $T'$ and $S$ in $W$ and $X$ in the case that $W$ is
a handlebody.

\begin{lemma}[Main construction lemma]\label{construction_lemma}
Let $\E$ be an end of the complete 
open orientable parabolic free hyperbolic
$3$--manifold $N$ with 
finitely generated fundamental group. Let $W \subset N$
be a submanifold such that $\partial W \cap \inte(N)$ separates $W$ from 
$\E$. Let $\Delta_1 \subset N \backslash \partial W$ be a finite 
collection of simple closed geodesics with 
$\Delta = \inte(W) \cap \Delta_1$ a non--empty
proper subset of $\Delta_1$. Suppose further 
that $\partial W$ is $2$--incompressible
rel. $\Delta$.

Let $G$ be a finitely generated subgroup of $\pi_1(W)$, and let $X$ be the 
covering space of $W$ corresponding to $G$. Let $\Sigma$ be the 
preimage of $\Delta$ in $X$, and $\hat{\Delta} \subset \Sigma$ a 
subset which maps homeomorphically onto $\Delta$
under the covering projection, and let $B \subset \hat{\Delta}$ be a nonempty
union of geodesics. Suppose there exists an embedded closed surface
$S \subset X \backslash B$
that is $2$--incompressible rel. $B$ in $X$, which 
separates every component of $B$ from $\partial X$.

Then $\partial W$ can be homotoped to a $\Delta_1$--minimal surface which, by abuse 
of notation, we call $\partial W'$, and the map of $S$ into $N$ given 
by the covering projection is homotopic to a map whose image $T'$ is 
$\Delta_1$--minimal. Also, $\partial W'$ (resp. $T'$) can be 
perturbed by an arbitrarily small perturbation to be
an embedded (resp. smoothly immersed) surface $\partial W_t$ (resp. $T_t$) 
bounding $W_t$ with the following properties:
\begin{enumerate}
\item{There exists an isotopy from 
$\partial W$ to $\partial W_t$ which never crosses
$\Delta_1$, and which induces an isotopy from $W$ to $W_t$, and a 
corresponding deformation of hyperbolic manifolds $X$ to $X_t$ which 
fixes $\Sigma$ pointwise.}
\item{There exists an isotopy from $S$ to $S_t \subset X_t$ 
which never crosses $B$, such that $T_t$ is the 
projection of $S_t$ to $N$.}
\end{enumerate}
\end{lemma}
\begin{proof}
\begin{figure}[ht]
\centerline{\relabelbox\small 
\epsfxsize 5.0truein
\epsfbox{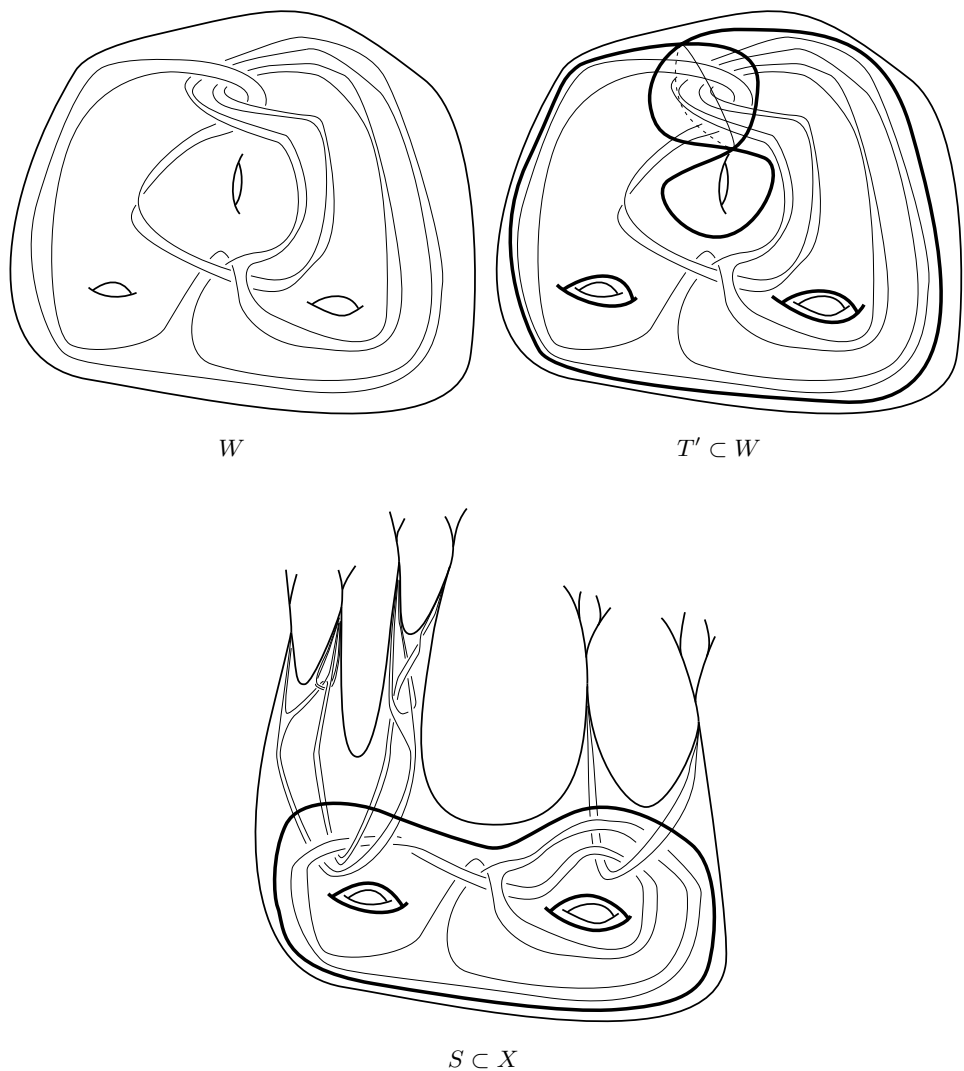}
\extralabel <-3.9truein, 3.3truein> {$W$}
\extralabel <-1.5truein, 3.3truein> {$T' \subset W$}
\extralabel <-2.7truein, 0.1truein> {$S \subset X$}
\endrelabelbox}
\caption{The surfaces $T'$ and $S$ in $W$ and $X$ respectively}
\label{21}
\end{figure}
The proof is reasonably straightforward, 
given the work in \S \ref{shrink} and \S \ref{cover_subsection}.
First, we obtain $\partial W'$ from $\partial W$ by shrinkwrapping rel. 
$\Delta_1$. Since $\Delta = \inte(W) \cap \Delta_1$ is a nonempty and 
proper subset of $\Delta_1$,
$\partial W$ satisfies the hypotheses of Theorem~\ref{shrinkwrap_exists}, and
therefore $\partial W'$ exists, and satisfies the desired properties.

For each $t = t_i$ in our approximating sequence, the metric $g_t$ on $W_t$
lifts to a metric on $X_t$ which, by abuse of notation, we also 
call $g_t$. 

Then with respect to the $g_t$ metric, $\partial X_t$ acts as a barrier
surface, and we can find a $g_t$ locally least area surface 
$S_t \subset \inte(X_t)$ which is $g_t$ globally least area
in the isotopy class of $S$ in $X_t\backslash N_{f(t)}(B)$,
by \cite{MSY}, just as in the proof of Lemma~\ref{surface_exists_for_each_t}.
Note that $S$ is necessarily homologically essential in $X_t\backslash N_{f(t)}(B)$,
since it separates each component of $B$ from $\partial X_t$ by hypothesis,
and therefore any surface isotopic to $S$ in $X_t\backslash N_{f(t)}(B)$
must intersect a fixed compact arc $\alpha$ from $B$ to $\partial X_t$.

The immersed surfaces $T_t \subset N$ 
are obtained by mapping $S_t$ to $W_t$ by the covering 
projection. After passing to a further subsequence of values
$t = t_i$, the limit of the $T_t$ surfaces exists as a map from $S$ to $N$, with 
image $T'$, by the argument of Lemma~\ref{limit_exists} applied 
locally. The regularity of $T'$ locally along $\Delta_1$ follows from 
the argument of Lemma~\ref{continuity_of_du}, since that argument is 
completely local. 

As in Lemma~\ref{intrinsically_negatively_curved}, we can repeat the
construction above with the $g_{s,t}$ metrics and obtain a limit $T'$ with
the desired regularity.

It follows that $T'$ is $\Delta_1$--minimal. Notice 
that some local sheets of $T'$ are actually {\em minimal} (in the 
hyperbolic metric) near geodesics
in $\Delta_1$, corresponding to subsets 
of $S_t$ in $X_t$ crossing components
of $\Sigma \backslash B$, by Lemma~\ref{geodesics_invisible}.
In any case, $T'$ is intrinsically $\CAT(-1)$, and the
theorem is proved.
\end{proof}

\subsection{Nonsimple geodesics}

When we come to consider hyperbolic manifolds with 
parabolics, we need to treat the
case that the geodesics $\Delta_1$ might not be simple. But there 
is a standard trick to reduce this case to the simple case, at the 
cost of slightly perturbing the
hyperbolic metric.

Explicitly, suppose $\Gamma \subset M$ is as in the 
statement of Theorem~\ref{shrinkwrap_exists}
except that some of the components 
are possibly not simple. Then for every
$\epsilon >0$ there exists
a perturbation $g$ of the hyperbolic metric on $M$ 
in a neighborhood of $\Gamma$ with
the following properties:

\begin{enumerate}
\item{The new metric $g$ agrees with the hyperbolic 
metric outside $N_\epsilon(\Gamma)$}
\item{With respect to the metric $g$, the 
curves in $\Gamma$ are homotopic to
a collection of simple geodesics $\Gamma'$}
\item{The metric $g$ is hyperbolic (i.e. has constant curvature $-1$)
on $N_{\epsilon/2}(\Gamma')$}
\item{The metric $g$ is $1+\epsilon$--bilipschitz 
equivalent to the hyperbolic metric, and 
the sectional curvature of the $g$ metric is pinched between 
$-1-\epsilon$ and $-1 + \epsilon$.}
\end{enumerate}
The existence of such a metric $g$ follows from 
Lemma 5.5 of \cite{Ca}. To make an orthopedic comparison: think
of the nonsimple geodesics $\Gamma$ as a collection of unnaturally fused
bones in $N_\epsilon(\Gamma)$; the bones are broken, reset, 
and heal as simple geodesics in the new metric.

It is clear that the methods of \S \ref{shrink} apply equally 
well to the metric $g$,
and therefore shrinkwrapping can be done with respect to the 
metric $g$, producing a surface which is intrinsically 
$\CAT(-1+\epsilon)$.

In fact, since such a metric exists for each 
$\epsilon$, we can take
a sequence of such metrics $g_\epsilon$ for 
each small $\epsilon > 0$,
produce a shrinkwrapped surface 
$T_\epsilon$ for each such $\epsilon$, and take
a limit $T$ as $\epsilon \to 0$ which is 
intrinsically $\CAT(-1)$, and which can be 
approximated by embedded surfaces, isotopic 
to $S$, in the complement of
$\Gamma\backslash C$ where $C$ is a finite subset 
of geodesics whose cardinality
can be {\it a priori} bounded above in terms of the 
genus of $S$. We will not be
using this stronger fact in the sequel, 
since the existence of a $\CAT(-1+\epsilon)$
surface is quite enough 
for our purposes.\section{Asymptotic Tube Radius and Length}

By \cite{Bo} an end of a complete hyperbolic $3$--manifold $N$ is 
geometrically infinite if and only if there exists an exiting 
sequence of closed geodesics.
In this chapter we show that if $\pi_1(N)$ is parabolic free, then the geodesics
can be chosen to be $\eta$--separated for some $\eta$; in particular, all are simple.

\begin{definition}
Let $N$ be a complete hyperbolic 
3-manifold with geometrically infinite end $\mE$. Define the
$\mE$-\emph{asymptotic tube radius} to be the supremum over all 
sequences $\lbrace \gamma_i \rbrace$ of  closed geodesics exiting 
$\mE$, of 
$$\limsup_{i \to \infty} \; \text{tube radius}(\gamma_i)$$
Similarly define the $\mE$-\emph{asymptotic length} to be infimum
over all sequences $\lbrace \gamma_i \rbrace$ 
as before of
$$\liminf_{i \to \infty} \; \text{length}(\gamma_i)$$
We will drop the prefix $\mE$ when the end in question 
is understood.\end{definition}

\begin{proposition}\label{bonahon} If 
$\mE$ is a geometrically infinite end of the complete hyperbolic 
3-manifold $N$ without parabolics, then asymptotic tube radius $>1/4$ 
asymptotic length. If asymptotic length = 0, then asymptotic tube 
radius = $\infty$. There exists a uniform lower bound
$\eta=0.025$ to the asymptotic tube radius of a geometrically infinite end
of a complete parabolic free hyperbolic $3$--manifold.
\end{proposition}

\begin{proof}  Meyerhoff \cite{Me} defines a monotonically
decreasing function $r:(0,0.1]\to [0.3,\infty)$ such that if 
$\gamma$ is a closed geodesic in $N$ and $\length(\gamma)\le  t$, then
tube radius$(\gamma)\ge r(t)$.  Furthermore lim$_{t\to 0}
r(t)=\infty$.  Therefore, the second statement of Proposition
\ref{bonahon} follows from \cite{Me} and the third  follows from
the first statement and \cite{Me}. (Actually, Proposition
\ref{log3} will show that $\log(3)/2$ is a lower bound.)

Now suppose that asymptotic length = $L\in [0.1,\infty)$.  Then there exists a 
sequence $\{\gamma_i\}$ exiting $\mE$ such that
length$(\gamma_i)\to L$. As in \S5 \cite{G2}, if 
tube radius $(\gamma_t) \le \frac{1}{4} \length(\gamma_i)$, then there exists a geodesic 
$\beta_i$ homotopic to a curve which is a union of a segment of 
$\gamma_i$ and an orthogonal arc from $\gamma_i$ to itself, and each 
of these segments has length $\le \length(\gamma_i)/2$.
By straightening these segments and using the law of cosines, we see
that if length$(\gamma_i)\ge 0.099$, then 
$\length(\beta_i) < \length(\gamma_i)-0.02$. Thus if
$\lim \sup$ tube radius $(\gamma_i) < L/4$, 
there exists a sequence $\{\beta_i\}$
such that $\lim \inf \length(\beta_i) \le L-0.02$ where $\beta_i$ is as above. 
Since $\infty>L$, $\{\beta_i\}$ must exit the same end as 
$\{\gamma_i\}$, which is a contradiction.

Now suppose that asymptotic length is infinite and
$\{\gamma_i\}$ is an exiting sequence such that length$(\gamma_i)\to
\infty$.  Given $R\ge 10$ we produce a new exiting sequence
$\{\sigma_i\}$ with tube radius 
$(\sigma_i)>R$ for all $i$.  If possible let $\alpha_i$ be a
smallest segment of $\gamma_i$ such that there is a geodesic path $\beta_i$ 
connecting $\partial\alpha_i$, $\length(\beta_i)\le 10R$ and 
$\beta_i$ is not homotopic to $\alpha_i$ rel
endpoints.  If  $\alpha_i$ does not exist, then tube
radius$(\gamma_i)\ge 5R$.  So let
us assume that for all $i$, $\alpha_i$ exists.  Note that 
$\length(\alpha_i)\to\infty$ or else the concatenations $\{\alpha_i*\beta_i\}$ are
homotopic to an exiting sequence 
of bounded length geodesics.  Also, asymptotic length infinite implies that as
$i\to\infty$, the injectivity radius of points of $\gamma_i$ (and hence
$\alpha_i\cup\beta_i$) $\to \infty$.  Therefore for
$i$ sufficiently large we can 
assume that $\length(\alpha_i)>10R$, $\length(\beta) = 10R$, 
and both of the angles between $\beta_i$ and $\alpha_i$ are at least
$\pi/2$.  The geodesic $\sigma_i$ homotopic to the curve 
obtained by concatenating $\alpha_i$ and $\beta_i$ lies within 
distance 2 of $\alpha_i\cup \beta_i$ and for the most part lies 
extremely close.  Indeed, if $A$ is an immersed least area annulus in $N$ with
$\partial A=\alpha_i*\beta_i\cup\sigma_i$, then the Gauss--Bonnet formula 
(Lemma~\ref{Gauss_Bonnet_corners})
implies that area$(A)\le \pi$.  Since the intrinsic curvature
of $A$ is $\le -1$, it follows that for $i$ sufficiently large, no point $a$ of
$A$ can be at distance $1$ from $\alpha_i\cup \beta_i$ and distance at least $1$
from $\sigma_i$, for the area of the disc of radius 1 about $a\in A$ would be 
$\ge \pi$. 

If tube radius$(\sigma_i)\le R$, then there would  be an arc
$\tau_i$ connecting points of $\sigma_i$ such that 
$\length(\tau_i)\le 2R$ and
$\tau_i$ cannot be
homotoped rel endpoints into $\sigma_i$.  Thus, for $i$ sufficiently large 
one finds new essential geodesic
paths $\beta_i^\prime$ of length $\le 10R $ with endpoints
in $\alpha_i-\partial \alpha_i$.  This contradicts the minimality
property of $\alpha_i$. 
\end{proof}

\begin{proposition} \label{log3} If $\mE$ is an end of the
complete,  orientible, hyperbolic 3-manifold $N$ and
$\pi_1(N)$ has no parabolic elements, then the $\mE$-asymptotic tube 
radius$>\log(3)/2$.\end{proposition}

\begin{remark}  We will not be using Proposition \ref{log3} in 
this paper.\end{remark}

\begin{proof}  Let $\{\gamma_i\}$ be a sequence of geodesics exiting $\mE$
such that
$$\lim_i \text{length}(\gamma_i)=l=\mE\text{--asymptotic length}.$$

If $l$ is small, i.e. $l\le 0.978$, then for $i$ sufficiently large,
tube radius $(\gamma_i) \ge \log(3)/2$ by \cite{Me} as explained in
Proposition 1.11\cite{GMT}.   If
$l$ is large, then for $i$ sufficiently large tube radius$({\gamma_i})\ge \log(3)/2$ by
Proposition \ref{bonahon}.  A hyperbolic geometry argument,  slightly more
sophisticated than the one cited above shows that
$l \ge 1.289785$ suffices.  Indeed, the proof of Proposition 1.11
\cite{GMT} shows that there exists $\epsilon > 0$ such that if  
$\length(\gamma_i) \ge 1.289785$, then either 
tube radius$(\gamma_i)\ge\log(3)/2$ or there 
exists an essential closed curve
$\kappa_i$ such that length$(\kappa_i)\le \length(\gamma_i)-\epsilon$ and
$d(\gamma_i,\kappa_i)\le l$.  If $\kappa^*_i$ denotes the geodesic homotopic to $\kappa_i$,
then $\{\kappa^*_i\}$  exits $\mE$ and lim inf length$(\kappa_i^*)\le l-\epsilon$ which
contradicts the fact that $l$ is asymptotic length.

It follows from \cite{GMT},
\cite{GMT}, \cite{KJ}, \cite{Li} and \cite{CLLM} that if $\delta$ 
is a shortest geodesic in a complete hyperbolic
3-manifold $N$, then either tube radius $\delta\ge \log(3)/2$ or $N$ is a closed hyperbolic
3-manifold.  (See Conjecture 1.31 \cite {GMT}.)  Therefore, if each 
$\gamma_i$ is a shortest length geodesic in $N$, then the proof of
Proposition \ref{log3} follows.

Assuming that asymptotic tube radius
$<\log(3)/2$, we will derive a contradiction using techniques which
require an understanding of
\S 1 \cite{GMT}.  Nevertheless,  the punchline follows exactly as in two
paragraphs above. Here is the idea.  Associated to each
$\gamma_i$ there is a 2--generator subgroup of $\pi_1(N)$ defined as
follows.  When viewed as acting on
$\BH^3$, one generator $f_i$ is a shortest translation along  a lift of
$\gamma_i$ and the other generator
$w_i$ takes that lift  to a nearest translate.  After
passing to a subsequence,  there exists 
$\epsilon>0$, $K<\infty$ such that for each $i$,   
there exists a closed curve
$\kappa_i$ such that $d(\kappa_i,\gamma_i)\le K$ and length$(\kappa_i)<\length\gamma_i-\epsilon$. 
Here $\kappa_i$ represents an element in the group generated by $f_i$ and $w_i$.  
   
Here are the details.   Given $\{\gamma_i\}$ there exist sequences 
$\{A_i\}, \{A'_i\}$ where $A_i$ is a lift of
$\gamma_i$ to $\BH^3$ and $A'_i=w_i(A)$
is a nearest $\pi_1(N)$ translate of 
$A$, where $w_i\in \pi_1(N)$.  By Definition 1.8 \cite{GMT} associated to
$f_i$ and $w_i$ there is a triple of complex numbers $(L_i,D_i,R_i)$ where
$\length(f_i)=\text{Re}(L_i)$ and
$\text{Re}(D_i)=d(A,A')$.  By compactness, after passing to a subsequence,
$\{(L_i,D_i,R_i)\}$ converges
to $(L,D,R)$, where  $\text{Re}(L)=l$.  Again by Definition 1.8,
$(L,D,R)$ gives  rise to a marked
2--generator group $\langle f,w \rangle$ where $f_i\to f$ and $w_i\to
w$. By Lemma 1.13 \cite{GMT}
we can assume that $(L,D,R)$  and the various $(L_i,D_i,R_i)$ lie in
the parameter space 
$\mathcal{P}$ defined in 1.11 \cite{GMT}. It
cannot lie in one of the 7 exceptional regions given in Table 1.2
\cite{GMT}, or  else by Chapter 3 
\cite{GMT}, \cite{KJ}, \cite{Li} and \cite{CLLM}, it 
and $(L_i,D_i,R_i)$ correspond to a closed hyperbolic
3-manifold for $i$ sufficiently large, for it is shown in these papers that a
neighborhood of each exceptional region corresponds to a unique
closed hyperbolic 3-manifold as conjectured in 1.31 \cite{GMT}. 
This implies that $N$ is covered by a closed 3-manifold, which is a
contradiction.  

The proof of Proposition 1.28 \cite{GMT} shows that
if $(L,D,R)$ does not lie in an exceptional region, then there
exists a killer word
$u(f,w)$ in $f$ and $w$ as defined in 1.18 \cite{GMT}.  This means that
$\length(u(f,w))<\length(f)$ or if 
$A=\axis(f)$, then $d(A,u(f,w)A)<\text{Re}(D)=d(A,w(A))$.   Therefore, for
$i$ sufficiently large either $\length(u(f_i,w_i))<\length(f_i)$ or
$d(A_i,u(f_i,w_i)A_i)<\text{Re}(D_i)$.  The latter  cannot happen since
$w_i$ was chosen to take $A_i$ to a nearest translate.

Since the nonexceptional points
of $\mathcal{P}$ are covered by finitely many compact regions and
each region has a killer word, it follows that for correct choice of
killer word, reduction of length is uniformly bounded
below  by some constant $\epsilon$. 

Since $\pi_1(N)$ has no parabolics, $u(f_i,w_i)$  corresponds to a
hyperbolic element and hence a geodesic 
$\sigma_i \subset N$.  If $u(f,w)$ is loxodromic, then the 
corresponding geodesic $\tilde\sigma\subset \BH^3$ is of bounded
distance from $A$.  Therefore, for all $i$, $d(\sigma_i, A_i)$
is uniformly bounded and hence $\{\sigma_i\}$ 
exits $\mE$.  Thus, asymptotic length $\le l-\epsilon$ which
is a contradiction. If $u(f,w)$ is parabolic,  then
$\length(\sigma_i)\to 0$,  $\{\sigma_i\}$ exits
the same end as $\{\gamma_i\}$ and hence asymptotic length equals zero.   To see that 
$\{\sigma_i\}$ exits $\mE$, note that  in
$\BH^3,  u(f,w)$ takes
a point $x$ to $y$ where $d(x,y)<l/4$ and hence for $i$ 
sufficiently large there are essential closed curves of length $<l/4$ 
passing within $2d(x,A)$ from $\gamma_i$.   
\end{proof} 

\begin{question}  What is the maximal lower bound for 
the asymptotic tube radius of a
geometrically infinite end $\mE$ of a 
complete, orientable, hyperbolic manifold with
finitely generated fundamental group, both in the cases that 
$\mE$ is parabolic free or not?\end{question}

\begin{question}  What is the upper bound for 
asymptotic length of a geometrically infinite end $\mE$?  It follows 
from Theorem \ref{main} that there is an upper bound which is
a function of rank$(\pi_1(\mE))$.  \end{question}

\section{Canary's Theorem}\label{canary section}
In this section we give a proof 
of Canary's theorem (Theorem~\ref{canary})
when $N$ is parabolic 
free.  
Our proof of Theorem \ref{main} will closely parallel this 
argument.

\begin{theorem}[Canary]\label{canary} If $\mE$ is a 
topologically tame end of the complete, orientable, hyperbolic 
3-manifold $N=\BH^3/\Gamma$, where $\Gamma$ has no 
parabolic elements, then either $\mE$ is geometrically finite 
or there exists a sequence of $\CAT(-1)$  surfaces exiting
the end.  If $\mE$ is parametrized by $S\times 
[0,\infty)$, then these surfaces are homotopic to
surfaces of the form $S\times t$, via a homotopy supported in
$S\times [0,\infty)$.
\end{theorem}

\begin{figure}[ht]
\centerline{\relabelbox\small \epsfxsize 4.0truein
\epsfbox{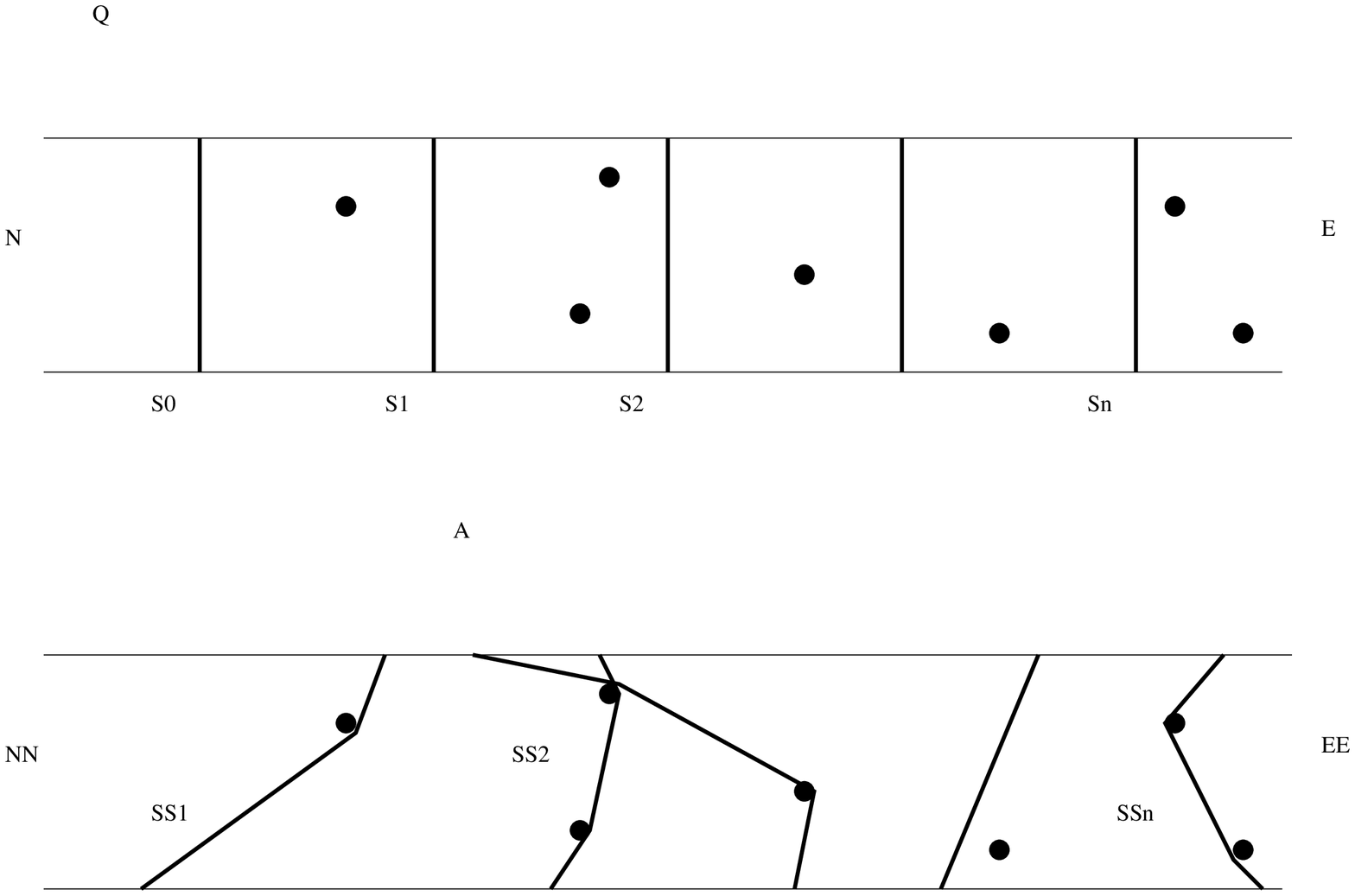}
\relabel {N}{$N$}
\relabel {NN}{$N$}
\relabel {E}{$\mE$}
\relabel {EE}{$\mE$}
\relabel {S0}{$S\times 0$}
\relabel {S1}{$S\times 1$}
\relabel {S2}{$S\times 2$}
\relabel {Sn}{$S\times n$}
\relabel {SS1}{$S_1$}
\relabel {SS2}{$S_2$}
\relabel {SSn}{$S_n$}
\relabel {Q}{Q: how can one find an exiting sequence of $\CAT(-1)$ surfaces?}
\relabel {A}{A: shrinkwrap!}
\endrelabelbox}
\caption{A schematic depiction of shrinkwrapping in action}
\label{41}
\end{figure}

\begin{proof}  It suffices to consider 
the case that $\mE$ is
geometrically infinite.   By Proposition~\ref{bonahon} there exists a
sequence of pairwise disjoint 
$\eta$-separated simple closed
geodesics $\Delta=\{\delta_i\}$ 
exiting $\mE$.   Assume that $\Delta$ and the
parametrization of $\mE$ are chosen so that for all $i\in \BN,\ $
$\delta_i\subset S\times (i-1,i)$.  Let $g$=genus$(S)$,
$\Delta_i=\{\delta_1,\cdots, \delta_i\}$ and let $\{\alpha_i\}$ be
a locally finite collection of embedded proper rays in $\mE$ such 
that $\partial \alpha_i \in \delta_i$.

An idea used repeatedly, in various guises, throughout 
this paper is the following.  If $R$ is a closed oriented surface and 
$T$ is obtained by shrinkwrapping $R$ rel the geodesics $\Delta_R$, 
then $R$ is homotopic to $T$ via a homotopy which does not meet 
$\Delta_R$, except possibly at the last instant.
Therefore, if $\delta_i\subset \Delta_R$ and $\langle R,\alpha_i \rangle =1$, then 
$T\cap \alpha_i\neq\emptyset$ and if $T\cap \delta_i=\emptyset$, 
then $\langle T,\alpha_i \rangle =1$.
Here $\langle \cdot , \cdot \rangle$ denotes algebraic intersection number.

\vskip 12pt

\noindent{\bf Warm up Case.}  Each $S\times i$ is 2-incompressible in
$N\backslash \Delta$. (E.g. $N= S\times \BR$)
\newline \newline
\emph{Proof.}  Apply Theorem~\ref{sample theorem} to
shrinkwrap $S \times i$
rel. $\Delta_{i+1}$ to a $\CAT(-1)$ surface $S_i$.  Since
$\langle S\times i,\alpha_i \rangle=1$,
$S_i\cap \alpha_i\neq\emptyset$.    Since
$\{\alpha_i\}$ is locally finite,
the Bounded  Diameter Lemma implies
that the $S_i$'s must exit $\mE$.   Therefore
for $i$ sufficiently large $S_i\subset \mE$ and $\langle S_i,\alpha_1 
\rangle=1$;
and hence the projection of $S_i$ into
$S\times 0$, (given by the product structure on $\mE$) is a degree-1
map between surfaces of the same genus.  Since such maps are
homotopic to homeomorphisms, we see that $S_i$ can be homotoped
within $\mE$ to a homeomorphism onto $S\times 0$.  See Figure~\ref{41} for
a schematic view.
\end{proof}

\noindent{\bf General Case.} (E.g. $N$ is an open 
handlebody.) \newline
\newline
\emph{Proof.}  Without loss of generality we can assume
that every closed orientable
surface separates $N$, (see Lemma~\ref{end-manifold} and
Lemma~\ref{core}).  We  use a purely
combinatorial/topological argument to find a particular sequence
of smooth surfaces exiting $\mE$.  We then shrinkwrap these surfaces and
show that they have 
the desired escaping and homological properties.

Fix $i$.  If possible, compress
$S\times i$, via a compression which either misses 
$\Delta$ or crosses $\Delta$
once say at $\delta_{i_1}\subset \Delta_i$.  If possible, compress
again via a compression meeting
$\Delta\backslash \delta_{i_1}$ at most once say 
at $\delta_{i_2}\subset \Delta_i$. After at most $n\le 2g-2$ such 
operations and deleting 2-spheres we  obtain embedded connected
surfaces $S^i_1, \cdots, S^i_{i_r}$, none of  which is a
2--sphere and each of which is 2--incompressible 
rel $\Delta_{i+1}\backslash \{\delta_{i_1}\cup\cdots\cup\delta_{i_n}\}$.
For each fixed $i$, each $\delta_j$ $(j \le i)$ with at most $2g-2$ exceptions
is separated from $\mE$ by exactly one surface $S^i_k$.  Call \emph{Bag}$^i_k$ the 
region separated from $\mE$ by $S^i_k$.   
Note that all compressions in the passage of $S_i$ to
$\{S^i_1,\cdots,S^i_{i_r}\}$ are on the non $\mE$-side.

\begin{figure}[ht]
\centerline{\relabelbox\small \epsfxsize 3.0truein
\epsfbox{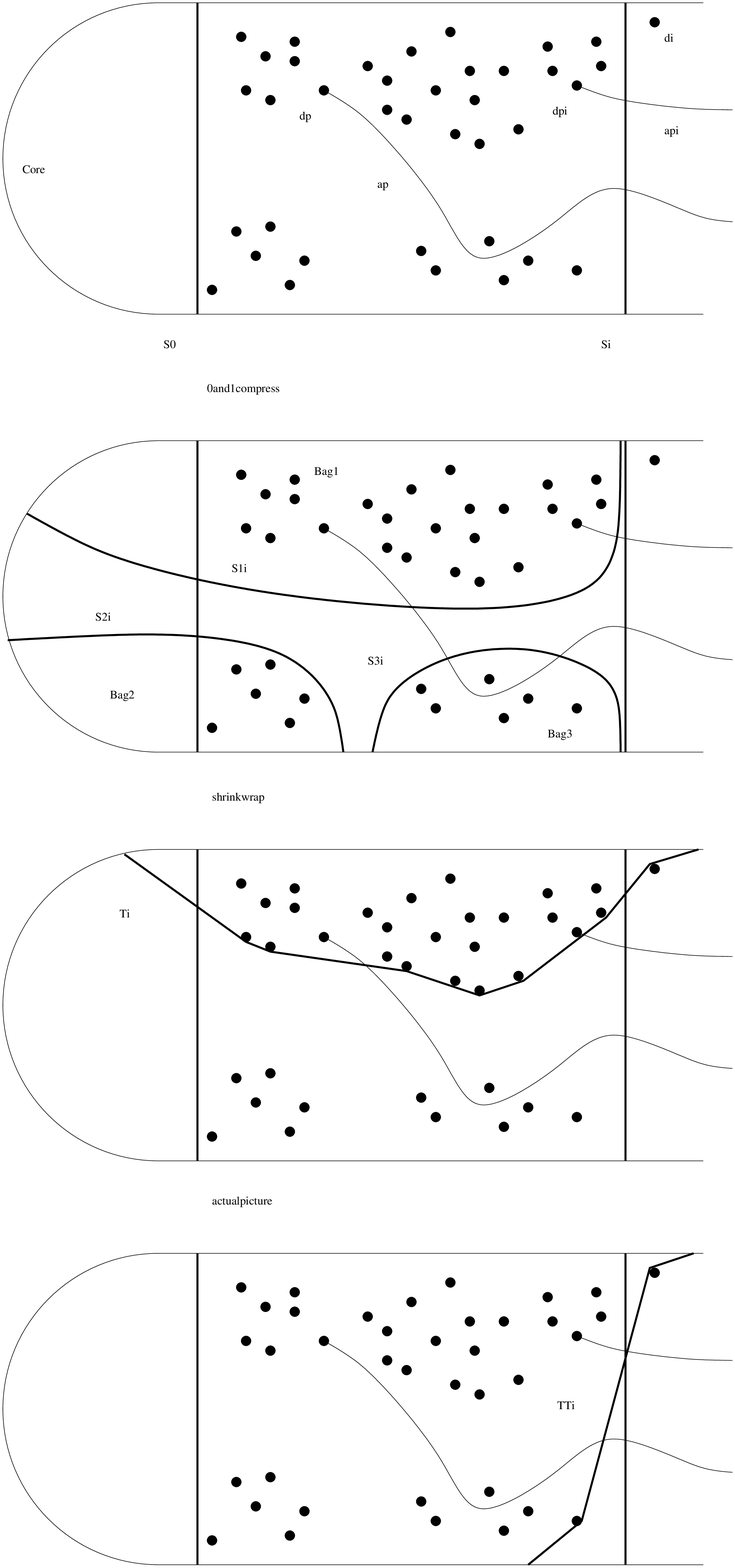}
\relabel {Core}{Core of $N$}
\relabel {S0}{$S\times 0$}
\relabel {Si}{$S\times i$}
\relabel {dp}{$\delta_p$}
\relabel {dpi}{$\delta_{p(i)}$}
\relabel {di}{$\delta_{i+1}$}
\relabel {ap}{$\alpha_p$}
\relabel {api}{$\alpha_{p(i)}$}
\relabel {0and1compress}{0 and 1 compress $S\times i$ to get the $S_j^i$}
\relabel {S1i}{$S_1^i$}
\relabel {S2i}{$S_2^i$}
\relabel {S3i}{$S_3^i$}
\relabel {Bag1}{$\text{Bag}_1^i$}
\relabel {Bag2}{$\text{Bag}_2^i$}
\relabel {Bag3}{$\text{Bag}_3^i$}
\relabel {Ti}{$T_i$}
\relabel {TTi}{$T_i$}
\relabel {shrinkwrap}{Shrinkwrap $S_1^i$ to get $T^i$}
\relabel {actualpicture}{Actual picture}
\endrelabelbox}
\caption{The Bounded Diameter Lemma and the intersection number argument show that
$S\times i$ undergoes no compression, and $T_i$ actually separates all of $\Delta_i$ from 
$\mE$.}
\label{42}
\end{figure}

Since each $i_r\le g$, we can find 
a $p\in \BN$ and for each $i$, a
reordering of the $S^i_j$'s  (and their bags) so 
that for infinitely many
$i\ge p$, $\delta_p\in \Bag^i_1$; furthermore, if for each $i$ such that 
$\delta_p\in \Bag^i_1$, we denote by $p(i)$ the maximal index such that
 $\delta_{p(i)}\in \Bag^i_1$,  then the set $\{p(i)\}$ is unbounded.  By
Theorem \ref{sample theorem},
$S^i_1$ is homotopic rel 
$\Delta_{i+1}\backslash\{\delta_{i_1},\cdots, \delta_{i_n}\}$ 
to a  $\CAT(-1)$ surface $S_i$.
Since the collection $\{\alpha_{p(i)}\}$ is infinite and locally finite, the Bounded
Diameter Lemma implies that a subsequence of these
$S^i_1$'s must exit $\mE$. Call this subsequence $T_1,T_2,\cdots$, where $T_i$ is
the shrinkwrapped $S^{n_i}_1$. Therefore, for $i$ sufficiently large, $T_i$ must
lie in $S\times (p,\infty)$ and   
$\langle T_i,\alpha_p \rangle=\langle S^{n_i}_1,\alpha_p,\rangle =1$. 
Therefore, projection of $T_i$ to $S\times p$ is degree 1.  
This in turn implies that genus $T_i=g$ and $T_i$ 
can be homotoped within $\mE$ to a homeomorphism
onto $S\times 0$. See Figure~\ref{42} for another schematic view.
\qed\newline

\begin{remark}  This argument shows that for $i$ sufficiently
large, $S\times i$ is already 2--incompressible in 
$N\backslash \Delta_i$.  Also, given
any $\eta$-separated collection of exiting geodesics a
sufficiently large finite subset is 2-disc 
busting. Actually using the technology of the last chapter, 
this statement holds for any sequence of exiting closed geodesics. 
\end{remark}

The proof of Theorem \ref{main} follows a
similar strategy.  Here is the outline
in the case that $N$ has a single end 
$\mE$ and no parabolics.  Given
a sequence of $\eta$-separated exiting simple closed 
geodesics $\Delta=\{\delta_i\}$ we pass to  subsequence (and possibly 
choose $\delta_1$ to have finitely many components) and find a 
sequence of connected embedded surfaces denoted
$\{\partial W_i\}$ such that for each $i,\ \partial W_i$ 
separates $\Delta_i=\delta_1\cup\delta_2\cup\cdots\cup\delta_i$ 
from $\Delta-\Delta_i$ and is  2-incompressible rel $\Delta_i$. It is 
{\it a priori} possible that the $\partial W_i$'s do not exit 
$\mE$. If $W_i$ denotes the compact region split off by 
$\partial W_i$, then after  possibly deleting
an initial finite set of $W_i$'s 
(and adding the associated $\delta_i$'s to $\delta_1$)
we find a compact 3-manifold $D\subset W_1$ which is a core for
$\mW=\cup W_i$.

We next find an immersed genus $\le g$ surface $T_i$,  which
homologically separates off a subset $B_i$ of $\Delta_i$ from $\mE$.
For infinitely many $i$,  $B_i$ includes a fixed $\delta_p$ 
and for these $i$'s the set $\{p(i)\}$ is unbounded, where $p(i)$ is 
the largest index of a $\delta_k\subset B_i$.
The surface $T_i$ separates $B_i$ from the rest in the sense that
$T_i$ lifts to an embedded surface $\hat T_i$ in the 
$\pi_1(D)$-cover $\hat W_i$ of $W_i$ and
in that cover $\hat T_i$ separates a lift $\hat B_i$ from 
$\partial \hat W_i$, the preimage
of $\partial W_i$.  The argument to 
this  point is purely topological and applies
to any 3-manifold with finitely generated fundamental group. In the 
general case,
$\{\partial W_i\}$ will not be an exiting sequence.

Next we shrinkwrap $\partial W_i$ rel $\Delta_{i+1}$ to a
$\CAT(-1)$ surface which we continue to call $\partial W_i$.   Then we
homotop $\hat T_i$ rel $\hat\Delta_i$ to a $\CAT(-1)$ surface in the induced
$\hat W_i$ and let $T_i$ denote the projection of $T_i$ to $N$.
The point of shrinkwrapping $\partial W_i$ is  that
$\partial \hat W_i$ is now a barrier which prevents $\hat T_i$ from 
popping out of $\hat W_i$ during the subsequent shrinkwrapping (compare with \S 2).
We use the $\delta_{p(i)}$'s to show
that, after passing to a subsequence, the
$T_i$'s exit $\mE$. We use $\delta_p$ to show that for $i$ sufficiently large,
$T_i$ homologically separates $\mE$ from a Scott core of $N$.

We have outlined the strategy.  For purposes of exposition, the construction 
of the $T_i$'s in \S 6 is slightly different from the sketch above.

In \S \ref{parabolic} we make the necessary embellishments 
to handle the parabolic case.

The next chapter develops the theory of end reductions which enables
us to define the submanifolds $W_i$.

\section{End Manifolds and End Reductions}\label{end reductions}

In this section, we prove a structure theorem for the 
topology of an end of a $3$--manifold with finitely generated 
fundamental group. A reference for basic $3$--manifold topology is
\cite{Hempel}.

The first step is to replace our original 
manifold with a $1$--ended manifold $M$ with the homotopy type of a 
bouquet of circles and closed orientable surfaces. We then prove 
Theorem~\ref{refined}, the {\em infinite end engulfing theorem}, 
which says that given an exiting seqence of homotopically non
trivial simple closed curves we can  pass to a subsequence
$\Gamma$ and find a submanifold
$\mathcal{W}$, with finitely generated fundamental 
group, containing $\Gamma$ which has the following 
properties:
\begin{enumerate}
\item{$\mW$ can be exhausted by codimension 0 compact submanifolds
$W_i$ whose boundaries are $2$--incompressible rel  $\Gamma\cap W_i$.}
\item{$\mW$ has a core which lies in $W_1$.}
\end{enumerate}

This completes the preliminary step 
in the proof of  Theorem~\ref{main}, as explained at the end of
\S 4.  The proof of Theorem~\ref{main} itself is in \S 6.

In what follows we will assume that all 3-manifolds are
orientable and irreducible.

\begin{lemma}\label{end-manifold}  If $\mathcal{E}$ is 
an end of an open Riemannian $3$-manifold $M^\prime$ with finitely 
generated fundamental group, then $\mathcal{E}$ is isometric to the 
end of a 1-ended 3-manifold $M$ whose (possibly
empty) boundary is a finite union of closed orientable surfaces.
A core of $M$ is obtained by attaching 1-handles to the components of
$\partial M$, unless $\partial M=\emptyset$, in which case a core is a 
1-complex and $M$=$M'$.
\end{lemma}

\begin{proof}  A thickened Scott core $C$ \cite{Sc} of $M^\prime$ is a
union of 1-handles (possibly empty) attached to a compact 3-manifold $X$ 
with incompressible boundary. Split $M^\prime$ along all the 
boundary components of $X$ and let $M$ be the component which 
contains $\mathcal{E}$. 
\end{proof}

\begin{remark}\label{end-manifold remark}$M$ is a 
submanifold of $M'$. $M$ is isometric to a submanifold 
$\hat M$ of the covering of $M^\prime$ corresponding to the inclusion
$\pi_1(M) \to \pi_1(M')$, and the inclusion $M\to \hat M$ is a 
homotopy equivalence.
\end{remark}

\begin{definition} Call a finitely generated  group a {\em free/surface
group} if it is a free product of orientable surface groups and a free group.
Call a 1-ended, irreducible, orientable, 3-manifold $M$ 
an \emph{end--manifold} if it has a compact (possibly empty) boundary 
and a compact core of the
form $\partial M\times I\cup 1$-handles if $\partial M\neq\emptyset$
or a handlebody if $\partial M=\emptyset$. 

\end{definition}
Note that $\pi_1(M)$ is a free/surface group for 
$M$ an end--manifold.

\begin{lemma}  \label{rank equality}If $G$ is 
a subgroup of a free/surface group, then
its $\pi_1$-rank equals its 
$H_1$-rank, both in $\BZ$ 
and $\BZt$-coeficients.\qed\end{lemma}
\begin{proof}
A finitely generated subgroup of a free/surface group is a free/surface group, 
and equality holds in that case. An infinitely generated subgroup of 
a free/surface group contains an infinitely generated free summand. 
Consequently, both $\pi_1$ and $H_1$ rank are infinite for such subgroups.
\end{proof}

\begin{lemma}\label{finitely generated}  An 
$H_1$-injective subgroup $G$ of a free/surface group $K$ is finitely 
generated.\end{lemma}

\begin{proof} Rank $\pi_1(G)=\rank(G/[G,G])\le \rank(K/[K,K])=\rank \pi_1(K)<\infty$.
\end{proof}

\begin{lemma} 
\label{core} A $1$-ended, orientable, irreducible 3-manifold $M$ 
with compact boundary is an end-manifold if and only if
$\pi_1(M)$ is a free/surface group, $H_2(M,\partial M)=0$ and $\partial M$ is 
$\pi_1$-injective.

Every closed embedded $\pi_1$-injective surface in an end-manifold is boundary parallel.
\end{lemma}

\begin{proof} 
Let $M$ be an end-manifold with core $C$ of the form
$\partial M\times I\cup 1$-handles or  handlebody if $\partial M=\emptyset$. 
Since the inclusion $C\to M$ is a homotopy equivalence, $\partial M $ 
is incompressible and $\pi_1(M)$ is a free/surface group.  If 
$T\subset M$ is a compact properly embedded
$\pi_1$-injective surface, then $T$ can be homotoped rel
$\partial T$ into $C$.  The cocores $D_i$ of the $1$-handles are properly
embedded disks whose boundary misses $\partial T$. Since $T$ is homotopically
essential, it follows that each intersection $T \cap D_i$ is 
homotopically inessential in $T$, and therefore $T$ can be homotoped 
off the cocores of the $1$--handles. Once this is done, $T$ can be 
further homotoped rel. boundary into $\partial M$, since $C$ 
deformation retracts to $\partial M$ in the complement of the
cocores of the $1$--handles. This implies that $H_2(M,\partial M)=0$.

If $H_2(M,\partial M)=0$ and $M$ has incompressible
boundary, a connected, closed orientable 
incompressible surface $R$ must separate off a connected,
compact Haken manifold $X$ with
incompressible boundary.  If $\pi_1(M)$ is also
a free/surface group, then $\pi_1(X)$ is a closed
orientable surface group and using \cite{St} we conclude that
$X=N(T)$ for some component $T$ of $\partial M$, so
$R$ is boundary parallel.  Therefore, if
$\partial M=\emptyset$, then any core is a handlebody.  If 
$\partial M \ne \emptyset$, then $M$ has a core $C$ which 
contains $\partial M$ \cite{Mc}.  If $C'$ is obtained
by maximally compressing $C$, then  each component of 
$\partial C'$ is boundary parallel and
hence $C=\partial M\times I\cup 1$-handles.\end{proof}

\begin{corollary}\label{submanifold 
end-manifold}  If $\mW$ is a 1-ended, $\pi_1$-injective
submanifold of the end-manifold $M$ such that
$\pi_1(\mW)$ is finitely generated 
and $\partial \mW$ is a union of components of $\partial M$,
then $\mW$ is an end-manifold.\end{corollary}

\begin{definition}   Given a 
connected compact subset $J$ of an open
irreducible 3-manifold $M$, the \emph {end reduction} $\mW_J$
of $J$ to $M$ is to first approximation the 
smallest open submanifold of $M$ which can engulf, up
to isotopy, any closed surface in $M\backslash J$ which is incompressible
in $M\backslash J$. End reductions were introduced by Brin-Thickstun 
\cite{BT1, BT2}.   Their basic properties were developed by 
Brin-Thickstun \cite{BT1, BT2} and Myers \cite{My}.   In particular 
\cite{BT1} show that $W_J$ can be created via the following 
procedure.  If $V_1\subset V_2\subset \cdots$
is an exhaustion of $M$ by compact connected codimension-0 submanifolds
such that $J\subset V_1$, then one inductively obtains an exhaustion 
$W_1\subset W_2\subset\cdots$ of $\mW_J$ by compact sets as follows.
Transform $V_1$ to $W_1$ through a maximal series of intermediate
manifolds $U_1=V_1, U_2, \cdots, U_n=W_1$ where $U_{k+1}$ is obtained from $U_k$ by  one 
of the following 3 operations.

\begin{enumerate}
\item{Compress along a disc disjoint from $J$.}
\item{Attach a 2-handle to $U_k$ which lies in $M\backslash \inte(U_k)$, and
whose attaching core circle is essential in $\partial U_k$.}
\item{Delete a component of $U_k$ disjoint from $J$.}
\end{enumerate}

Having constructed $W_i$, pass to a subsequence of the $V_j$'s and reorder so that
$W_i\subset \inte(V_{i+1})$.  Finally pass from $V_{i+1}$ to $W_{i+1}$ via a 
maximal sequence of the above operations.  Since $\partial W_i$ is 
incompressible in $M-J$, an essential compression of
$U_k$ can be isotoped rel boundary to one missing $W_i$.  Therefore,
we will assume such operations miss $W_i$ and 
hence $W_i\subset\inte(W_{i+1})$.  Brin and Thickstun \cite{BT1} show 
that $\mW_J$ is up to isotopy independent of all choices.

We  say that $\{W_i\}$ is a \emph{standard exhaustion} of $\mW_J$ if
$W_1\subset W_2\subset \cdots$ and $\mW_J=\cup_i W_i$, where for each $i$,
$W_i$ arises from $V_i$  via a sequence of the  three
end-reduction operations and $V_1\subset V_2\subset\cdots$ is an
exhaustion of $M$ by compact submanifolds.
\end{definition}

\begin{remark}
Note that operations 
(1) and (2) reduce the sum of the ranks of $\pi_1$ of the boundary
 components. It follows that the transition from $V_i$ to $W_i$ is obtained by 
a {\em finite} sequence of operations.
\end{remark}

\begin{remark} 
(Historical Note)  Brin and Thickstun \cite{BT1},
\cite{BT2} study end reductions to develop a necessary and
a sufficient condition, \emph{end 1-movability}, for taming an end of a
3-manifold.   More recently, Myers \cite{My}  has promoted the
use of end reductions to address both the $\BR^3$-covering space
conjecture and the Marden conjecture.\end{remark}

\begin{lemma}\label{injectivity}  The 
inclusion $i_J:\mW_J\to N$ induces
$\pi_1$ and $H_1$-injections.  The 
latter in both $\BZ$ and $\BZt$ homology.\end{lemma}

\begin{proof} 
The $\pi_1$-injectivity was first proven in \cite{BT2} and rediscovered in
\cite{My}.   Our proof of $H_1$-injectivity mimics 
the proof of $\pi_1$-injectivity in \cite{My}. Suppose
$C\subset W_i$ is a union of oriented simple
closed curves bounding the surface $S\subset M$. 
Note that by elementary
$3$--manifold topology, we can assume $S$ is 
embedded.

By choosing $n$ sufficiently large we can assume
that $W_i\cup S\subset V_n$.  If $V^1_n$ is obtained by adding a
2-handle to $V_n$, then $S\subset V^1_n$.
If $V^1_n$ is obtained by compressing $V_1$, via a compression missing $J$,
then by modifying $S$ near the compressing disc we obtain a surface $S_1$ spanning
$C$ (orientably, if need be) with $S_1\subset V^1_n$.  If $V^1_n$ is 
obtained by deleting components of $V_1$ which miss $C$, 
then $S_1=S\cap V^1_n$ still spans $C$.  Since $W_n$ is obtained from 
$V_n$ by a finite sequence of such operations it
follows that $C$ bounds in $W_n$ and hence in $\mW_J$. \end{proof}

$H_1$-injectivity of $\mW_J$ in $N$ gives us the following crucial 
corollary:

\begin{corollary}  \label{end-reduction is finitely 
generated}  An end-reduction in an end-manifold has finitely 
generated fundamental group.\end{corollary}

\begin{proof} Combine Lemma~\ref{injectivity} with Lemma~\ref{finitely generated}.\end{proof}

\begin{definition}  If $\mW_J$ is an end reduction of the codimension-0
submanifold $J$ in $N$, then we 
say that $\mW_J$ is \emph{trivial} if $\mW_J$ is isotopic to an open 
regular neighborhood of $J$ or equivalently $\mW_J$
is isotopic to $\inte(J)$.  $\mW_J$ is \emph{eventually trivial} if it has 
an exhaustion $W_1\subset W_2\subset
\cdots$ such that $\partial W_i$ 
is parallel to $\partial W_j$ for all $i,j$.
\end{definition}

We now study end reductions of disconnected spaces $J$.  While the
following technology and definitions can be given for more general objects we 
restrict our attention to a finite unions of pairwise disjoint
closed  (possibly non simple) curves none of which lie in a 3-cell.
Ultimately we will address end reductions of infinite sequences of exiting curves.

\begin{definition} \label{house} If $J$ 
is a finite union of pairwise disjoint closed curves in
an open irreducible 3-manifold $M$, we say  that $J$ is 
\emph{end nonseparable} if there is a compact connected submanifold 
$H$ such that $J\subset \inte(H)$ and $\partial H$ is incompressible 
in $M\backslash J$.  Such an $H$ is called a
\emph{house} of $J$.  If $J$ is
end nonseparable, then define $\mW_J$ to be an end-reduction of 
$H$, and call $\mW_J$ the {\em end reduction of $J$}.
\end{definition}

\begin{lemma}  \label{well defined}  The end-reduction 
$\mW_J$ of an end nonseparable union $J$ of  closed curves is 
well defined up to isotopy.\end{lemma}

\begin{proof} Let $H$ and 
$H^\prime$ be two houses for $J$. We want to show
that if $\mW_H$ is an end reduction of $H$, then there is an isotopy
of $H'$ to $H'_1$ fixing $J$, so that the end reduction $\mW_{H'_1}$ of $H'_1$
is equal to $\mW_H$. By the definition of a house for $J$, both $H$ and
$H'_1$ satisfy the property that they are connected submanifolds of 
$M$ whose boundaries are incompressible in $M\backslash J$.

Let $\{W_i\}$ be a standard exhaustion of $\mW_H$ arising from 
the exhaustion $\{V_i\}$ of $M$.  By passing to a subsequence we can 
assume that  $H^\prime\cup H\subset V_2$.   By considering the 
passage of $V_2$ to $W_2$, we observe that $H^\prime$ can
be isotoped to $H^\prime_1$ rel $J$ to lie in $\inte(W_2)$ and that
$\partial W_2$ is incompressible in $M\backslash H^\prime_1$.  Thus
$\mW_H$ is also an end-reduction of $H'$.  Since
end reductions are unique up to isotopy the result follows, and
we may unambiguously denote $\mW_H$ by $\mW_J$.\end{proof}

\begin{lemma}\label{partition}  Let $A$ be a 
finite union of pairwise disjoint  closed
curves in the open irreducible 3-manifold $M$. Then $A$ canonically
decomposes into finitely many maximal pairwise
disjoint end nonseparable 
subsets $A_1,\cdots,A_n$.   Indeed, if $B$ is a
maximal end nonseparable
subset of $A$, then $B=A_i$ for some $i$. 
\end{lemma}

\begin{proof}  Since each element of $A$ is end nonseparable, it  suffices
to show  that if $B$ and $C$ 
are end nonseparable subsets  of $A$, then either
$C\cup B$ is end nonseparable or $C\cap B=\emptyset$.  Let $H_B$ and $H_C$
be houses for $ B$ and $C$ respectively.  Let $V\subset N$ be a compact 
submanifold containing $H_B \cup H_C$.  By considering the passage of 
$V$ to $W$ by a maximal sequence of compressions, 2-handle
additions, and deletions which are taken with respect to $B\cup
C$, one sees that $H_B$ (resp. $H_C$)
can be isotoped to lie in $W$ via an isotopy 
fixing $B$ (resp. $ C$). If $B\cap C\neq\emptyset$,
then $W$ is connected and hence is a house for $B\cup C$.
\end{proof}

\begin{lemma} \label{disjointness} If $A_1,\cdots 
A_n$ are the maximal end nonseparable components of
a finite set $A$ of pairwise disjoint  closed curves in an open
irreducible 3-manifold $M$, then they have pairwise disjoint end reductions.  In 
particular they have pairwise disjoint houses. 
\end{lemma}

\begin{proof} Let $A_1, A_2, \cdots ,A_n$ be the 
maximal end nonseparable
subsets of $A$.  Let $\{V_k\}$ be an exhaustion 
of $M$ with $A \subset V_1$.  Consider 
a sequence $V_1=U_1,\cdots,U_n=W_1$ where the passage from one to the 
next is isotopy, compression, 2-handle addition or deletion, where 
the compressions or deletions are taken with respect
to $A$. By passing to a subsequence of the exhaustion we can assume that
$W_1\subset V_2$, and in the above manner pass from $V_2$ to 
$W_2$.  In like manner construct $W_3,W_4,\cdots$.   By deleting 
finitely many of the first $W_i$'s from the sequence
and reindexing, we can assume that all the $W_i$'s have the same number of 
components.

It suffices to show that if $W$ is a component of $W_k$, 
then $W$ contains a unique $A_i$ and that $\partial W$
is incompressible in $M\backslash A_i$.  Indeed, it suffices to 
prove incompressibility of $\partial W$ in $M\backslash (W\cap A)$, 
for then $W$ is a house and can only contain one $A_i$ by maximality. 
If $\partial W$ is compressible in $M\backslash (W\cap A)$ it
must compress to the outside via some compressing disc $D$.  Consider 
a term $V_n$ in the exhausting sequence with $W_k\cup D\subset V_n$. 
By considering the passage of $V_n$ to $W_n$ we can rechoose the disc 
spanning $\partial D$ to obtain a new compressing
disc $E\subset W_n$.   Since $\partial W$ is incompressible in
$M\backslash \inte(W_k)$, it follows that $E$ must hit a component of $W_k$ 
distinct from $W$.  This imples that  $W_n$ contains
fewer components than $W_k$, which is a contradiction.
\end{proof}

\begin{lemma}  \label{multi injectivity}If $A_1, A_2, \cdots ,A_n$
are as in Lemma~\ref{disjointness}, with pairwise 
disjoint end  reductions $\mW_{A_1},\mW_{A_2} ,\cdots , \mW_{A_n}$, 
then $\mW_{A_1}\cup \mW_{A_2} \cup\cdots \cup \mW_{A_n}$ is
$H_1$-injective in $M$, in both $\BZ$ and $\BZt$ coefficients.\end{lemma}

\begin{proof}  Repeat the proof of Lemma 
\ref{injectivity}.\end{proof}

\begin{corollary} \label{maximality bounded} Let $A$ be a union of
finitely many pairwise  disjoint   closed
curves in the end-manifold $M$.  If each component of $A$ is
homotopically nontrivial,
then $A$ breaks up into at most rank$(\pi_1(M))$ maximal non separable
subsets.\end{corollary}

\begin{proof}  If $A$ partitions into maximal non separable subsets
$A_1,\cdots, A_n$, then the
$H_1$-rank of $\mW_{A_i}$ is non trivial, since $\pi_1(\mW_{A_i})$ is a
nontrivial subgroup of a free/surface group.  
Now apply the previous lemma.\end{proof}

\begin{lemma}\label{infinite partition} 
Let $\gamma_1,\gamma_2,\cdots$ be a sequence of homotopically 
nontrivial, pairwise disjoint closed curves in the 
end-manifold $M$.  Then we can group together finitely many
of the curves into $\gamma_1$, and pass to a subsequence so 
that
\begin{enumerate}
\item  Any finite subset of $\{\gamma_1,\gamma_2,\cdots\}$ which contains
$\gamma_1$ is end nonseparable.
\item  Each component of $ \gamma_1$, and each $\gamma_i$, $i\ge 2$
represent the same element of $H_1(M, \BZt).$\end{enumerate}\end{lemma}

\begin{proof} By passing to a subsequence we can assume that each
$\gamma_i$ represents the same 
element of $H_1(M,\BZt)$.  By Lemma
$\ref{partition}$, if $T$ is a 
finite subset of $\Gamma$, then
$T$ canonically  partitions into 
finitely many end nonseparable
subsets $S_1,\cdots, S_n$
with corresponding pairwise disjoint end reductions $\mW_1,\cdots, \mW_n$. 
Define 
$$C(T)=\sum_{i=1}^n \rank(H_1(\mW_i,\BZt))= \rank(H_1(\cup_{i=1}^n \mW_i,\BZt))\le
\rank  (H_1(M,\BZt))$$ 
where the last inequality follows from Lemma~\ref{multi injectivity}. 
Define 
$$C(\Gamma)=\max\{C(T)\, |\, T \ \textrm{is a finite subset of}\ \Gamma\}$$  
Now pass to an infinite subset of $\Gamma$ with
$C(\Gamma)$  minimal.  By Lemma~\ref{multi injectivity},  if 
$T\subset \Gamma$ with $C(T)=C(\Gamma)$, then adding a new element to 
the $T$ does not increase the number of end nonseparable subsets in 
its canonical partition.  Since $C(\Gamma)$ is minimal, we 
can enlarge $T$ by adding finitely many elements so that the enlarged
$T$, which by abuse of notation we still call $T$, is end 
nonseparable.  Again by maximality of $C(T),\ T$
together with any finite subset of $\Gamma$ is still end non separable.  Now 
express $\Gamma$ as $\cup \gamma_i$ with $\gamma_1=T$. 
\end{proof}

\begin{theorem}[Infinite end engulfing theorem]\label{refined}  If
$\gamma_1, \gamma_2,\cdots$ is a locally 
finite sequence of pairwise disjoint, homotopically non trivial, 
closed curves in the end-manifold $M$, then after passing 
to a subsequence,  allowing $\gamma_1$ to have multiple components and
fixing a base point  $x\in \gamma_1$, there exist compact 
submanifolds $D\subset W_1\subset W_2\subset\cdots$ of $M$ such 
that 
\begin{enumerate}
\item $\partial W_i\cap\partial M$ is a union of components of
$\partial  M$ and $\partial W_i-\partial M$ is
connected.
\item If $\Gamma_i=\cup_{j=1}^i\gamma_j$, then $\Gamma_i\subset
W_i$,  and $\Gamma_i$ can be homotoped into $D$ via a homotopy
supported in $W_i$.
\item $\partial W_i$ is 2-incompressible rel $\Gamma_i$.
\item If $\mW=\cup W_i$, then $\mW$ is $\pi_1$ and
$H_1$--injective in both $\BZ$ and $\BZt$ coeficients.
\item $D$ is a core of $\mW$ and is of the form 
$\partial \mW\times I \, \cup$ 1-handles.
\end{enumerate}
\end{theorem}

\begin{figure}[ht]
\centerline{\relabelbox\small 
\epsfxsize 3.0truein
\epsfbox{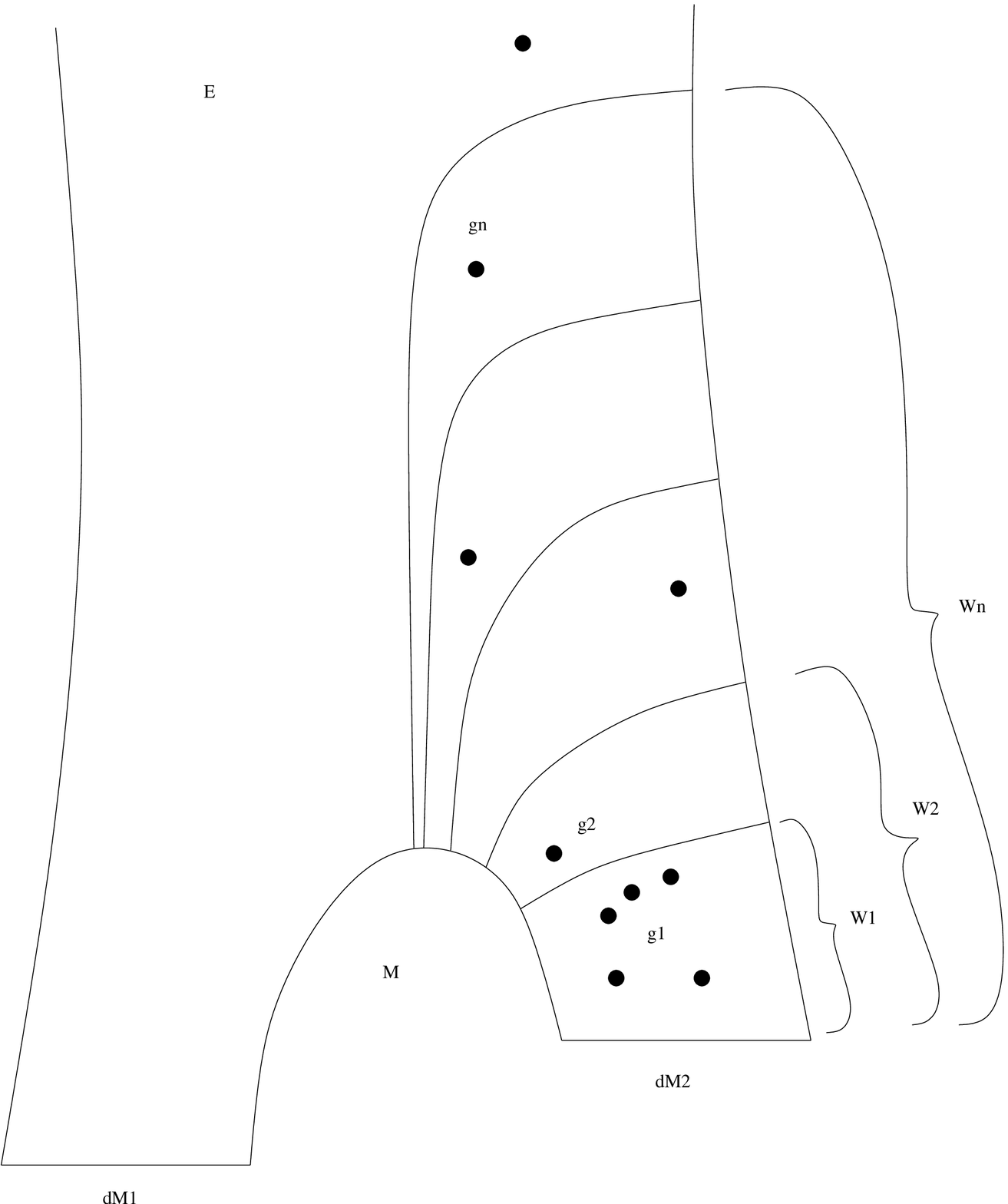}
\relabel {M}{$M$}
\relabel {dM1}{$\partial M$}
\relabel {dM2}{$\partial M$}
\relabel {E}{$\mE$}
\relabel {g1}{$\gamma_1$}
\relabel {g2}{$\gamma_2$}
\relabel {gn}{$\gamma_n$}
\relabel {W1}{$W_1$}
\relabel {W2}{$W_2$}
\relabel {Wn}{$W_n$}
\endrelabelbox}
\caption{A schematic view of $\lbrace W_i \rbrace$, $\lbrace \gamma_i \rbrace$ and $M$.}
\label{51}
\end{figure}

The conclusion of this theorem is schematically depicted in Figure~\ref{51}.

\begin{proof}    By passing to a subsequence and allowing  $\gamma_1$ to
have multiple components we can assume that $\Gamma=\{\gamma_i\}$
satisfies the conclusions of Lemma \ref{infinite partition}.  Assume that
$\gamma_1$ has at least two elements. 

Let $\mW_1$ be an end reduction to $\gamma_1$ with standard exhaustion
$W_{1,1}\subset W_{1,2}\subset \cdots$.  Let $D_1$ be a core for $\mW_1$ with
$x\in D_1$. Since $\mW_1$ is $\pi_1$-injective in $M$, by passing to a
subsequence we can choose $W_{1,1}$ so that each component of
$\gamma_1$ can be homotoped into $D$ via a homotopy in $W_{1,1}$.  Furthermore,
$H_1$-injectivity allows us to assume that within $W_{1,1}$ each 
$\gamma\in \gamma_1$ is
$\BZ_2$-homologous to a component $\gamma'$ of  $\gamma_1$ with
$\gamma'\neq\gamma$.  If
$\gamma$ represents the trivial class, it should be homologically trivial in
$W_{1,1}$.  Finally $W_{1,1}$ should be sufficiently large so that $\partial
W_{1,1}$ is a union of components of $\partial M$ and a single component
disjoint from $\partial M$.  Note that $\partial W_{1,1}$ is 2-incompressible
rel $\gamma_1$.  Indeed, by construction $\partial W_{1,1}$ is
incompressible in $M/\gamma_1$, hence $\partial W_{1,1}$ is incompressible to the
outside and any essential compressing disc $D$ for $W_{1,1}$ must intersect
$\gamma_1$ at least once.  If $D$ meets the component $\gamma$ of $\gamma_1$,
then since $\gamma$ is either $\BZ_2$-homologically trivial in $W_{1,1}$ or
$\BZ_2$-homologous to a $\gamma'\neq\gamma$ it follows that $|D\cap\gamma_1|\ge 2$.  

By passing to a subsequence of $\{\gamma_i\}$ we can assume that 
$\gamma_2\cap W_{1,1}=\emptyset$.   By Lemma \ref{infinite partition},
$\Gamma_2=\{\gamma_1,\gamma_2\}$ is end nonseparable.  Let $\mW_2$ be an end
reduction for $\Gamma_2$ with standard exhaustion $W_{2,2}\subset
W_{2,3}\subset W_{2,4}\subset\cdots$ where $W_{1,1}\subset \inte(W_{2,2})$. 
As above let $D_2$ be a core for $\mW_2$ with $x\in D_2$ and by choosing 
$W_{2,2}$ sufficiently large we can assume that it supports a homotopy of
$\Gamma_2$ into $D_2$ as well as homologies between elements of $\Gamma_2$. 
Finally, $\partial \mW_{2,2}$ is a union of components of $\partial M$ and a single
component disjoint from $\partial M$.  As above, $\partial \mW_{2,2}$ is
2-incompressible rel $\Gamma_2$.  

Having inductively constructed $\Gamma_{i-1}, W_{i-1,i-1}$ and $D_{i-1}$ pass to
a subsequence of $\{\gamma_j\}$ so that
$\Gamma_{i-1}\subset\Gamma_i=\{\gamma_1,\cdots,\gamma_i\}$ with 
$\gamma_i\cap W_{i-1,i-1}=\emptyset$.  Let $\mW_i$ be an end reduction
of $\Gamma_i$ with standard exhaustion 
$W_{i,i}\subset W_{i, i+1}\subset W_{i,i+2}\subset\cdots$ where 
$W_{i-1,i-1}\subset \inte(W_{i,i})$.  Let $D_i$ be
a core of $\mW_i$, but if possible let $D_i= D_{i-1}$.  Finally  $W_{i,i}$ 
should be chosen sufficiently large to support homotopies of $\Gamma_i$ into
$D_i$ and homologies as described in the previous paragraphs and so that
$\partial W_{i,i}$ is a union of components of $\partial M $ and a single other
component.  As above $\partial W_{i,i}$ is 2-incompressible rel $\Gamma_i$.

Let $\mW=\cup_i W_{i,i}$.  The proof of Lemma \ref{injectivity} shows that $\mW$
is $\pi_1$ and $H_1$-injective in M in $\BZ$ and $\BZ_2$ coeficients.  
By Lemma~\ref{finitely generated}
$\mW$ has finitely generated fundamental group and hence for some $n$,
the inclusion $W_{n,n}\to \mW$  is $\pi_1$-surjective.   Since 
$D_n$ and $ W_{n,n}$ have the same $\pi_1$-image in $\pi_1(\mW_n)$ and
hence in $\pi_1(M)$, they have the same image in $\pi_1(\mW)$ and hence $D_n\to\mW$
is $\pi_1$-surjective.  Since the inclusions $D_n\to \mW_n\to M$ are
$\pi_1$-injective it follows that $D_n$ is also $\pi_1$-injective in $\mW$ and
hence is a core of $\mW$.  Therefore, if $m\ge n$, then $D_n$ and $W_{m,m}$ have
the same $\pi_1$-image in $\pi_1(\mW)$ and hence in $\pi_1(M)$.  Since the $\mW_m$ is
$\pi_1$-injective in $M$, it follows that $D_m$ and $D_n$ have the same image in
$\pi_1(\mW_m)$ and hence $D_n$ is a core of $\mW_m$ for all $m\ge n$.  

To complete the proof of the Theorem reorganize $\{\gamma_i\} $ so that
$\gamma_1$ is now the old $\{\gamma_1,\gamma_2,\cdots,\gamma_n\}$, and
for all $i\in \BN,\ \gamma_{1+i}=$ old $\gamma_{n+i}$. Let $D=D_n$. Finally, for
$i\in \BN$, let $W_i=W_{n+i-1,n+i-1}$.   

By Corollary \ref{core} each $\mW_i$ is an end-manifold and hence $D$
could have been taken to be of the form $N(\partial \mW)\times I \cup 1$-handles, 
if $\partial \mW\neq\emptyset$ and a handlebody otherwise.\end{proof}

\begin{definition}  Call the $\mW$ constructed in Theorem~\ref{refined}  
an \emph{end-engulfing} of $\Gamma$.\end{definition}

The material in the rest of this chapter will not be used elsewhere
in this paper; in particular, it is not used to prove any of the results
of \S 0.  

\begin{lemma}\label{engulfing}  If $J\subset J'$ are finite, end
nonseparable unions of homotopically essential, pairwise disjoint, 
closed curves with end-reductions $\mW$ and $\mW'$, then $\mW$ is
isotopic rel $J$ to $\mW_1$, where
$\mW_1\subset \mW'.$\end{lemma}

\begin{proof} Let $W_1\subset W_2\subset\cdots$ be a standard  
exhaustion of $\mW$.  Let $Z_1\subset Z_2 \subset
\cdots$ be a standard exhaustion of
$\mW'$  arising from the  exhaustion
$\{V_i\}$ of $M$.  By passing to subsequence we can assume that
$W_1\subset V_1$.  By considering the passage of $V_1$ to $Z_1$ we
can isotope $W_1$ rel $J$
to lie in $Z_1$.  Proceeding by induction and passing to subsequence,
we can assume that $W_k\subset V_k$ and
$W_{k-1}\subset Z_{k-1}$.  By considering the passage of $V_k$ to $Z_k$
(which fixes $Z_{k-1}$) we can isotope $W_k$ rel $W_{k-1}$ to lie
in $Z_k$.  The isotoped $W_i$'s give rise to an isotopy of $\mW$ to
$\mW_1$ with $\mW_1\subset \mW'$. 
\end{proof}

\begin{remark} Given $J\subset J'$ with end-reductions $\mW$ and
$\mW'$ one can isotope $\mW \rel J$  to $\mW_1$ so that
$\mW_1\subset \mW'$ (Lemma~\ref{engulfing}). On the other hand one cannot in
general isotope $\mW'$ to contain $\mW$.  One need only look at the case 
of $J\subset J'$ being nested balls in the Whitehead manifold to
find examples. Such considerations make it challenging to find
nested end--reductions $\mW_1\subset\mW_2\subset \mW_3 \subset\cdots\ $. 
\end{remark}

\begin{theorem}[Finite end reduction theorem]\label{finite non separable} 
Let $M$ be an end-manifold. 
If $\Gamma=\{\gamma_i\}$ is an end nonseparable union of finitely 
many homotopically essential, pairwise disjoint,  closed curves, then 
an end reduction $\mW_\Gamma$ of $\Gamma$ has
finitely generated fundamental group and given a standard exhaustion
$\{W_i\}$, by passing to a subsequence, for all 
$i,j<k,$
$$\text{in}_*(\pi_1(W_i))=\text{in}_*(\pi_1(W_j)) \subset \pi_1(W_k)$$ 
and the map
$\text{in}_*:\pi_1(W_k)\to \pi_1(\mW_\Gamma)$ restricted to 
$\text{in}_*(\pi_1(W_i))$ induces an isomorphism
onto $\pi_1(\mW_\Gamma)$. 
Here in$_*$ denotes the map induced by inclusion.
\end{theorem}

We first prove a topological lemma.

\begin{lemma}\label{scott}  If $M$ 
is an end-manifold, then $M$ has an exhaustion by
compact manifolds $V_1\subset V_2\subset \cdots$, such that for each
$i>1$ either $V_i$ is a handlebody, in which case
$\partial M=\emptyset$, or $V_{i}$ is 
obtained by attaching 1-handles to a $N(\partial M)$.
\end{lemma}

\begin{proof} If $\partial M=\emptyset$, then 
$\pi_1(M)$ is free and this result follows
directly from \cite{FF}. 
If $\partial M\neq\emptyset$, it suffices
to show that if $X$ is any 
compact submanifold of $M$, then $X\subset V$ where
$V$ is obtained by thickening $\partial M$ and attaching
1-handles.  We use the standard argument, e.g. see \cite{BF},
\cite{BT2} or \cite{FF}.
Using the loop theorem we can pass from $X$ to a submanifold $Y$, 
with incompressible boundary  via a  sequence of
compressions and external 2-handle additions.  By appropriately
enlarging $X$ to $X_1$,
so as to contain these 2-handles, we can pass from $X_1$ to 
$Y$ by only compressions. By enlarging $Y$, and hence $X_1$, we 
can assume that $\partial M\subset Y$ and no component
of $M\backslash \inte(Y)$ is compact.  By Lemma~\ref{core} each 
component of $\partial Y$ is boundary parallel
and hence $Y$ is of the form $N(\partial M)\cup 1$-handles.\end{proof}

\noindent\textit{Proof of Theorem~\ref{finite non separable}}.    Let $V_1\subset 
V_2\subset\cdots$ be an  exhaustion of
$M$ as in Lemma~\ref{scott} so that $\Gamma\subset V_1$.   
Let $W_1\subset W_2\subset\cdots$ be a standard exhaustion of 
$\mW_\Gamma$ arising from the exhaustion $\{V_i\}$ of $M$.

By Definition~\ref{house} and 
Lemma~\ref{end-reduction is finitely generated},
$\pi_1(\mW_\Gamma)$ is finitely
generated, so we can pass to a subsequence and assume
that the induced map $\pi_1(W_1)\to \pi_1(\mW_\Gamma)$ is 
surjective.

Let $H_i=\text{in}_*(\pi_1(W_1))$ where $\text{in}:W_1\to W_i$ is inclusion.  We
now show that after passing to a subsequence 
of the $W_i's, i\ge 2,$  the induced maps 
$$H_2\to H_3\to \cdots\to \pi_1(\mW_\Gamma)$$
are all isomorphisms.

For $j\ge 1$, let $G_j=\alpha^j_*(\pi_1(W_1))$, where $\alpha^j:W_1\to V_j$
is inclusion.  Each $G_j$ is a finitely generated subgroup of 
$\pi_1(V_j)$ and hence is a  free product of finitely many closed 
orientable surface groups and a finitely generated free group.  Since 
for all $j$, rank$(G_j)\le \rank(G_1)$,
there are only finitely many possibilities for such groups and
hence by passing to a subsequence  we can assume that for $j,k >1$, the
groups  $G_j$ and $G_k$ are abstractly isomorphic. Free/surface groups are
obviously linear, hence residually finite by Malcev \cite{Mv}. Further,
Malcev \cite{Mv} went on to show that finitely generated residually
finite groups are Hopfian, i.e. surjective self maps are isomorphisms.
This implies that the induced maps $G_{2}\to G_{3}\to G_{4}\to \cdots$ are 
all isomorphisms.   If 
$$K=\kernal (\text{in}_*:\pi_1(W_1)\to \pi_1(\mW_\Gamma))$$ 
then
$$K=\kernal(\pi_1(W_1)\to\pi_1(M))=\kernal(\pi_1(W_1)\to G_{2})$$
We now show that
$K=\kernal(\pi_1(W_1)\to \pi_1(W_2))$.  One 
readily checks that if $W_1\subset V$ 
and $K=\kernal(\pi_1(W_1)\to \pi_1(V))$, and $V^\prime$ is obtained 
from $V$ by compression, 2-handle addition or deletion where these 
operations are performed in the complement of $W_1$, then
$K\subset \kernal(\pi_1(W_1)\to \pi_1(V^\prime))$.  This
implies that if $K_2=\kernal(\pi_1(W_1)\to \pi_1(W_2))$, then
$K\subset K_2$.  On the other hand $K_2\subset K$ since
$K=\kernal(\pi_1(W_1) \to \pi_1(\mW_\Gamma))$.  Therefore, the
induced maps $H_2\to H_3\to \pi_1(\mW_\Gamma)$ are isomorphisms.

Apply the argument of the previous paragraph to obtain a subsequence
of $\lbrace W_i \rbrace$
which starts with $W_1$ and $W_2$ such that 
the $\pi_1$-image of $W_2$ in $W_j, j>2$,  maps isomorphically 
to $\pi_1(W_\Gamma)$, via the map induced by inclusion.
Continue in this manner to construct $W_3, W_4, \cdots$.
\qed\newline

\noindent\textbf{Addendum to Theorem \ref{refined}} 
\emph{We can obtain the following additional property. 
If} $i,j<k$, \emph{then}
$$\text{in}_*(\pi_1(W_i))=\text{in}_*(\pi_1(W_j))
\subset \pi_1(W_k)$$ \emph{where in denotes inclusion.   The 
map} $\text{in}_*:\pi_1(W_k)\to \pi_1(\mathcal{W})$ \emph{restricted to} 
$\text{in}_*(\pi_1(W_i))$
\emph{induces an isomorphism onto}
$\pi_1(\mathcal{W})$.\vskip 12 pt

\begin{proof}  Apply Theorem~\ref{refined} to produce the space
$D$ as well as the sets
$\{\gamma_i\}, \{\Gamma_i\},$ $\{W_i\}$ which we now relabel as 
$\{\gamma'_i\}, \{\Gamma'_i\}, \{W'_i\}$.  Define 
$\gamma_1= \Gamma_1=\gamma'_1$,
$W_1=W'_1$, $\gamma_2=\gamma_2'$ 
and
$\Gamma_2=\Gamma_2'$.  Let
$W'_2=W^1_2\subset W^2_2\subset W^3_2\subset \cdots$ be a standard
exhaustion of an end reduction 
$\mW_2$ of $\Gamma_2$, which we can assume
satisfies  the conclusions of Theorem~\ref{finite non separable}.  
Defining $W_2 = W^2_2$, we  see that the
restriction of $\text{in}_*:\pi_1(W_2)\to \pi_1(\mW_2)=\pi_1(D)$ to
in$_*(\pi_1(W_1))\subset \pi_1(W_2)$ is an 
isomorphism.  Choose $\gamma_3=\gamma'_{i_3}\in \{\gamma'_i\}$ so 
that $\gamma_3\cap W_2=\emptyset$ and define $\Gamma_3=\Gamma_2\cup \gamma_3$.  Let
$W^1_3\subset W^2_3\subset  W^3_3\subset\cdots$ be a 
standard exhaustion of an end-reduction $\mW_3$ of $\Gamma_3$
which satisfies the conclusions of Theorem~\ref{finite non separable} 
and has $W_2\subset W^1_3$.  Defining $W_3 = W^2_3$, we see that
the restriction of $\text{in}_*:\pi_1(W_3)\to \pi_1(\mW_3)=\pi_1(D)$ to 
in$_*(\pi_1(W_2))\subset\pi_1(W_3)$ is an isomorphism.  Now define 
$\gamma_3=\gamma_{i_3}'$ and $\Gamma_2=\Gamma_1\cup \gamma_{i_3}$. 
In a similar manner construct 
$\gamma_4,\gamma_5\cdots\ ,\  \Gamma_4, \Gamma_5,\cdots\ ,\  W_4, W_5,\cdots$ 
and finally define $\mW=\cup W_i$.
\end{proof}

\begin{remarks}  If one allows each $\gamma_i$ to be a finite set of
elements, then we can obtain the conclusion (in Theorem~\ref{refined} and its
addendum) that each $\gamma_i$ is $\BZt$-homologically trivial. 
\end{remarks}

\begin{question}  Let $M$ be a connected, compact, 
orientable, irreducible 3-manifold such that $\chi(M)\neq 0$ and let 
$G$ be a subgroup of $\pi_1(M)$.  If the induced map
$G/[G,G]\to H_1(M)$ is injective, is $G$ finitely generated?\end{question}

\begin{question}  Let $\Gamma$ be a locally 
finite collection of pairwise disjoint
homotopically essential closed curves such that
$C(\Gamma)=C(\Gamma')$ for any infinite 
subset $\Gamma'$ of $\Gamma$.  Is it true, that given
$n\in \BN$, there exists an end-engulfing of $\mW=\cup W_i$ of $\Gamma$
such that for all $i$,
$|E\cap \Gamma_i|\ge n$ for all essential compressing discs $E$ of $W_i$?
\end{question}

\section{Proof of Theorems \ref{main}, 	\ref{marden} and \ref{clean}:
Parabolic Free Case}

\noindent\textbf{Proof of Theorem \ref{main}.} By Lemma~\ref{end-manifold} and
Remark \ref{end-manifold remark} it suffices to consider the case 
that $\mE$ is the end of an end-manifold 
$M\subset N$ such that the inclusion $M\to N$ is a
homotopy equivalence.  By Lemma \ref{bonahon} there exists an
$\eta$-separated collection $\Delta=\{\delta_i\}$ of closed geodesics 
which exit $\mE$.  We let $\Delta_i$ denote the union 
$\Delta_i = \cup_{j \le i} \delta_j$.
 Apply Theorem~\ref{refined} to $\Delta$ and $M$ to 
pass to a subsequence also called
$\Delta$  where we allow $\delta_1$ 
to have finitely many components.  Theorem~\ref{refined} also 
produces  a manifold
$\mW$ open in $M$ and exhausted by compact 
manifolds $\{W_i\}$ having the following
properties.

\begin{enumerate}

\item $\mW$ is $\pi_1$ and $H_1$ injective (in $\BZ$ and $\BZt$ coeficients)
in $M $ and hence $\pi_1(\mW)$ is a free/surface group.

\item   $\partial W_i\backslash \partial M$ is a closed connected surface which
separates $\Delta_i$ from $\mE$ and is 2-incompressible in $N$ rel. $\Delta_i$.

\item There exists a compact submanifold core $D\subset W_1$ of $\mW$ such that
for each $i,\ \delta_i$ can be homotoped into $D$ via a
homotopy supported in $W_i$.   $D$ is either of the form
$\partial \mW\times I$
with 1-handles attached to the 1-side, if $\partial \mW\neq\emptyset$ or a
handlebody, otherwise.
\end{enumerate}

Let $G_i$ denote in$_*(\pi_1(D))\subset \pi_1(W_i)$.  Let $\partial_e D$ denote
$D\cap\partial M=\partial \mW$.  Let $ X_i$ denote the covering space
of $W_i$ with group $G_i$ and let $\hat D$ denote the lift of $D$.   Pick a homotopy of
$\Delta_i$ into $D$ supported in $W_i$.  This homotopy lifts to a
homotopy of $\hat \Delta_i$ into
$\hat D$, thereby picking out the closed preimages $\hat \Delta_i$ of
$\Delta_i$ which are in 1-1 correspondence with $\Delta_i$.  Let
$\{\hat\delta_1,\cdots,\hat\delta_i\}$
denote these elements.

\vskip 12pt
\noindent\emph{Claim.}    Each $W_i$ is a compact atoroidal Haken
manifold and $\partial W_i$ contains a surface of genus $\ge 2$.  
\vskip 12pt

\noindent\emph{Proof of Claim.}  Each  embedded torus in $N$ is compressible, since $N$ is
parabolic free.  A compressible torus in an irreducible 3-manifold is either a
\emph{tube}, i.e. bounds a solid torus or a \emph{convolutube}, i.e. bounds a
cube with knotted hole $X$, which is a 3-ball with an open regular neighborhood
of a properly embedded arc removed.  Furthermore,  $X$ lies in a 3-ball.  Therefore, if some
component of $\partial W_i$ is a torus, $W_i$ is either a solid torus or a cube with knotted
hole.  The former can contain at most one closed geodesic and the latter none.  Since $W_i$
contains at least two closed geodesics $\partial W_i$ cannot contain a torus.

If $W_i$ contained  an embedded incompressible torus $T$, then the compact region bounded by
$T$ would lie in $W_i$.  This implies that $T$ is a convolutube.  In $\BH^3$, the universal
covering of $N$,  let $\tilde\Delta_i$ denote the preimage of
$\Delta_i$ and $\tilde W_i$ the preimage of $W_i$.  Since $T$ lies in a 3-cell in $N$, $T$
lifts to a torus $\tilde T$ isometric to $T$.  Using the loop theorem, it follows that
$\tilde\partial W_i$ is incompressible in $\BH^3\backslash\tilde \Delta_i$ and $\tilde T$ is
incompressible in $\tilde W_i$.  We will show that  after an isotopy of $\tilde T$ supported in
$\tilde W_i$, there exists an embedded 3-ball $F\subset \BH^3$ such that $T\subset F$ and
$F\cap \tilde \Delta_i=\emptyset$.  This implies that $T$ is compressible in $F$, via a
compressing disc $D$ disjoint from $\tilde \Delta_i$.  Since $\tilde W_i$ is incompressible in
$\BH^3\backslash\tilde \Delta_i$ it follows that $D$ can be isotoped rel $\partial D$ so that
$D\subset \tilde W_i$.  This contradicts the fact that $\tilde T$ is incompressible in $\tilde W_i$.  

Here is how to find $F$.   Let $E\subset \BH^3$ be a large round ball transverse to $\tilde
\Delta_i$ which contains $\tilde T$ in its interior.  $\tilde \Delta_i\cap
E=\{\alpha_1,\cdots,\alpha_n\}$ is a finite union of unknotted arcs, i.e. there exists pairwise
disjoint embedded discs $\{D_1,\cdots,D_n\}$ such that for each $k, D_k\subset E$ and
$\partial D_k$ consists of $\alpha_k$ together with an arc lying in $\partial E$.  For each
$k$, either $\alpha_k\cap \tilde W_i=\emptyset$ or $\alpha_k\subset \tilde W_i$.  Since
$\partial \tilde W_i$ is incompressible in $\BH^3\backslash\tilde\Delta_i$, $E$ and the $D_k$'s
can be isotoped, via an isotopy which fixes $\tilde \Delta_i$ pointwise, to $E'$ and $D'_k$'s
so that $E'$ is a 3-ball containing $\tilde T$, 
the $D'_k$'s are unknotting discs for the $\alpha_k$'s
and for each $k$, either $D'_k\cap \tilde W_i=\emptyset$ or $D'_k\subset \tilde W_i$.  After an
isotopy  of $\tilde T$ supported in $E'\cap \tilde W_i$, for each $k$, $\tilde T\cap
D'_k=\emptyset$.  Finally let $F$ equal $E'-\cup_{j=1}^n \inte(N(D'_k))$, where $N(D'_k)$ is a
very small regular neighborhood of $D'_k$. \qed

\vskip 12pt

It now follows from Thurston (see Proposition 3.2 \cite{Ca} or \cite{Mo}) 
that $\inte(X_i)$ is topologically tame.  Let $\bar X_i$
denote its manifold compactification.  Since $\partial_e D$ is a union of  components
of $\partial W_i$ and $X_i$ is the $\pi_1(D)$ cover, there is a canonical identification of 
$\partial_eD$ with some set of components of $\partial \bar X_i$. 
Let $\partial_e\hat D$ denote these components.
Having the same homotopy type as $D$, it follows by the usual group theoretic
reasons that $\bar X_i$ either compresses down to a $3$-ball, or to a 
possibly disconnected  (closed orientable surface)$\times I$.  In 
the former case $\bar X_i$ is a handlebody which, for
reasons of Euler characteristic,  is of the same genus as $D$.   In 
the latter case, since $\partial_e\hat D$ is an incompressible surface,  
$\bar X_i$ is topologically $\partial_e\hat D\times I$ with 1-handles 
attached to the 1-side. Let $\bar S_i$ denote $\partial \bar X_i-\partial_e\hat D$.
Again by reason of Euler characteristic, 
$\genus(\bar S_i)=\genus(\partial_{\mE}D)$, where
$\partial_{\mE}D=\partial D-\partial_eD$.

Define $g'=\genus(\bar S_i)=\genus(\partial_\mE D)$.  We show that 
$g'\le g=\genus(\partial_\mE C)$, where $C$ is the original core of $M$.  By
construction $D\cap \partial M$ is a union of components of
$C\cap \partial M$, therefore it suffices to show 
that the number of 1-handles attached to $N(D\cap \partial M)$ is not 
more than the number of 1-handles attached to $N(C\cap \partial M)$ 
in the constructions of $D$ and $C$ respectively.
If $D\cap \partial M=C\cap \partial M$, then this follows immediately 
from the fact that $D$ and $C$ are cores respectively of $\mW$ and $M$ and 
the $H_1$-injectivity of $\mW$ in $M$.  
Let $E= (C\cap \partial M)\backslash (D\cap \partial M)$.  
The $H_1$-injectivity of $C$ in $M$ implies that the inclusion
$H_1(C\cap \partial M)\to H_1(M)$ is injective. 
The $H_1$-injectivity of $D$ in $\mW$ and the $H_1$-injectivity of $\mW$ in $M$
implies that $D$ is $H_1$-injective in $M$.  If the kernel of
$\text{in}_*:H_1(D\cup E)\to H_1(M)$ is 
nontrivial, then a non trivial homology between $D$ and $E$ would lie 
in some $V_j$, where $j>1$ and $V_j$ is a term in the exhausting 
sequence of $M$ used for constructing $\{W_i\}$.  
Arguing as in the proof of Lemma~\ref{injectivity},
$W_j$ contains a non trivial homology between $D$ and $E$ and hence,
$E\cap W_j\neq\emptyset$.  This
contradicts the fact that  $W_j\cap \partial M=D\cap \partial M$.  Therefore

\begin{align*}~\label{long_formula}
(*) \; \; \; 
\text{number of 1-handles of } D &\le \rank H_1(M) -
(\rank H_1(D\cap \partial M) + \rank H_1(E)) \\
&= \rank H_1(M) - \rank(H_1(\partial C \cap M)) \\
&= \text{number of 1-handles of } C \\
\end{align*}

Isotope $\bar S_i\subset \bar X_i$  to an embedded 
surface  $\hat S^i\subset X_i$ via an isotopy which does not cross 
$\hat \Delta_i$. Next, if possible, compress $\hat S^i$ via  a 
compression either disjoint from $\hat \Delta_i$ or
crossing $\hat \Delta_i$ once, say at $\hat\delta_{i_1}$.  If possible,
compress the resulting surface via a compression crossing
$\hat \Delta_i\backslash \hat\delta_{i_1}$ 
at most once and so on. Since $\genus(\hat S^i)=g'$, there
is an {\it a priori} upper  bound on the number of compressions
we need to do. In the end we obtain connected surfaces
$\hat S^i_1,\cdots, \hat S^i_n$ in $X_i$ 
which are 2-incompressible rel 
$\hat \Delta \backslash \{\hat\delta_{i_1},\cdots, \hat\delta_{i_m}\}$ where 
$m<2g'-1$, and both $n$ and $\genus(\hat S^i_j)$ are $\le g'$. Since 
$X_i\backslash \hat\Delta_i$ is irreducible,
we can assume that no $\hat S^i_j$ is a 2-sphere. These
$\hat S^i_j$'s create a partition $B^i_1,\cdots B^i_n$ of
$\hat \Delta_i\backslash \{\hat\delta_{i_1},\cdots, \hat\delta_{i_m}\}$,
where $B^i_j$ is the subset of
$\hat\Delta_i\backslash \{\hat\delta_{i_1},\cdots, \hat\delta_{i_m}\}$ separated
from $\bar S_i$ by   $\hat S^i_j$. Each $\hat S^i_j$ is incompressible to the
$\bar S_i$ side, since the component of $\bar X_i$ split along 
$\hat S^i_1\cup\cdots\cup\hat S^i_n$ 
which contains $\bar S_i$ is homeomorphic to 
$N(\hat S^i_1\cup\cdots \cup\hat S^i_n)\cup 1$-handles.  
Therefore, each $\hat S^i_j$ is 2-incompressible rel $B^i_j$.

As in the proof of Canary's theorem, after appropriately reordering the
$B^i_j$'s we can find  a $p\in \BN$ and a sequence $k_1<k_2<\cdots$ such that
$\hat \delta_p\subset B^{k_i}_1$ and
if $p(i)$ denotes the largest index of a $\hat\delta_j \in B^{k_i}_1$, then
$\lim_{i\to \infty} p(i)=\infty$.  In general reorder the 
$\hat S^i_j$'s so that, if possible, $\hat\delta_p\subset B^i_1$.

Fix $i$.  Let $W^\prime_i$ be the union of $W_i$ together with the
components of $N\backslash \inte(M)$ which non
trivially intersect $\partial W_1$.  Let $ Y_i$ denote the covering
of $W^\prime_i$ with fundamental
group $G_i$.  View $X_i, \hat \Delta_i$, and the $\hat S^i_j$'s etc. as
sitting naturally in $Y_i$.  Let $\delta\in \Delta$ be disjoint from
$W_i$.  Apply Lemma~\ref{construction_lemma} to $W_i,\ \delta\cup \Delta_i$ and
$S^i_1$.

We have the following dictionary between terms appearing in
our current setup (on the left) and the terms appearing
in the hypothesis of Lemma~\ref{construction_lemma} (on the right):

\begin{eqnarray*}
\text{the geodesics }\delta \cup \Delta_i & \longleftrightarrow & 
\text{the geodesics }\Delta_1 \\
\text{the manifold }W_i' & 
\longleftrightarrow & \text{the manifold }W \\
W_i' \cap (\delta \cup 
\Delta_i) = \Delta_i &\longleftrightarrow & 
W \cap \Delta_1 = \Delta \\
\text{the subgroup }G_i\text{ of }\pi_1(W_i') & 
\longleftrightarrow & 
\text{the subgroup }G\text{ of }\pi_1(W) \\ 
\text{the cover }Y_i \text{ with }\pi_1(Y_i) = G_i& 
\longleftrightarrow & 
\text{the cover }X\text{ with }\pi_1(X) = G \\
\text{the lifted geodesics }B^i_1 & \longleftrightarrow & 
\text{the lifted geodesics }B \\
\text{the surface }\hat S^i_1 & 
\longleftrightarrow & \text{the surface }S \\
\end{eqnarray*}

Then Lemma~\ref{construction_lemma} constructs surfaces 
$T^i$ and $P^i$ where the correspondence is
\begin{eqnarray*}
\text{the shrinkwrapped surface }T^i & \longleftrightarrow & 
\text{the shrinkwrapped surface }T \\
\text{the approximating surface }P^i & \longleftrightarrow & 
\text{the approximating surface }T_t \\
\end{eqnarray*}

In more detail: $W'_i$ is  isotopic to a manifold $W_i^{\new}$, via 
an isotopy fixing $\Delta_i\cup \delta$ pointwise.   This isotopy 
induces a homotopy of the
covering projection $\pi:Y_i\to W'_i\subset N$ to a covering projection
$\pi^{\new}:Y_i\to   W_i^{\new}\subset N$.  Our $\hat S^i_1$ is 
isotopic to a surface $\hat P^i$ via an 
isotopy avoiding $B^i_1$ and the projection of $\hat P^i$ into $N$ is 
a surface $P^i$ which is
homotopic to a $\CAT(-1)$ surface $T^i$. 
Furthermore, $P^i$ and $T^i$ are at
Hausdorff distance $\le 1$ and 
the homotopy from $P^i$ to $T^i$ is 
supported within the 1-neighborhood of $P^i$.

We relabel superscripts, and by abuse of 
notation we let the sequence $\{T^i\}$ 
stand for the old subsequence
$\{T^{k_i}\}$, with $\delta_{p(k_i)}$ being denoted by
$\delta_{p(i)}$, etc. We also drop the superscript \emph{new} so in
particular the projection $\pi:Y_i\to W'_i$ now refers to $\pi^{\new}:Y_i\to W^{\new}_i$.

We use the $\delta_{p(i)}$'s to show that $\{T^i\}$ exits
$\mE$. Let $\{\alpha_i\}$ be a locally finite collection of properly embedded  rays
from $\{\delta_{i}\}$ to $\mE$. For each $i$,  
$\hat S^i_1$ intersects some component $\omega$ of
$\piinv(\alpha_{p(i)})$
with algebraic intersection number 1, so $\hat P^i\cap\piinv(\alpha_{p(i)})\neq\emptyset$.   
Therefore  for all $i$, we have an inequality
$\dist(T^i,\alpha_{p(i)})\le 1$.   Our assertion now follows from the Bounded
Diameter Lemma.

\begin{figure}[ht]
\centerline{\relabelbox\small 
\epsfxsize 3.0truein
\epsfbox{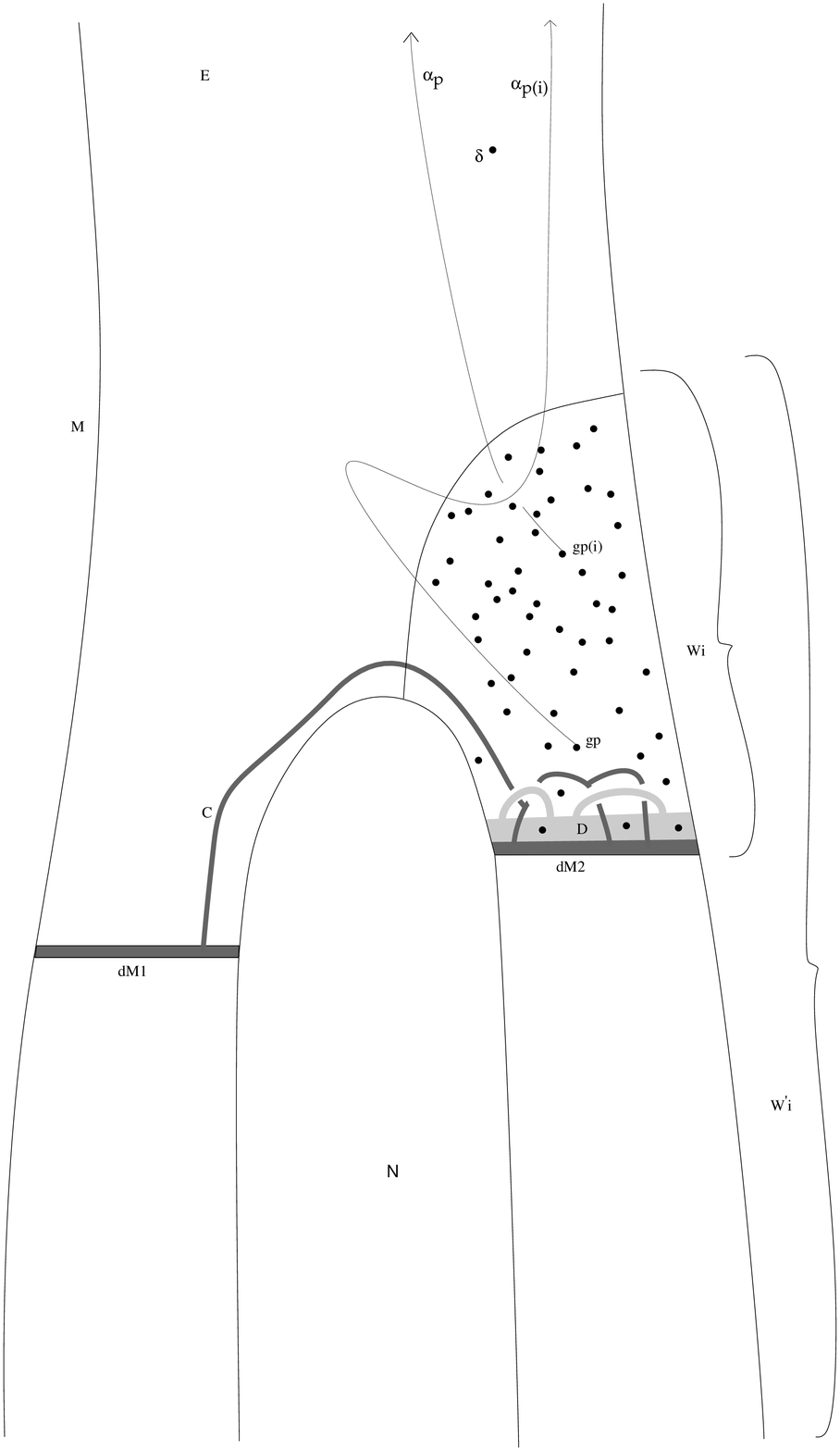}
\relabel {E}{$\mE$}
\adjustrelabel <0pt, -8pt> {dM1}{$\partial M$}
\adjustrelabel <0pt, -8pt> {dM2}{$\partial M$}
\adjustrelabel <-6pt, 0pt> {C}{$C$}
\relabel {D}{$D$}
\adjustrelabel <14pt, 0pt> {Wi}{$W_i$}
\extralabel <0.0truein, 1.65truein> {$W_i'$}
\adjustrelabel <-8pt, 0pt> {N}{$N$}
\relabel {gp}{$\delta_{p}$}
\relabel {gp(i)}{$\delta_{p(i)}$}
\adjustrelabel <-3pt, -3pt> {d}{$\delta$}
\extralabel <-1.9truein, 5.0truein> {$\alpha_{p(i)}$}
\extralabel <-1.3truein, 5.0truein> {$\alpha_{p}$}
\endrelabelbox}
\caption{}
\label{61}
\end{figure}

\begin{lemma}\label{norm}  Let $\mE$ be an end of $M$  an orientable,
irreducible  3-manifold with finitely generated fundamental group.  If $C$ is a
3-manifold compact core of $M$
and $Z$ is the component of $M\backslash C$ which contains $\mE$, then
$[\pC]$ generates $H_2(Z)$ and is Thurston norm minimizing. Here $\pC$ is the
component of $\partial C$ which faces $\mE$.\end{lemma}

\begin{proof}  First,  $\pC$ is connected, or else there exists a closed
curve $\kappa$ in $M$ intersecting a component of $\pC$ once, hence $\kappa$
is not homologous to a cycle in $C$, contradicting the fact 
that $C$ is a core.  That $[\pC]$ generates
$H_2(Z)$ follows from the fact that any 2-cycle $w$  in $Z$ is
homologous to one in $C$, so the restriction of that homology to $Z$ 
gives a homology of $w$ to $n[\pC]$ for some
$n$.  Equivalently, observe that the inclusion $C \to M$ is a homotopy
equivalence, and use excision for homology.

Let $Q\subset \inte(Z)$ be a Thurston norm minimizing surface 
representing $[\pC]$. We can choose $Q$ to be connected since $H_2(Z)=\BZ$.  
Let $V\subset Z$ be the submanifold between $\pC$
and $Q$.  If $\genus(\pC)>\genus(Q)$, then there exists a nonzero $z$
in the kernel of  $in_*:H_1(\pC)\to H_1(V)$.
This follows from the well known fact that for any
compact orientable 3-manifold $V$, the rank of the kernel of the map
$in_*:H_1(\partial V)\to H_1(V)$ is $\frac{1}{2}
\rank (H_1(\partial V))$.  Since
$C$ is a core, $z$ is in the kernel of the map $in_*:H_1(\pC)\to H_1(C)$.
This gives rise to a class $w'\in H_2(M)$ and dual class
$z'\in H_1(\pC)$  with $\langle z',w' \rangle \neq 0$, which is again 
a contradiction.
\end{proof}

Using this lemma, we now complete the proof of Theorem~\ref{main}.

Let $Z$ denote the component of $N$ split open along $\pC$ which
contains $\mE$. By Lemma \ref{norm}, $\pC$ generates $H_2(Z)$.  Next,
observe that if $\beta$ is any ray in $Z$ from
$\partial_\mE C$ to $\mE$ and $R$ is any immersed closed 
orientable surface in $Z$, then $[R]=n[\partial_\mE C]\in H_2(Z)$ 
where $\langle R,\beta \rangle=n$.  
To see this, note that $\langle n[\pC],\beta \rangle =n$ by
considering $n$ copies of $\pC$ slightly 
pushed into $Z$ and algebraic intersection number is independent of
representative of homology class.

We now use $\delta_p$ to show that for
$i$ sufficiently large $[T^i]$ is homologous in $Z$ 
to $[\partial_{\mE}C]\in H_2(Z) \cong \Z$; see Figure~\ref{61}. Let
$\beta$ be the ray $\sigma*\alpha_p$ where $\sigma\subset Z$ is a path 
from $\partial_\mE C$ to $\partial \alpha_p$.  For what follows 
assume that $i$ is sufficiently large so 
that
$$N_2(P^i)\cap(\sigma\cup\delta_p\cup C)=\emptyset$$
where $N_2(P^i)$ denotes the $2$--neighborhood of $P^i$, 
and hence
$$\langle T^i,\beta \rangle = \langle P^i,\alpha_p \rangle$$ 

We now compute this value.  By  perturbing $P^i$, if
necessary, we can assume that $P^i$ is transverse to $\alpha_p$ and 
no intersections occur at double points of $P^i$.  There is a $1$--$1$ 
correspondence of sets
$$\{\alpha_p\cap P^i\} \longleftrightarrow \{\piinv(\alpha_p) \cap \hat P^i\}$$  
Let $\Bag_i$ denote the component of $Y_i$ 
split along $\hat P^i$ which is disjoint from $\bar S_i$.  Note that
$\piinv(\delta_p)\cap \partial \Bag_i=\emptyset$.  If
$\kappa$ is a component of $\piinv(\alpha_p)$, 
then $\langle \kappa,\hat P^i \rangle=0$ if no endpoints
lie in $\Bag_i$ while $\langle \kappa,\hat P^i \rangle =1$ 
if exactly one endpoint lies in $\Bag_i$. 
To see this, orient $\alpha_p$ so that the positive end escapes to $\mE$.
Then the positive end of each lift $\kappa$ is in $\partial Y_i$, which is outside
Bag$_i$. It follows that if $p$ is an endpoint of 
$\kappa$ in $\Bag_i$, then $p$ is the negative end of $\kappa$,
and $\langle \kappa,\hat P^i\rangle =1$.
Since $\hat\delta_p$ lies in $\Bag_i$, there is at least $1$ component
$\kappa$ of $\pi^{-1}(\alpha_p)$ with such an endpoint in $\Bag_i$, and 
therefore
$$\langle \piinv(\alpha_p),\hat  P^i \rangle \ge 1$$ and hence 
$$[T^i]=n[\partial_{\mE}C]\in H_2(Z)$$
for some $n\ge1$.

Therefore 
$$|\chi(\partial_\mE C)|\ge |\chi(T_i)| \ge 
x_s(n[\partial_\mE C])=x(n[\partial_\mE C])=nx([\partial_\mE C]) =
n|\chi(\partial_\mE C)|$$ 
and hence $n=1$ and $\genus(T^i)=\genus(\partial_\mE C)$.  
Here $x$ (resp. $x_s$) denotes the Thurston (resp. singular Thurston) norm
on $H_2(Z)$. The first inequality follows by construction, the second 
by definition, the third since $x_s=x$ (\cite{G1}), the fourth since 
$x$ is linear on rays \cite{T2} and the fifth by Lemma~\ref{norm}. 
This completes the proof of Theorem~\ref{main}.
\qed

\begin{remark}\label{barrier remark} 
Since for $i$ sufficiently large, $\genus(T^i)=g$,
it follows that for such $i$, no compressions occur 
in the passage from $\bar S^i$ to $\hat S^i_1$.
This mirrors the similar phenomenon seen in the proof of Canary's theorem.
If the shrinkwrapped $\partial W'_i$ is 
actually a $\Delta_i$-minimal surface
disjoint from $\Delta_i$, then $\partial X_i$ is
a least area minimal surface for the hyperbolic metric, and we
can pass directly from $\hat S^i_1$ to a 
$\hat \Delta_i$-minimal surface $\hat T^i$ by shrinkwrapping in
$X_i$.  Our $T^i$ is then the projection of $\hat T^i$ to $N$.

If the shrinkwrapped $\partial W'_i$ touches $\Delta_i$, then we can
still shrinkwrap $\hat S^i_1$ in $X_i$.  
In this case $X_i$  is bent and possibly squeezed along parts of
$\hat \Delta_i$ and it is cumbersome to discuss the geometry and topology
of $X_i$. Therefore we chose for the purposes of exposition to express $T^i$
as a limit of surfaces.  These surfaces are projections of $g_{t_k}$--minimal
surfaces in the smooth Riemannian manifolds $X_i$ with Riemannian
metrics $g_{t_k}$.  As metric spaces, the $(X_i, g_{t_k})$ converge
to the bent and squeezed hyperbolic ``metric'' on $X_i$.
\end{remark}

\noindent\textbf{Tameness Criteria}\label{criteria}  Let $\mE$
be an end of the  complete hyperbolic 3-manifold $N$ with finitely generated
fundamental group and compact core $C$.  Let $Z$ be the component of 
$N\backslash \inte (C)$ containing $\mE$ with $\pC$ denoting
$\partial Z$.   Let $T_1, T_2, \cdots$ be a sequence of surfaces mapped into $N$.  
Consider the following properties.
\begin{enumerate}\item $\genus(T_i)=\genus( \pC)$.
\item $T_i\subset Z$ and exit $\mE$.
\item Each $T_i$ homologically separates $C$ from $\mE$ (i.e.
$[T_i]=[\pC]\in H_2(Z))$.
\item Each $T_i$ is CAT(-1).
\end{enumerate}\vskip 12pt

\begin{theorem}[Souto \cite{So}]\label{souto} If $T_1, T_2, \cdots$
is a sequence of mapped surfaces in the
complete hyperbolic 3-manifold $N$ with core $C$ and end $\mE$ which
satisfies Criteria (1), (2) and (3), then $\mE$ is topologically tame.\qed\end{theorem}

Theorem \ref{souto} follows directly from the proof of Theorem 2,
\cite{So}.  That proof makes essential use of the work of Bonahon
\cite{Bo} and Canary \cite{Ca}.   We now show how Criterion (4) enables us to
establish tameness without
invoking the impressive technology of \cite{Bo} and \cite{Ca}.  Our argument,
inspired in part by Souto's work, requires only elementary hyperbolic geometry
and basic 3-manifold topology.\vskip 12pt

\noindent\textbf{A topological argument that criteria (1) - (4) imply
tameness.}\label{topological}  It suffices to
consider the case that
$\mE$ is the end of an end-manifold $M$ which includes by a homotopy
equivalence into $N$,
and that $C\subset M$ is of the form $\partial M\times I\cup 1$-handles, where
the 1-handles attach to $\partial M\times 1$ and 
$\partial M=\partial M\times 0$.  

Using standard arguments, we can replace
the $T_i$'s by \emph{simplicial hyperbolic surfaces} as defined in \cite {Ca}.
The idea of how to do this is simple: the $\CAT(-1)$ property
implies that each $T_i$ has an essential simple closed curve
$\kappa_i$ of length uniformly bounded above.  If $\kappa_i^*$
denotes the geodesic in $N$ homotopic to
$\kappa_i$, then either the $\kappa_i^*$ have length bounded below by some
constant, and are therefore contained within a bounded neighborhood 
of $\kappa_i$, or else the lengths of the $\kappa_i^*$ get arbitrarily short,
and therefore they escape to infinity. In either case, the
sequence $\kappa_1^*,\kappa_2^*,\cdots$ exits $\mE$. Then we can triangulate
$T_i$ by a $1$--vertex triangulation with a vertex on $\kappa_i^*$, 
and pull the simplices tight to geodesic triangles. This produces a 
simplicial hyperbolic surface, homotopic
to $T_i$, which is contained in a bounded 
neighborhood of $\kappa_i^*$ rel. the
thin part of $N$, and therefore 
these surfaces also exit $\mE$.  From now on we assume that each $T_i$ is a
simplicial hyperbolic surface.

Note that either $\pC$ is incompressible in $N$ and hence $M$ is homotopy 
equivalent to $\pC\times [0,\infty)$ or each $T_i$ is compressible in
$N$, i.e. there  exists an essential simple closed curve in $T_i$ that is 
homotopically trivial in $N$.  Indeed, using the $\pi_1$-surjectivity of $C$
and the irreducibility of $N,\ T_i$ can 
be homotoped into $C$.  If $T_i$ is
incompressible in $N$, then
$T_i$ can be homotoped off the 1-handles and then homotoped into
a component of $\partial M$.  Using Criterion (3), the degree of 
this map is one, which implies
that $T_i$ is homotopic to a homeomorphism onto  a 
component of $\partial M$.  Since
$\genus(T_i)=\genus(\pC)$, it follows that 
$C=\partial M\times I$ and hence $\pC$ is incompressible in $N$.

Since either $M$ is homotopy equivalent to $\partial M\times [0,\infty)$ or each
$T_i$ is compressible, it follows by Canary \cite{Ca} and Canary--Minsky
\cite{CM} (see also Proposition 3 \cite{So}) 
that there exists a compact set $K\subset Z$ such that  each $T_i$
can be homotoped within $Z$ to a simplicial hyperbolic surface $T^0_i$ which
nontrivially intersects $K$.  Here $Z$ is the closure of $M-\inte(C)$.
By the Bounded Diameter Lemma, there exists a compact set
$K_1\subset Z$ such that for each $i$, $T^0_i\subset  K_1$.   

Since  $x=x_s$ \cite{G1}, 
there exists a sequence of embedded genus-$g$ surfaces 
$A_1,A_2,\cdots$ such that for each
$i$,\ $A_i$ lies in a small 
neighborhood of $T_i$ and 
$[A_i]=[T_i]=[\pC]\in H_2(Z)$. 
By passing to subsequence we can assume that the $A_i$'s are pairwise 
disjoint and each $A_i$ is disjoint from $K_1$ and 
separates $\mE$ from $K_1$. Let $A_{[p,q]}$ 
denote the compact region between 
$A_p$ and $A_q$.   Since $\genus(A_p)=\genus(\pC)$, it follows by Lemma
\ref{norm} that $A_p$ is Thurston norm minimizing in $Z$ and hence is
$\pi_1$-injective in $Z$ and  $A_{[p,p+1]}$.   

To establish tameness it suffices to show that each $A_{[p,p+1]}$ is a 
product. Fix $p\in \BN$.  Let $j$ be
sufficiently large so that $T_j$ separates  $A_{[p,p+1]}$ 
from $\mE$. Let $T$ be a surface
of genus $g=\genus(\pC)$.   Using \cite{Ca}, \cite{CM}, let $F:T\times I\to Z$ 
be a homotopy such that $F|T\times 1=T_j$ and  $F(T\times 0)\subset K_1$.  
By Stallings and Waldhausen, after a  homotopy of  $F$ rel $\partial F$ 
we can assume that $F^{-1}(A_p\cup A_{p+1})$ are $\pi_1$-injective surfaces in
$T\times (0,1)$.  See p. 60 \cite{Wa}.  
Since $F(T\times \partial I)\cap A_{[p,p+1]}=\emptyset$, these surfaces are disjoint from 
$T\times \partial I$, and by arguing as in \S3 \cite{Wa}, they are isotopic to surfaces
of the form $T\times t$, $t\in (0,1)$.  Therefore, after a further homotopy we
can assume that $F^{-1}(A_p\cup A_{p+1})=T\times B$, where $B\subset (0,1)$ is a
finite set of points. Since each $F|T\times t$ homologically separates $\mE$
from $\pC$, each $F|T\times b$ is a degree-1 map onto either $A_p$ or $A_{p+1}$
and hence after another homotopy we can assume that for each
$b\in B$, $F|T\times b$ is a homeomorphism onto its image.  
Therefore there exists $b,b'\in B$ such that
$F|T\times[b,b']$ maps degree-1 onto $A_{[p,p+1]}$ and the restriction of $F$
to $T\times \partial [b,b'] $ is a homeomorphism. 
Therefore  $F:T\times[b,b']\to A_{[p,p+1]}$ is a
$\pi_1$-injective, degree-1 map whose restriction to $\partial(T\times[b,b'])$
is a homeomorphism onto $\partial A_{[p,p+1]}$.   Since both domain and range
are irreducible, such a map is homotopic rel boundary to a homeomorphism, by
Waldhausen \cite{Wa}.
\qed\vskip 12pt

\begin{remarks} 
\begin{enumerate}
\item{In the presence of an escaping sequence of CAT(-1) surfaces, 
\emph{hyperbolic surface interpolation} and the \emph{bounded diameter lemma} is all the
hyperbolic geometry needed to establish tameness.}
\item{This argument makes crucial use of the fact that the homotopy $F$ is
supported in $Z$ and each $A_i$ is incompressible in $Z$.}
\end{enumerate}
\end{remarks}

\noindent\textbf{Proof of Theorem \ref{marden}}  It suffices to
consider the case that $N$ is orientable, since it
readily follows using \cite{Tu}, that $N$ is tame if and only if its
orientable cover is tame. If $\mE$ is geometrically finite, then by \cite{EM} $\mE$ is tame.
Now assume that $\mE$ is geometrically infinite.  Theorem~\ref{main}
provides us with a collection $\{T_i\}$ which satisfies Tameness 
Criteria (1)-(4). Now apply Theorem~\ref{souto}. \qed\newline

\noindent\textbf{Proof of Theorem~\ref{clean}}.  It suffices to prove Theorem
\ref{clean} for the geometrically infinite ends of orientable
manifolds.  It follows from  Theorems
\ref{main} and \ref{marden} that  $\mE$ is topologically of the
form $T\times [0,\infty)$,
where $T$ is a surface of genus $g$.  Theorem \ref{main} provides for
us a sequence $\{T_i\}$ of surfaces satisfying Tameness Criteria (1)-(4).  Since
for $i$ sufficiently large
$T_i\subset T\times[0,\infty)$ and  homologically separates $T\times 0$ from
$\mE$, it follows that the projection
$T_i$ to $T\times 0$ is a degree 1 map of a genus $g$ surface to itself and
hence is homotopic to a homeomorphism.
\qed

\section{The  Parabolic Case}\label{parabolic}

Thanks to the careful expositions in \cite{Bo}, \cite{Ca} and 
\cite{So} it is now routine to obtain general theorems in the 
presence of parabolics
from the corresponding results in the 
parabolic free case.

We now give  the basic definitions and provide 
statements
of our results in
the parabolic setting.

The following is 
well known, e.g. see \cite {Ca} for an expanded
version of more 
or less the following discussion. Let $N$ be a complete 
hyperbolic 3-manifold, then for
sufficiently small $\epsilon$, the 
$\epsilon$-thin part, $\Nthin$ of $N$ is a union of
solid tori (Margulis tubes), rank-1 cusps and rank-2 cusps. Let $\Nthick$ 
denote $N\backslash \inte(N_{\le \epsilon})$.
The space $N_0^\epsilon = \Nthick\cup$ Margulis tubes is called the
\emph{neutered space} of $N$, though we often drop the $\epsilon$.
The parabolic locus
$\partial N_0^\epsilon=P^\epsilon$ (usually just denoted $P$) 
is a finite union of tori $T_1,\cdots,T_m$ and open 
annuli $A_1,\cdots A_n$.  Each annulus $A_i$ is of the
form $S^1\times \BR$ such that for $t\in \BR$, each $S^1\times t$ bounds a 
standard 2-dimensional cusp in $\Nthin$. Let $N^0_{\le \epsilon}$ denote
the cusp components of $N_{\le \epsilon}$.
By \cite {Mc}, $N$ has a compact core $C\subset N_0$
which is also a core of $N_0$ and the 
restriction to each component
$P'$ of $P$ is a core
of $P'$.   Such a core for $N_0$ is called a \emph{relative core}. In
particular if $P'$ is an annulus, then we can assume that
$C\cap P'=S^1\times[t,s].$  By Bonahon \cite{Bo}, the ends of $N_0$ are in 1-1
correspondence with components of $\partial C/P$.   If $N=\BH^3/\Gamma$, then an
end of $N_0$ is \emph{geometrically finite} if it has a neighborhood disjoint 
from $C(\Gamma)/\Gamma$, the convex core of $N$.  Such an end has an 
exponentially flaring geometry similar to that of a geometrically 
finite end of a parabolic free manifold.  The end $\mE$ of $N_0$ 
is \emph{topologically tame} if
it is a relative product, i.e. there 
is a compact surface $S$ and an
embedding $S\times [0,\infty)\to N_0$ 
which parametrizes $\mE$.   If $U$ is a
neighborhood of $\mE$,
then by passing to a smaller neighborhood we can assume that $U\cap A_i$ 
is either $\emptyset$ or of the form
$S^1\times (t,\infty)$ or 
$S^1\times ((-\infty,s)\cup (t,\infty))$.  Adding the
corresponding 2-dimensional cusps to $S^1\times$pts. we obtain
$U_{P}$ the $\emph{parabolic extension}$ of $U$.  So if $\mE$ is 
topologically tame, $U_{P}$ is topologically $S^P\times [0,\infty)$, 
where $S^P$ is topologically $\inte(S)$ 
and geometrically $S$ with cusps added.

Following \cite{Bo} and \cite{Ca} we say that the end
$\mE$ of $N_0$ is \emph{simply degenerate} if it is topologically
tame, has a neighborhood $U$ with 
a sequence $f_i:S^P\to U_{P}$
such that $f_i$ induces a
$\CAT(-1)$ structure on $S^P$, the $f_i$'s eventually miss given 
compact sets and each $f_i$ is properly homotopic in $U_{P}$ to a 
homeomorphism of $S^P$
onto $S^P\times 0$.  We say that $\mE$ is 
\emph{geometrically tame} if it is simply degenerate or
geometrically finite.  The manifold $N$ is \emph{geometrically tame}
if each end of $N_0$ is geometrically tame.

Francis Bonahon showed that if 
$\epsilon$ is sufficiently small, then an end
$\mE$ of $N^\epsilon_0$ 
is geometrically infinite if and  only if
there exists a
sequence $\Delta=\{\delta_i\}$ of closed geodesics lying in
$N^\epsilon_0$ and exiting $\mE$.

We can now state the general version of the results stated in
the introduction.

\begin{theorem}\label{clean parabolic}  Let $N$ be a 
complete hyperbolic
3-manifold with finitely generated fundamental 
group with neutered space $N_0$.
The end
$\mE$ of $N_0$ is simply 
degenerate if there exists a sequence of closed
geodesics exiting the 
end.\end{theorem}

\begin{theorem}\label{geometrically tame 
parabolic}  A complete
hyperbolic 3-manifold with finitely generated 
fundamental group is
geometrically 
tame.\end{theorem}

\begin{theorem} \label{marden parabolic}  If $N$ 
is a complete hyperbolic
3-manifold with finitely generated 
fundamental group, then each end
of $N_0$  is
topologically tame.  In 
particular, each end of $N$ is topologically 
tame.
\end{theorem}

\begin{theorem}\label{ahlfors parabolic}  If 
$N=\BH^3/\Gamma$ is a complete
hyperbolic 3-manifold with finitely 
generated fundamental group, then the
limit 
set $L_\Gamma$ of $\Gamma$ is either $\Sinfty$ or has Lebesgue measure zero.  
If $L_\Gamma=\Sinfty$, then $\Gamma$ acts ergodically 
on
$\Sinfty$.\end{theorem}

\begin{theorem}[Classification Theorem] 
\label{classification
theorem parabolic} If $N$ is a complete 
hyperbolic 3-manifold with finitely
generated fundamental group, then 
$N$ is determined up to isometry
by its topological type, its 
parabolic structure, the
conformal boundary of $N_0$'s geometrically 
finite ends and the ending
laminations of $N_0$'s 
geometrically
infinite ends. \end{theorem}

\begin{theorem}[Density Theorem]\label{density 
parabolic}  If
$N=\BH^3/\Gamma$ is a complete finitely generated 
3-manifold with finitely
generated fundamental group, then
$\Gamma$ 
is the algebraic limit of geometrically finite Kleinian 
groups.\end
{theorem}

\begin{theorem}\label{main parabolic} Let 
$N$ be a complete
hyperbolic 3-manifold
with finitely generated 
fundamental group and with  associated
neutered space $N_0$.   Let 
$\mE$ be an end of $N_0$ with relative compact core
$C$.  Let $S$ be a compact surface with the topological type of
$\delta_{\mE}C$, the component of the frontier of $C$ which faces
$\mE$.   Let $U_\mE$ denote a parabolic extension of 
a neighborhood of $\mE$.  If there exists a sequence of
closed geodesics exiting $\mE$, then there 
exists a sequence $\{S_i\}$ of proper
$\CAT(-1)$ surfaces in 
$U_{P}$ homeomorphic to $\inte(S)$
which eventually miss every 
compact set and
such that each $S_i\cap N_0$  homologically separates 
$C$ from $\mE$. Furthermore, if $S_i\cap N_0$ lies to the $\mE$ side
of $C$, then no accidental parabolic $\alpha\subset S_i\cap N_0$ can be
homotoped into a cusp via a homotopy disjoint from $C$.
\end{theorem}

\begin{remarks} \begin{enumerate}
\item Theorem \ref{main parabolic} is the main technical result of this section
and at the end of this section we will deduce from it Theorems \ref{clean
parabolic} and \ref{marden parabolic}.
\item Theorem~\ref{marden parabolic} has been independently proven by Agol
\cite{Ag}.
\item Theorem \ref{clean parabolic} implies Theorem \ref{geometrically tame
parabolic} as follows.  A complete hyperbolic 3-manifold is geometrically tame if
each end of
$N_0$ is either geometrically finite or simply degenerate.  By definition, ends
of $N_0$ are either geometrically finite or geometrically infinite.  Using
Bonahon's characterization of geometrically infinite ends and Theorem \ref{clean
parabolic} it follows that geometrically infinite ends are simply degenerate.
\item Theorem  \ref{ahlfors parabolic} immediately follows from Theorem
\ref{geometrically tame parabolic} by the work of Thurston
\cite{T} and Canary \cite{Ca}.  It also follows from \cite{Ca} that the  various
results of \S 9
\cite{Ca} hold  for $N$.
\item Theorem~\ref{marden parabolic} is the last step needed to  prove
the monumental classification theorem, the other parts being established by
Alhfors, Bers, Kra, Marden, Maskit, Mostow, Prasad, Sullivan,
Thurston, Minsky, Masur--Minsky, Brock--Canary--Minsky, Ohshika,
Kelineidam--Souto, Lecuire, Kim--Lecuire--Ohshika, Hossein--Souto and Rees. See
\cite{Mi} and \cite{BCM}.
\item The Density Theorem was conjectured by Bers, Sullivan and Thurston.
Theorem \ref{marden parabolic} is one of very many results, many of them recent,
needed to build a proof.  Major contributions were made by Alhfors, Bers, Kra,
Marden, Maskit, Mostow, Prasad, Sullivan, Thurston, Minsky, Masur--Minsky,
Brock--Canary--Minsky, Ohshika, Kelineidam--Souto, Lecuire,
Kim--Lecuire--Ohshika, Hossein--Souto, Rees, Bromberg and Brock--Bromberg.
\item The rest of this section is devoted to proving Theorems
\ref{main parabolic}, \ref{clean parabolic} and \ref{marden parabolic}.
\end{enumerate}\end{remarks}

Given the manifold $N$ with neutering $N_0$ and end $\mE$ of $N_0$ 
we explain how to find a \emph{relative end-manifold} $M$ containing $\mE$.

\begin{definition} If $A$ is a cod-0 submanifold of a 
manifold with boundary, then the
\emph{frontier} $\delta A$ of $A$ is 
the closure of $\partial A-\partial M$.  If $(R,\partial R) \subset 
(N_0, P)$ is a mapped surface (resp. $R$ is a properly mapped surface 
in $N$ whose ends exit the cusps), then a 
$P$-\emph{essential annulus} for $R$ is annulus (resp. half open annulus) $A$   
with one component mapped to an essential simple curve of $R$ 
which cannot be homotoped within $R$ into
$\partial R$ (resp. an end of $R$) and another
component (resp. the end of $A$) mapped into $P$ 
(resp. properly mapped into a cusp). Let
$C_0$ be a 3-manifold relative core of
$N_0$.  Using \cite{Mc} we can assume that $C_0$ is 
of the form $H_0\cup$ 1-handles where $P_0 =C_0\cap 
P=H_0\cap P$ is a core of $P$ consisting of annuli and tori and 
$H_0$ is a compact 3-manifold with incompressible frontier.  Furthermore
$\delta H_0$ has no  $P$-essential
annuli disjoint from $\inte(H_0)$.  Define
$\delta_{\mE}C_0$ to be the 
component of $\delta C_0$ which faces
$\mE$ and  $\delta_{\mE}H_0$ 
to be the components of $\delta H_0$ which face
$\mE $.  Define $M$ to be the closure of the component of
$N_0$ split along $\delta_{\mE}H_0$ which contains
$\mE$. Define $\partial_pM=P\cap M$ and
$\partial_hM=\delta_{\mE}H_0$.  We call 
$M$ a \emph{relative end-manifold}.  By construction $\partial_hM$ 
has no $P$-essential annuli lying in $M$.  By slightly 
thickening $\partial_hM$ and retaining the 1-handles of $C_0\cap M$ 
we obtain a core $C$ of $M$.   If $W$ is a codimension-0 submanifold of $M$, 
then $\partial_pW$ (resp. $\partial_hW$) denotes $W\cap \partial_pM$ 
(resp. $W\cap \partial_hM$).
\end{definition}

By passing to the $\pi_1(M)$ cover of $N$ we reduce to the case that
$\text{in}:M\to N_0$  is a homotopy equivalence;  
furthermore, for each component $R$ of $\partial_hM$, the
inclusion of $R$ into the corresponding component of 
$N_0\backslash \inte(M)$ is a homotopy equivalence.  

\vskip 12 pt

By passing to a subsequence 
we can assume that $\Delta =\{\delta_i\}$ is a
collection of geodesics escaping $\mE$ and is weakly 1000-separating.
As in  Lemma 5.5 \cite{Ca} we slightly perturb the hyperbolic metric
in the $1$-neighborhood of $\Delta$ to a metric $\mu$ such that for each 
$i$, $\delta_i$ is $\epsilon$-homotopic to a simple geodesic 
$\gamma_i$  and $\mu$ has pinched
negative curvature in $(-1.01,-.99)$ and is $1.01$-bilipshitz equivalent to
the hyperbolic metric.  Let
$\Gamma$ be the resulting collection of simple closed 
curves.

\begin{lemma}  \label{end-engulfing parabolic} 
Let $M$ be a relative end-manifold
in the complete hyperbolic $3$--manifold $N$.  Given a 
sequence $\Gamma$ of homotopically essential  closed curves we 
can pass to an infinite subsequence also called $\Gamma$ which is 
the disjoint union $\gamma_1\cup\gamma_2\cup\gamma_3\cup\cdots$ 
where $\gamma_1$ has finitely many components and the other 
$\gamma_i$'s have one component.  If
$\Gamma_i$ denotes $\cup_{j=1}^i\gamma_j$, then there exists a
manifold $\mW$ open in $M$, exhausted by a sequence of compact manifolds
$\{W_i\}$ with the following properties,

\begin{enumerate}

\item $\mW$ is $\pi_1$ and 
$H_1$ injective (in $\BZ$ and $\BZt$ coeficients) in $M$ and 
hence $N$.

\item   For all $i, \partial_h W_1 = \partial_h W_i$ 
and is a union of components of $\partial_hM$. At most 
one component of $\partial_pW_1$ can lie in a component of 
$\partial_pM$.  For all $i, \partial_pW_i$
is a union of essential annuli, each of which 
contains a component of $\partial_pW_1$.
The frontier $\delta W_i$ is connected, 
separates $\Gamma_i$ from $\mE$ and is 2-incompressible rel 
$\Gamma_i$.

\item There exists a compact submanifold core 
$F\subset W_1$ of $\mW$ such that
each $\Gamma_i$ can be homotoped 
into $F$ via a homotopy supported in $W_i$.   $F$ is of the 
form $(W_1\cap\partial M)\times I$
with 1-handles attached to the 1-side.  Finally
$|\chi(\delta F)|\le|\chi(\delta C)|$.

\end{enumerate}\end{lemma}

\begin{proof}  Except for the last inequality, this lemma is just 
the relative form of that part of 
Theorem~\ref{refined} which was 
used to prove Theorem~\ref{main}.  
Let $J$ be a
connected compact 
set and $V_1\subset V_2\subset V_3\cdots$ an
exhaustion of $M$ such 
that
$\partial_h M \subset
\partial V_1$ and the $\partial_p V_i $ 
are the tori  of $M$ and essential
annuli  which meet each annular 
component of
$\partial_p M $ in exactly one component.   Define the 
\emph{relative
end reduction} $\mW_J$ of $J$ to
be the manifold 
exhausted by submanifolds $\{W_i\}$ where $V_i$
passes to $W_i$ via 
the
operations of compression, 2-handle addition, deletion and 
isotopy, where the
compressions and
2-handle additions are done only 
to $\delta V_i $ and its successors.
The same
arguments as before 
show that $\mW_J$ is both $\pi_1$ and
$H_1$-injective and as before 
we
can define   relative notions of end non separable and 
end-engulfing respectively for
finite and locally finite infinite 
collections of homotopy essential pairwise
disjoint closed 
curves. Similarly, since $\mW$ is a relative end manifold its core can be taken
to be of the form stated in 3.   The last inequality is the relative version of 
inequality $(*)$ from the proof of Theorem~\ref{main}.
\end{proof}

We now show that each $W_i$ is an atoroidal Haken manifold with negative Euler characteristic
by showing that every embedded torus incompressible in $W_i$ is boundary parallel and some
component of $\partial W_i$ is not a 2-sphere.  If $T$ is an embedded torus in $N$, then either
$T$ cuts off a rank-2 cusp or
$T$ is compressible.  Therefore, if each component of $\partial W_i$ is a torus incompressible
in
$N$, then $N$ has finite volume and Theorem \ref{main parabolic} holds.  The proof of the Claim
of
\S6 shows that no component of $\partial W$ is a torus compressible in $N$.  Therefore, some
component of $\partial W_i$ has genus $\ge 2$.  Since   tori incompressible in $N$ cut off
rank-2 cusps, any nonboundary parallel torus $T$ in $W_i$ must be compressible in $N$.  To show
that
$T$ is compressible in $W_i$, it suffices to show that it is compressible in $W_i^P$ the
parabolic extension of $W_i$ in $N$.  Now $\partial W_i^P$ is a finite union of properly
embedded surfaces in $N$ which are incompressible in $N\backslash \Gamma_i$ and $N$ has pinched
negative sectional curvature.  The proof of the Claim now applies to show that every embedded
torus in $W_i^P$ which is compressible in $N$ is also compressible in $W_i^P$. 

\begin{lemma} \label{controlled spilling} If $(E,\partial E)\subset
(N_0,P)$ is a compact
$\Gamma$-minimal surface (possibly non embedded), then $E$ cannot be
homotoped rel $\partial E$ into $P$.

Suppose $f:R\to N$ is a properly mapped
$\Gamma$-minimal surface such that for each $\epsilon>0,\
\finv(N_0^\epsilon)$ is compact and $R$ has no $P$-essential annuli
disjoint from $\Gamma$.  If $f$ is transverse to $N^\epsilon_0$, then each
component of
$R\backslash \inte(\finv(N^\epsilon_0))$ is either a compact disc or
a half open
annulus.\end{lemma}

\begin{proof} If such a homotopy exists, then   the lift
$\tilde E$ of $E$ to $\BH^3$ has the property that there exists a
closed horoball
$H$ with $\partial \tilde E\subset \textint(H)$ and $E\cap \partial 
H\neq\emptyset$.
This violates the maximum principle.

Therefore if $\sigma$ is a component of $\finv(P)$ which bounds a
disc $D$ in $R$, then
$f(D)\cap N_0\subset P$.  If $\sigma$ is essential in
$R$, then $\sigma$ can be homotoped into an end of $R$, since there
are no
$P$-essential annuli for $R$ disjoint from $\Gamma$.   Again the
maximum principle implies that the entire annular region bounded by
$\sigma$ is mapped into a component of $N\backslash \inte(N_0)$.\end{proof}

If $W \subset N$ is a codimension $0$ submanifold, we say that $W$
has {\em standardly embedded cusps} if the restriction of $W$ to each Margulis
tube neighborhood of a cusp of $N$ is either the entire cusp, or is a finite
union of products $\text{annulus} \times \R^+$ where the product structure
is compatible with the product structure on the cusp. If $N$ is obtained from
a complete hyperbolic manifold by neutering, then this neutering should restrict
to a neutering of $W$.

Here is the parabolic version of Lemma~\ref{construction_lemma}. The reader
may want to refresh their memory by first rereading 
Lemma~\ref{construction_lemma}:

\begin{lemma}[Parabolic construction lemma]\label{construction_lemma_parabolic}
Let $\mE$ be an end of the complete open orientable irreducible Riemannian 
$3$--manifold $N$ with metric $\mu$,
with finitely generated fundamental group, and neutering $N_0$ with 
parabolic locus $P$. Let $W \subset N$ be a submanifold such that $\partial W \cap 
\inte(N)$ separates $W$ from $\mE$, and whose ends are standardly embedded cusps in the 
cusps of $N$. Let $\Delta_1 \subset N_0\backslash \partial W$ be a finite 
collection of simple closed geodesics with $\Delta = W \cap \Delta_1$ a non--empty proper
subset of $\Delta_1$. Suppose futher that $\partial W$ is 
$2$--incompressible rel. $\Delta$ and has no $P$--essential annuli disjoint from $\Delta_1$.

Let the Riemannian metric $\mu$ on $N$ agree with a hyperbolic metric outside
tubular neighborhoods $N_\epsilon(\Delta_1)$ and inside tubular neighborhoods
$N_{\epsilon/2}(\Delta_1)$, having $\Delta_1$ as core geodesics, and such that
$\mu$ is a metric with sectional curvature pinched between $-1.01$ and $-0.99$.

Let $G$ be a finitely generated subgroup of $\pi_1(W)$, and let $X$ be the
covering space of $W$ corresponding to $G$. Let $\Sigma$ be the preimage of
$\Delta$ in $X$ with $\hat{\Delta} \subset \Sigma$ a subset which maps
homeomorphically onto $\Delta$ under the covering projection, and let
$B \subset \hat{\Delta}$ be a nonempty union of geodesics. Suppose there exists
a properly embedded surface $S \subset X \backslash B$ of finite topological
type, whose ends are standard cusps in the cusps of $X$ such that $S$ is
$2$--incompressible rel. $B$ in $X$ and has no $P$-essential annuli disjoint
from $B$, and which separates every component of $B$ from $\partial X$.

Then $\partial W$ can be properly homotoped to a $\Delta_1$--minimal 
surface which, by abuse of notation, we call $\partial W'$, 
and the map of $S$ into $N$ given by the covering projection is
properly homotopic to a map whose image $T'$ is $\Delta_1$--minimal 
and whose ends exit the cusps of $N$. 

Also, $\partial W'$ (resp. $T'$) can be perturbed by an
arbitrarily small perturbation to be an embedded (resp. smoothly 
immersed) surface $\partial W_t$ (resp. $T_t$) bounding $W_t$ with 
the following properties:
\begin{enumerate}
\item{There exists a proper isotopy from $\partial W$ to $\partial W_t$ which
never crosses $\Delta_1$, and which induces a proper isotopy from $W$ to 
$W_t$, and a corresponding deformation of pinched negatively curved 
manifolds $X$ to $X_t$ which fixes $\Sigma$ pointwise.}
\item{There exists a proper isotopy from $S$ to $S_t \subset X_t$ which never 
crosses $B$, such that $T_t$ is the projection of $S_t$ to $N$.}
\item{Each of the limit surfaces $F \in \lbrace\partial 
W',T'\rbrace$ {\em relatively exits the manifold as its restriction 
exits the neutered part.} That is, if $\C$ is a rank $1$ cusp 
foliated by totally geodesic $2$ dimensional cusps $C \times \R$
perpendicular to the boundary annulus $S^1 \times \R$, then if 
the intersection of $F$ with $\partial \C$ is contained in the region 
$S^1 \times [t,\infty)$, the intersection of $F$ with $\C$ is 
contained in the region $C\times [t,\infty)$, and similarly if the 
intersection is contained in $S^1 \times (-\infty,t]$.}
\end{enumerate}
\end{lemma}
\begin{proof}
The essential differences between the statements of 
Lemma~\ref{construction_lemma} and 
Lemma~\ref{construction_lemma_parabolic} are firstly that the metric 
in the parabolic case is pinched, so that the geodesics can be chosen to 
be simple; secondly that the surfaces in question are all properly 
embedded, and the isotopies and homotopies
are all proper; and thirdly that the limit surfaces relatively exit the
manifold as their restriction exits the neutered part.

These issues are all minor, and 
do not introduce any real complications in the
proof. The only question whose answer might not be immediatly apparent is how 
to perturb the metric $\mu$ to the $g_t$ metrics near cusps; it 
turns out that this is straightforward to do, and technically easier 
than deformations along geodesics, since the perturbed metrics 
actually have curvature bounded above by $0$.

We will find an exhaustion of $N$ by increasingly larger neutered spaces
$N^t_0$, each endowed with a metric $g_t$,  which is obtained from
the $\mu$-metric by deforming it along the geodesics $\Delta_1$ and along 
$\partial N^t_0$.  Our $\partial W_t$ will restrict to $g_t$-area 
minimizing representatives of the isotopy class of
$\partial W\cap N^t_0$.
The convergence and regularity of the limit surface $\partial W'$
near the geodesics will proceed  exactly as in \S 1 and \S 2. 
The convegence and regularity in the cusps will
follow from \S 1 using the absence of $P$-essential annuli disjoint
from $\Delta_1$.

To describe the deformed geometry along the cusps, we 
first recall the usual
hyperbolic geometry of the (rank 1) cusps. We 
parameterize a rank $1$
cusp $\C$ as $S^1 \times [1,\infty) \times \R$,
where the initial $S^1 \times [1,\infty)$ factor is 
a $2$--dimensional cusp $C$.
With the hyperbolic metric, the three 
co--ordinate vector fields are
orthogonal; we denote these by $\frac 
{\partial} {\partial \theta}, \frac {\partial} {\partial z}$ and 
$\frac {\partial} {\partial y}$ respectively, so
that $\theta \in S^1$, $z \in [1,\infty)$ and $y \in \R$. An
orthonormal basis in the hyperbolic metric is 
$z \frac {\partial} {\partial \theta}, z \frac{\partial}
{\partial z}, z \frac{\partial} {\partial y}$. 
Let $h:\R^+ \to \R^+$ be a
monotone increasing function with $h(z) = z$ 
for $z < 1$, and $h(z) = 2$ for $z \ge 3$.
Then let $$h_t(z) = \frac 1 {1-t} h((1-t)z)$$
and define $g_t$ on $\C$ to be the metric with 
orthonormal basis $h_t(z) \frac {\partial} {\partial \theta}, h_t(z) 
\frac {\partial} {\partial z}, h_t(z) \frac {\partial} {\partial y}$. 
Notice that the group of Euclidean
symmetries of the boundary $\partial \C$ extends to an isometry of $\C$ for the
$g_t$ metric, for all $t$. In particular, the surface
$$H_s = S^1 \times [1,\infty) \times s$$
is totally geodesic for the $g_t$ metric, and therefore 
acts as a barrier surface {\em for all $t$}.

Moreover, as $t \to 1$, the $g_t$ metrics converge
to the hyperbolic metrics on compact 
subsets, and in fact for every compact
$K \subset \C$, there is an 
$s>0$ such that the $g_t$ and the
hyperbolic metrics
agree for $t\le s$. Finally, for each $t>0$, the subset
$S^1 \times [3/(1-t),\infty) \times \R \subset \C$ is isometric to a Euclidean
product, for the $g_t$ metric, and therefore the surface
$$F_t = S^1 \times \frac 3 {1-t} \times \R$$
is totally geodesic for the $g_t$ metric, and also 
acts as a barrier surface.

Finally, notice that the $g_t$ metrics 
lift to a family of {\em isometric} metrics on $\H^3$, and by the 
symmetries above, therefore have uniformly pinched sectional 
curvatures, and are uniformly bilipschitz to the hyperbolic metric in 
the region bounded away from the cusps by $F_t$.

Let $N^t_0$ be the neutered space whose boundary consists of the
surfaces of type $F_t$ constructed above.  
Endow $N^t_0$ with the $g_t$ metric.  Now apply \cite {MSY},
as in Lemma \ref{construction_lemma}, to the surface 
$\partial W\cap N^t_0$ to obtain
the surface $\partial W^1_t$ which 
is $g_t$-least area among all surfaces properly
isotopic to $\partial W\cap N^t_0$.  By extending $\partial W^1_t$
\emph{vertically} we obtain the
surface $\partial W_t\subset N$ which is properly isotopic 
to $\partial W$. As in Lemma \ref{construction_lemma} these surfaces 
weakly converge geometrically to a surface $\partial W'$.   We will 
show that there is a proper isotopy of $\partial W$ to
$\partial W'$.

\vskip 12pt
Let $N^\epsilon_0$ denote a fixed neutered space 
transverse to $\partial W'$
and  countably many
$\partial W_t$'s 
which converge to $\partial W'$.   Define $\partial W^\epsilon_t$ to 
be $\partial W_t \cap N^\epsilon_0$ together with the disc components 
of $\partial W_t\backslash N^\epsilon_0$.
Since $\area_{g_t}(\partial W^\epsilon_t)$ is uniformly bounded, and
the hyperbolic area form is dominated on all $2$--planes by $g_t$, the
hyperbolic area of $\partial W^\epsilon_t$ is uniformly bounded.  We
show that a disc $D$ of $\partial W_t\backslash N^\epsilon_0$ cannot stray too far 
into the cusp  and hence for all $t$, $\partial W^\epsilon_t\subset 
N^\eta_0$ for some sufficiently small
$\eta$.  Indeed, the lift $\tilde D$ to the universal cover $\tilde
N$ of $N$ is an embedded disc of uniformly bounded area.   
If $t$ is very close to 1, then $d_\rho(\partial N^\epsilon_0,\partial N^t_0)=d_t>0$.  
If $x\in \tilde D$ and $d_\rho(x,\partial N^\epsilon_0\cup N^t_0)=d_t/2$, 
then $\area_\rho(\tilde D\cap (N^\epsilon_0-N^t_0))>\pi d_t^2/4$. 
Therefore, for $t$ sufficiently large, $\tilde D$ and
hence $D$ has uniformly bounded $\rho$-diameter.

\vskip 12pt

Therefore, if $\epsilon < \eta$, then the $\partial W^\epsilon_t$'s converge
weakly to the surface $N_0^\epsilon\cap W'$ which we define to be $\partial
W'_\epsilon$. For $t$  sufficiently large
$\partial W^\epsilon_t$ and
$\partial W'_\epsilon$ 
are of the same topological type and very close geometrically. 
By Lemma~\ref{controlled spilling},
$\partial W'_\epsilon\cap N^\epsilon_0$ has no components which can be
homotoped rel boundary into
$\partial N^\epsilon_0$, hence $\partial W^\epsilon_t$ shares 
the similar property.  Since $\partial W$ has
no $P$-essential annuli disjoint from $\Delta_1$ it follows that each
component of $\partial W^\epsilon_t\backslash \textint N_0$  can be properly homotoped 
in $\partial W^\epsilon_t$ into an end of that surface.  Therefore, 
if some non disc component of $\partial W_t-\textint N^\epsilon_0$ 
was not a half open annulus, then one can find a
component $E$ of $\partial W_t\cap N^\epsilon_0$ which can be homotoped rel
$\partial E$ into $\partial N^\epsilon_0$, which is a contradiction.  Note 
that $\partial W'_\epsilon$ is of the same topological type as 
$\partial W\cap N_0$ and and that the $\partial W'_\epsilon$'s form 
an exhaustion of $\partial W'$.  By arguing as in the proof of Lemma 
\ref{approximate_by_isotopies},
there exists a homotopy $F:\partial W\times I\to N$ 
with the property that $F(\partial W\times 0)=\partial W$, 
for infinitely many $t<1$,
$F(\partial W\times t)=\partial W_t$ and 
$F(\partial W\times 1)=\partial W'$.

If $\partial W'$ intersects 
$\partial \C$ in the subset $S^1 \times
[s,\infty)$ for some $s$, then for $t$ sufficiently large $\partial
W_t$ must  intersect
$\partial \C$ in the subset $S^1 \times [s-\epsilon,\infty)$.  
Since projection to the barrier surface $H_{s-\epsilon}$ along
horoannuli and horotori is area reducing, this  implies that
$\partial W_t\cap \C$ is contained in
$S^1 \times [1,\infty) \times [s-\epsilon,\infty)$, which in turn
implies that $W'\cap \C\subset S^1 \times [1,\infty) \times [s,\infty)$.  As in
Lemma~\ref{construction_lemma}, the surfaces $\partial W_t$
converge on compact subsets to $\partial W'$. 
The main results of \S 1 imply that $\partial W'$ is $\Delta_1$--minimal.

A similar argument proves similar facts about $S_t$, $T_t$  and $T'$.
\end{proof}

Let $G_i$ denote in$_*(\pi_1(F))\subset \pi_1(W_i)$.  Fix a basepoint
$f\in F$.  Let $ X_i$ denote the 
covering space of $W_i$ (based at $f$) with group $G_i$.
The homotopy of
$\Gamma_i$ into $F$ supported in $W_i$ lifts to $X_i$, hence 
provides us with a canonical
$\hat\Gamma_i$ of closed lifts of 
$\Gamma_i$ in 1-1 correspondence
with $\Gamma_i$.
Since $W_i$ is an 
atoroidal Haken manifold with nonzero Euler characteristic, 
it follows by Thurston that $\inte(X_i)$ is topologically tame
(see Proposition 3.2 \cite{Ca}).  By \cite{Tu2} a
compactification $\bar X_i$ of $\inte(X_i)$ extends 
$\inte(X_i)\cup \partial_h\hat F\cup\partial_p\hat F$, 
where $\hat F$ is the lift of $F$ to $X_i$. 
Since  $\hat F$ is a core of $\bar X_i$ it follows that
$\bar X_i$ is a union of a closed (possibly disconnected or empty)
orientable $\text{surface} \times I$ with
1-handles attached to the 
$\text{surface} \times 1$ side.
Let $\bar S_i$ denote the unique 
boundary component of $\bar X_i$ which
is not a closed component of 
$\hat F$. Push 
$\bar S_i\backslash \inte(\hat F\cap \partial \bar X_i)$ slightly to
obtain  a properly embedded surface $\hat S_i\subset X_i$ with
$\partial \hat S_i=\partial\delta \hat F$ via a 
homotopy disjoint from $\hat\Gamma_i$.  
Being connected with the same Euler characteristic and same 
number of boundary components, $\hat S_i$ is of 
the same topological type as $\delta F$.  Let $\hat Z_i''\subset \bar X_i$ be the 
compact region with frontier
$\hat S_i$.  Let $\chi:=|\chi(\hat S_i)|=|\chi(\delta F)|$.   Define
$\partial_p X_i$ and $\partial_h X_i$ the
respective  preimages of $\partial_pW_i$ and 
$\partial_hW_i$.

Let $W'_i$ denote $W_i$ together with the components of $N_0\backslash\inte(M)$ 
which hit $\partial
W_1$.  Let $Y_i$ be the cover of $W'_i$ with $\pi_1(Y_i)=G_i$.  As in the parabolic free case,
$X_i$ naturally embeds in $Y_i$ and the inclusion is a homotopy equivalence.  Define
$\partial_pY_i$ to be the preimage of $\partial_pW'_i$.  Note that $\delta W_i$, the frontier
of $W_i\subset M$ equals $\delta W'_i$, the frontier of $W'_i$ in $N_0$.  

If possible compress
$\hat Z_i''$ along 
$\delta\hat Z_i''=\hat S_i$
             via  compressions that hit
$\hat\Gamma_i$ at most once.  Continue in this manner to 
obtain the
region  $\hat Z_i'$ whose frontier is 2-incompressible rel 
$\hat B^i : = \hat \Gamma_i \backslash \{\gamma_{i_1},\cdots, \gamma_{i_m}\}$ 
where both $m$ and $|\chi(\delta \hat Z_i')| \le \chi$. 
Since $X_i\backslash \hat\Gamma_i$ is irreducible, we can assume
that no component of $\partial \hat Z_i'$ is a 2-sphere.

Before we shrinkwrap
$\delta W'_i$ and $\delta \hat Z_i'$ we need to \emph{annulate} them, i.e.
compress them
along essential annuli into $P$ and $\partial_pY_i$.
Geometrically we
are eliminating accidental parabolics so that we can invoke 
Lemma \ref{construction_lemma_parabolic}.

Let $L_1,\cdots, L_k$ be a
maximal collection of pairwise disjoint, embedded, essential  annuli in $N_0$
disjoint from
$\Gamma_i$ such that for each $j$,
$\partial L_j$ has one component on $\delta W'_i$ and one component on $P$.
Furthermore assume
that $\inte(L_i)\cap \delta W'_i=\emptyset$.  Now annulate $\delta W'_i$ along
each $L_i$ to obtain
the surface $\delta W^L_i$.  So if $L_i$ lies to the outside of
$W'_i$, then the effect
on $W_i$ is to add $N(L_i)$.  If
$L_i\subset W'_i$ and $L_i\times I$ is a product neighborhood, then
this annulation
deletes $L_i\times \inte(I)$
from $W'_i$.     There are  no
$P$-essential annuli
for $\delta W^L_i$ disjoint from $\Gamma_i$ and
$\partial W^L_i$ is
2-incompressible rel $\Gamma_i$. Indeed, since $\delta W_i^L$ is embedded and 2-incompressible,
we need only consider embedded $P$-essential annuli, by the generalized Loop theorem.  The
modification
$W'_i\to W^L_i$  induces a modification of $Y_i$ as follows.  If $L_j$
annulates $W'_i$
to the outside, then enlarge $Y_i$ in the natural way.  This will enlarge the
parabolic boundary $\partial_pY_i$.  If $\delta W'_i$ gets annulated to the
inside, then do not change $Y_i$. By abuse of notation, we relabel the space
obtained from $Y_i$ as $Y_i$.  Let $W^l_i$ denote $W_i$ modified only along outer
$P$-essential annuli.  Note that $Y_i$ is a covering space of $W^l_i$.  In like manner,
annulate
$\delta
\hat Z'_i$ along a maximal collection of pairwise disjoint annuli which are disjoint from
$\hat B_i$. Let $\hat Z_i\subset Y_i$ denote the result of annulating $\hat Z_i'$.
Note that $\hat Z_i$ can be constructed so that there are no $\partial_p Y_i$-essential
annuli for $\delta \hat Z_i$ disjoint from $\hat B_i$.
\vskip 12pt

Now fix $i$. Let $\gamma\in
\Gamma\backslash W^L_i$.  Apply
Lemma~\ref{construction_lemma_parabolic} using the following dictionary
between our setting and the setting of 
Lemma~\ref{construction_lemma_parabolic}:
$\delta \hat Z_i$ corresponds to the surface $S$, $ W^L_i$
corresponds to $W,\ Y_i$ corresponds to $X$,
$\gamma\cup\Gamma_i$ corresponds to $\Delta_1$,
$\Gamma_i$ corresponds to $\Delta$, and
$\hat B^i$ corresponds to $B$.
We conclude that if
$S^i$ denotes the projection of
$\delta \hat  Z_i$ into $N$, then $S^i$ is homotopic to a $\CAT(-1)$ 
surface
$T^i$ with the following
properties.  The surface
$T^i$ is homotopic 
to a surface $P^i\subset W^{\new}_i$ where $W^{\new}_i$ is isotopic to $W^l_i$ and $P^i$
lifts to 
an embedded surface
$\hat P^i\subset Y_i^{\new}$, where $Y_i^{\new}$ is the corresponding cover of $W^{\new}_i$. 
The isotopy of $W^{\new}_i$ to $W^l_i$ induces a deformation of spaces $Y_i^{\new}$ to
$Y_i$ which fixes
$\hat B^i$ pointwise.   Futhermore,
$\hat P^i$  is isotopic to the corresponding
$\delta \hat Z_i$ via an
isotopy disjoint
from
$\hat B^i$.  Given 
$\epsilon>0$, the $P^i$ can be chosen so that the
homotopy of $T^i$ to $P^i$ restricted 
to $
N^\epsilon_0$ lies in an $\epsilon$-neighborhood of $P^i\cap
N_0^\epsilon$.  By abuse of notation we will view
$\hat P^i$ as bounding the region $\hat Z_i\subset Y^{\new}_i$ and we will drop the
superscripts new, etc.

Let $\{ \alpha_i\}$ be a locally finite collection of embedded proper
rays in $N_0$ to $\mE$
              emanating from $\{\gamma_i\}$.

Let $\pi:Y_i\to N$ be the composition of the covering map to $W^l_i$ and
inclusion.   Let
$B^i=\pi(\hat B^i)$.  If $b\in B^i$ and is disjoint from $N(P^i,1)$, then
some component of $T^i$ homologically separates $b$ from $\mE$.
Indeed if
$\alpha_b$ is the ray from $b$ to $\mE$, then $\piinv(\alpha_b)\cap
\hat Z_i$ is a finite union of compact segments.  If both endpoints
lie in $\partial Y_i$, then it contributes nothing to the algebraic 
intersection number
$\langle \alpha_b,P^i \rangle$.  Otherwise it has one endpoint in 
$\piinv(\alpha_b)
$ and one in
$\partial Y_i$ and hence contributes +1.  Therefore
$$\langle \alpha_b,T^i\rangle = \langle \alpha_b,P^i \rangle >0.$$ 
Since $|B^i|\to \infty$, the $B^i$'s are weakly 1000 separating and the $T^i$'s have
uniformly bounded area, it follows that for $i$ sufficiently large, 
some $b\in B^i\backslash B^{j_i}$ is disjoint from
$N(P^i,1)$, where $j_i<i$ and lim$_{i\to \infty} j_i=\infty$.  Therefore some
subsequence of components of
$\{T^i\}$ exits
$\mE$.

By reducing $\epsilon$, if necessary, we can assume that
$\partial N^\epsilon_0$ is  transverse to all the
$T^i$'s.  By Lemma \ref{controlled spilling}, for each $i$, each
component of $T^i\cap (N\backslash \textint(N_0))$
is either a disc or a half open annulus.
Therefore, the restriction of
each component of $T^i$ to $N_0$ is a connected surface.

\begin{lemma}\label{norm parabolic}  Let $M$ be a relative end-manifold with core $C$ of the
form
$\partial_hM\times I
\cup 1-handles$.  Let $Z$ denote the closure of $M-C$ with $\partial_pZ=\partial_pM\cap Z$ and
$\partial_{\mE}Z=\partial_{\mE}C$.  \begin{enumerate} \item $Z$ is Thurston norm minimizing in
$H_2(Z,\partial_pZ)=\BZ$.
\item If $R$ is a Thurston norm minimizing surface (in either the singular or embedded norms),
representing $[\partial_{\mE}Z]\in H_2(Z,\partial_pZ)$, then for each component $Q$ of
$\partial_pZ$ we have $|R\cap Q|=1$.  In particular, $R$ has no $P$-essential annuli in
$Z$.\end{enumerate}\end{lemma}

\begin{corollary}  Let $N$ be a complete hyperbolic 3-manifold with neutering $N_0$ and
relative core $C$ for $N_0$.  Let 
$R\subset N_0-C$ be a primitive Thurston norm minimizing surface representing an element of
$H_2(N_0-C,\partial_pN_0)$.   Then every homotopy of an accidental parabolic of $R$ into
$\partial_pN_0$ must cross $C$.  \qed\end{corollary}

\noindent\emph{Proof of Lemma.}   The proof of 1 is similar to that of Lemma \ref{norm}.
Recall that since
$C$ is a core, the inclusion
$(\partial_{\mE}Z,\partial\partial_{\mE}Z)\to (Z,\partial_p(Z))$ is an isomorphism. 

Now let $R$ be a possibly singular Thurston norm minimizing surface representing
$[\partial_{\mE}Z]$.  By
\cite{G1}, $\chi(R)=\chi(\partial_{\mE}Z)$, so if $R$ hits $\partial_p(Z)$ in extra components,
then $\genus(R)<\genus(\partial_{\mE}Z)$.  Let $S\subset Z$ be an embedded surface representing
$[\partial_{\mE}Z]$ such that $R$ lies in $\inte(Z')$ where $Z'$ is the compact submanifold cut
off by $S$.  If $\{a_1, \cdots, a_{2g}\}$ is a basis of cycles in
$H_1(\partial_{\mE}Z,\partial
\partial_{\mE}Z)$ which are disjoint from $\partial\partial_{\mE}Z$, then there exists surfaces
$A_1,\cdots, A_{2g}$ with $\partial A_i\subset \partial_{\mE}Z\cup S$ and
$[A_i\cap\partial_{\mE}Z]=n_i[a_i]\in H_1(\partial_{\mE}Z,\partial
\partial_{\mE}Z)$ where $n_i\neq 0$.  For each $i$, let $[A_i\cap S]=b_i \in H_1(S,\partial
S)$.  Since the subgroup of $H_1(S,\partial S)$ which restricts trivially to $\partial S$ is of
rank
$<2g$, it follows that $b_1,\cdots, b_{2g}$ are linearly dependent.  This  implies that
inclusion 
$(\partial_{\mE}Z,\partial\partial_{\mE}Z)\to (Z,\partial_p(Z))$  is not $H_1$-injective, a
contradiction.  If $R$ had a $P$-essential annulus, then we can construct a norm minimizing
surface $R'$ with $|\partial R'|=|\partial R|+2$.\qed

\vskip 12pt
 
We next show that if some component $T$ of $T^i$ has the property that
$T\cap  N_0$ homologically
separates $C$ from $\mE$, then $T=T^i$  is homeomorphic to $\partial_{\mE}C$ and  represents
the class
$[\partial_{\mE}C]\in H_2(N_0,P)$. Suppose that
$[T\cap N_0]=n[\partial_{\mE}C]\in H_2(N_0,P)$.  By Lemma \ref{controlled spilling}, after a
homotopy supported in a small
neighborhood of the cusps we can push the disc components of $T\cap
N\backslash\inte(N_0)$ into
$N_0$ and get 
$|\chi(F)|\ge\chi\ge |\chi(T)|=|\chi(T\cap N_0)|$.  
By Lemma \ref{end-engulfing parabolic}, $|\chi(\partial_{\mE}C |\ge |\chi(F)|$.  
On the other hand 
$$|\chi(T)|\ge x_s(n[\partial_{\mE}C])=
x(n[\partial_{\mE}C])=nx(n[\partial_{\mE}C])=n|\chi(\partial_{\mE}C)|,$$
where the $x$ and $x_s$ respectively denote the Thurston and singular Thurston norms and the
inequality is by definition, the first equality by \cite{G1}, the second by
\cite{T2} and the third by Lemma \ref{norm parabolic}.    The
only possibility is that
$n=1$ and $|\chi(T)|=|\chi(\partial_{\mE}C)|=\chi$ and hence $T=T^i$.  By Lemma \ref{norm
parabolic},
$T$ and $\partial_{\mE}C$ have the same number of boundary components and hence $T=T^i$ is
homeomorphic to $\partial_{\mE}C$.  
 In particular no compressions or annulations occurred to $\hat S^i$.

We claim that the sequence $\{T^i\cap N_0\}$ exits $N_0$.  Otherwise,  there
exists an $m$ with $1\le m\le \chi$, a subsequence $T^{i_1}, T^{i_2},\cdots$
and a compact connected submanifold $K_1 \subset N_0$ such that 
$C\subset K_1$ and for each
$j$,\, $m$ components of
$T^{i_j}$ non trivially intersect $K_1$ and if $R^{i_j}$ are the
components of $T^{i_j}$ which miss $K_1$,
then $R^{i_j}\cap N_0$ is an exiting sequence.  
Since each component
$T$ of $T^{i_j}$ has
$T\cap N_0$ connected, it follows from the bounded diameter lemma that there exists a
compact set $K_2$ such that for all $j$, if $T$ is a component of
$T^{i_j}$ with $T\cap K_1\neq
\emptyset$, then $T\cap N_0\subset K_2$.   Let
$N$ be so large that $\gamma_N\cap\alpha_N\cap N_2(K_2)=\emptyset
$ and
$\gamma_N \subset B_{i_j}$ for infinitely many values of $j$.
Let $\beta_N$ be
a path from $\gamma_N$ to $K_2$.  Since $R^{i_j}$ exits $N_0 $ it follows
that for $j$ sufficiently large
$(\gamma_N\cup \beta_N) \cap N_2(T^{i_j})=\emptyset$.   This implies that
some component
             $T$ of $T^{i_j}$ homologically separates $\gamma_N$ and
hence $C$ from
$\mE$.  Therefore $|\chi(T)|=\chi$.  Since $T\cap \alpha_N\neq 
\emptyset$, this implies that
$T\subset R^{i_j}$ and hence $m=0$, which is a
contradiction.

Since the sequence $\{T^i\cap N_0\}$ exits $N_0$ it follows from the
previous paragraphs that for $i$ sufficiently large, $T^i$ is
homeomorphic to $\hat S^i$, and $T^i\cap N_0$ represents the class
$[\partial_{\mE}C]\in H_2(N_0,P)$.   Since $\{T^i\}$ exits $\mE$, if
$B$ is a cusp of $
N$ parametrized by $S^1\times [0,\infty)\times \BR$, then by Proposition
\ref{construction_lemma_parabolic},  given $n\in \BR,\  T^i\cap
B\subset S^1\times
[0,\infty)\times (n,\infty)$\qed

\begin{remark}\label{parabolic barrier remark}  Since for $i$ sufficiently
large, $T^i$ is of topological type of $\partial_{\mE}C$, it
follows {\it a posteriori} that no
compressions or annulations occurred in the passage from $\bar S^i$
to $\partial \hat Z_i$.
This mirrors the similar
phenomena seen in the proofs of Canary's theorem and Theorem
\ref{main}.\end{remark}

\noindent\textbf{Proof of Theorem \ref{marden parabolic}.} Tameness of
the ends of $N_0$
follows as in the proof of Theorem \ref{marden}.  In
particular if the end $\mE$ of $N_0$ is not geometrically finite,
then by applying  the
proof of  Theorem 2 \cite {So} to $\{T^i\}$ (with the disc components
of $\{T^i\}\cap$
(cusps) pushed into $N_0$) it follows that
$\mE$ is tame.  Alternatively, as in the proof that Criteria (1)-(4) 
implies tameness, we
can use the hyperbolic surface interpolation technique and basic  $3$-manifold
topology to prove that $\mE$ is tame.  Finally tameness of
$N_0$ implies tameness of $N$.
\qed

\vskip 12pt

\noindent\textbf{Proof of Theorem \ref{clean parabolic}.}  It
suffices to prove Theorem
\ref{clean parabolic} for orientable
manifolds which have the 
homotopy type of a relative end-manifold.
It follows 
from
Theorems
\ref{main parabolic} and \ref{marden parabolic} that  a 
parabolic extension
$U_{P}$ of a neighborhood $U$ of $\mE$ is 
topologically of the form
$\textint (T)\times
[0,\infty)$, where $T$ 
is a surface homeomorphic to
$\partial_{\mE}C$ and $C$ is a core of 
$N_0$.   By
Proposition
\ref{construction_lemma_parabolic}, if $(T^i\cap \partial N_0)\subset
\textint (T)\times [t,\infty)$, then
$T^i\backslash \textint N_0\subset \textint(T)\times [t,\infty)$.  Therefore
$\{T^i\}$ exits compact sets in
$\textint (T)\times [0,\infty)$.  Since
for $i$ sufficiently large,
$T^i$ is properly immersed in  $\textint(T)\times[0,\infty)$ and  homologically
separates $\textint(T)\times 0$ from
$\mE$, it follows that the projection
$T_i$ to $\textint(T)\times 0$ is a proper degree-1 map of a surface of finite
type to itself and hence is properly homotopic to a
homeomorphism.
\qed

\end{document}